\documentclass[12pt,a4paper,oneside]{amsart}
\usepackage[normalem]{ulem}
\usepackage{amsthm}
\usepackage{amsmath}
\usepackage{amssymb}
\usepackage{amscd}
\usepackage[utf8]{inputenc}
\usepackage{typearea}
\usepackage{eufrak}
\usepackage{yfonts}
\usepackage{textcomp}
\usepackage{mathrsfs}
\usepackage{hyperref}
\usepackage[draft]{fixme} % erase draft to disable fixme notes
\usepackage{pdfsync}
\usepackage[active]{srcltx}
\usepackage{xypic}
\usepackage{verbatim}
\usepackage{graphicx, array}
\usepackage{longtable}
\usepackage{epstopdf}
\usepackage{xcolor}
\usepackage{cite}
\usepackage{scalerel}

\newcommand{\N}{{\mathbb{N}}}

\newcommand{\Z}{{\mathbb{Z}}}
\newcommand{\Q}{{\mathbb{Q}}}
\newcommand{\C}{{\mathbb{C}}}
\newcommand{\Oo}{{\mathcal{O}}}
\newcommand{\OAB}{{\mathcal{O}_{A,B}}}
\newcommand{\OGE}{{\mathcal{O}_{G,E}}}

\newcommand{\OmA}{{\Omega_{A}}}

\newcommand{\Grpd}{\mathcal{G}}
\newcommand{\Grpdu}{\mathcal{G}^{(0)}}
\newcommand{\Grpdc}{\mathcal{G}^{(2)}}
\newcommand{\CG}{{\mathcal{G}_{\text{tight}}(\SGE)}}
\newcommand{\SGE}{{\mathcal{S}_{G,E}}}
\newcommand{\EGE}{{\mathcal{E}}}

\newcommand{\SE}{{\mathcal{S}_{E}}}
\newcommand{\Eu}{\widehat{\mathcal{E}}_\infty}
\newcommand{\Et}{\widehat{\mathcal{E}}_{\text{tight}}}
\newcommand{\Cat}{\Lambda}
\newcommand{\Obj}{\Lambda^\circ}
\newcommand{\PisoG}{\text{PIso}(\Cat,\dmap)}
\newcommand{\dom}{\operatorname{dom}}
\newcommand{\ran}{\operatorname{ran}}
\newcommand{\Semi}{\mathcal{S}}
\newcommand{\SemiT}{\mathcal{T}}
\newcommand{\Zig}[1][\Lambda]{\mathcal{Z}_{#1}}
\newcommand{\ZigM}[1][\Lambda]{\mathcal{Z}(#1)}
\newcommand{\Idem}[1][\Semi]{\mathcal{E}(#1)}
\newcommand{\osh}[1][]{\sigma_{#1}}
\newcommand{\elmap}[2]{
	\tau^{#1}\osh^{#2}}

\newcommand{\Filt}{\hat{\mathcal{E}}_0}
\newcommand{\FiltT}{\hat{\mathcal{E}}_{tight}}
\newcommand{\UlFilt}{\hat{\mathcal{E}}_{\infty}}
\newcommand{\Mtop}{\mathcal{M}}

\newcommand{\ideG}{\mathbf{1}_G}
\newcommand{\ideQ}{\mathbf{1}_Q}
\newcommand{\dmap}{\mathbf{d}}

\newcommand{\cU}{\mathcal U}
\newcommand{\cP}{\mathcal P}
\newcommand{\cS}{\mathcal S}
\newcommand{\cK}{\mathcal K}

\newcommand{\cI}{\mathcal I}
\newcommand{\cA}{\mathcal A}
\newcommand{\cZ}{\mathcal Z}
\newcommand{\bu}{\bullet}
\newcommand{\lra}{\longrightarrow}
\newcommand{\iso}{\overset{\sim}{\lra}}

\newcommand\n{\mathfrak{n}}
\newcommand{\sour}{\textbf{s}}
\newcommand{\ranT}{\textbf{r}}
\newcommand{\reg}{\mathrm{reg}}

\newcommand{\coker}{\operatorname{coker}}
\newcommand{\id}{\operatorname{id}}
\newcommand\sing{\operatorname{sing}}

\newcommand\triqui{\vartriangleleft}
\newcommand\fX{\mathfrak{X}}

\newtheorem{lemma}{Lemma}[section]
\newtheorem{corollary}[lemma]{Corollary}
\newtheorem{theorem}[lemma]{Theorem}
\newtheorem{proposition}[lemma]{Proposition}

\newtheorem{remark}[lemma]{Remark}
\newtheorem{remas}[lemma]{Remarks}
\newtheorem{definition}[lemma]{Definition}

\newtheorem{example}[lemma]{Example}

\newtheorem{noname}[lemma]{}
\newtheorem{notation}[lemma]{Notation}

\usepackage{tikz}
\tikzstyle{vertex}=[circle]
\tikzstyle{goto}=[->,shorten >=1pt,>=stealth,semithick]
\usetikzlibrary{graphs}
\usepackage{tikz-cd}
\usetikzlibrary{arrows.meta}
\usepackage{wrapfig}
\usetikzlibrary{shapes,decorations.pathreplacing}
\usetikzlibrary{calc}
\usetikzlibrary{backgrounds}
\usetikzlibrary{automata}
\usetikzlibrary{positioning}
\usetikzlibrary{decorations.markings}
\usetikzlibrary{arrows}
\usetikzlibrary{scopes}
\usetikzlibrary{intersections}
\usetikzlibrary{fit}
\tikzset{
  curarrow/.style={
  rounded corners=8pt,
  execute at begin to={every node/.style={fill=red}},
    to path={-- ([xshift=60pt]\tikztostart.center)
    |- (#1) node {}
    -| ([xshift=-60pt]\tikztotarget.center)
    -- (\tikztotarget)}
    }
}
\tikzset{
  labcurarrow1/.style={
  rounded corners=8pt,
  execute at begin to={every node/.style={fill=red}},
    to path={-- ([xshift=60pt]\tikztostart.center)
    |- (#1) node[fill=white] {$\scriptstyle{\id - H_n(\mathcal X_\reg,A)}$}
    -| ([xshift=-60pt]\tikztotarget.center)
    -- (\tikztotarget)}
    }
}
\tikzset{
  labcurarrow2/.style={
  rounded corners=8pt,
  execute at begin to={every node/.style={fill=red}},
    to path={-- ([xshift=60pt]\tikztostart.center)
    |- (#1) node[fill=white] {$\scriptstyle{\id - \Phi_n}$}
    -| ([xshift=-60pt]\tikztotarget.center)
    -- (\tikztotarget)}
    }
}
\tikzset{
  labcurarrow3/.style={
  rounded corners=8pt,
  execute at begin to={every node/.style={fill=red}},
    to path={-- ([xshift=60pt]\tikztostart.center)
    |- (#1) node[fill=white] {$\scriptstyle{\id - \Phi_1}$}
    -| ([xshift=-60pt]\tikztotarget.center)
    -- (\tikztotarget)}
    }
}
\tikzset{
  labcurarrow4/.style={
  rounded corners=8pt,
  execute at begin to={every node/.style={fill=red}},
    to path={-- ([xshift=60pt]\tikztostart.center)
    |- (#1) node[fill=white] {$\scriptstyle{\id - \Phi_0}$}
    -| ([xshift=-60pt]\tikztotarget.center)
    -- (\tikztotarget)}
    }
}

\begin{document}

\title[Self-similar graphs]{An overview on self-similar graphs, their generalizations, and their associated algebras}%

\author{Enrique Pardo}
\address{Departamento de Matem\'aticas, Facultad de Ciencias\\ Universidad de C\'adiz, Campus de
Puerto Real\\ 11510 Puerto Real (C\'adiz)\\ Spain.}
\email{enrique.pardo@uca.es}\urladdr{https://sites.google.com/gm.uca.es/webofenriquepardo/inicio}

\dedicatory{Dedicated to Anna, my light and my strength.}

\thanks{The author was partially supported by PAIDI grant FQM-298 of the Junta de Andaluc\'{\i}a,  by the Spanish State Research Agency (through grant number PID2023-147110NB-I00), and by the grant ``Operator Theory: an interdisciplinary approach'', reference ProyExcel 00780, a project financed in the 2021 call for Grants for Excellence Projects, under a competitive bidding regime, aimed at entities qualified as Agents of the Andalusian Knowledge System, in the scope of the Plan Andaluz de Investigaci\'on, Desarrollo e Innovaci\'on (PAIDI 2020), Consejer\'ia de Universidad, Investigaci\'on e Innovaci\'on of the Junta de Andaluc\'ia.}

\subjclass[2020]{46L05, 46L55, 46L80, 16S99, 16S10, 16S88, 16S36, 18B40, 20F65, 20J05, 22A22, 20L05, 20M18, 20M25}

\keywords{Self-similar group, Nekrashevych algebras, Katsura algebra, Self-similar graph, Inverse semigroup, Tight groupoid, groupoid $C^*$-algebra, Steinberg algebra, $k$-graph, Twisted groupoid, Left cancellative small category, Zappa-Sz\'ep product, Groupoid Homology, K-theory}

%\date{\today}

\begin{abstract}
In these notes, we introduce the concept of self-similar graph, associated with groups acting on graphs. We define the corresponding $C^*$-algebra using different complementary approaches, to understand its basic properties. We also analyze various generalizations that appear in the literature and, in particular, review the relationship of this construction with Zappa-Sz\'ep products. Finally, we present very recent results on homology and $K$-theory for these algebras.
\end{abstract}

\maketitle

\tableofcontents

\section* {Introduction}

These notes were originally intended to be complementary material for an introductory course on self-similar graphs and their algebras, presented by the author at the CIMPA School ``K-theory and Operator Algebras'', held in La Plata and Buenos Aires (Argentina) from 28 July to 8 August 2025.\vspace{.2truecm} 

Self-similar graphs and their algebras were introduced in 2017 by Ruy Exel and the author \cite{EP2} as a generalization of two previous constructions: self-similar groups \cite{N2} and Katsura algebras \cite{Kat2}. Although it is a recent topic, in the last years several authors have studied various generalizations, extending the scope of these objects to a huge variety of algebras. However, so far no survey has given a general overview of this topic and its results.

Given the extent of the existing material and since we did not pretend to write an exhaustive monograph, we had to make some choices to keep the notes manageable enough. Thus, we decided:\begin{enumerate}
\item Do not include proofs of the results. Instead, we included precise references of all the results stated, to facilitate this information to the interested reader.
\item To avoid some topics and so reduce the size of the notes. For example, we do not talk about self-similar actions of compact groups on topological graphs \cite{BKQ}, bicategorical descriptions of self-similar groups \cite{AKM22}, $kk$-theoretical results \cite{Willie}, among others.
\item In each section, the original notation used by the authors of the corresponding reference is kept, even if it does not fit the standard introduced in \cite{EP2}, to facilitate the work for potential readers who browse the original sources.
\end{enumerate}

Hence, we concentrate ourselves in studying the original concept from an inverse semigroups/groupoid approach, and to look at its natural generalizations from the most elementary situation to the most abstract one: from actions of groups on finite graphs to Zappa-Sz\'ep products of groupoids on left cancellative small categories.

We hope these notes will help the reader to have a global vision of the concept.

\section{The origins}

In this section, we revisit two previous constructions, which were source of inspiration for self-similar graph definition.

\subsection{Self-similar groups and Nekrashevych algebras}\label{Subsect:Nekrashevych}

Self-similar groups appeared as a class of groups whose presentation is easily encoded by automata, providing examples of exotic behaviors (nontrivial torsion, intermediate growth, etc.). The simplest way of describing these groups is through their action on $n$-regular rooted trees (equivalently, on words written on an alphabet of $n$ symbols). we follow \cite{GNSV, N-book} for definitions and results.

To be precise: given a (finite) set $X$, regarded as an alphabet, we can define $X^*$ to be the set of finite words, and $X^\omega$ to be the set of infinite words on this alphabet. The set $X^\omega$ becomes an ultrametric space when equipped with the metric $d(u,v)=2^{-\vert u\wedge v \vert}$, where $u\wedge v$ denotes the longest common prefix of $u$ and $v$.

Since every closed ball of $X^\omega$ equals the set $uX^\omega$ for some given prefix $u\in X^*$, every isometry $f$ of $X^\omega$ induces a permutation of $X^*$. As this permutation preserves inclusion of balls, it is an automorphism of the Cayley graph of $X^*$ that preserves lengths (and so levels). This means that, if $X$ is a finite set, and $G$ is a discrete group acting on $X^\omega$ by isometries, then it induces an action of $G$ on $X^*$ by $\vert X\vert$-ary tree automorphisms.

Now, let $g:X^*\rightarrow X^*$ be an automorphism of the rooted tree $X^*$. Given a vertex $v\in X^*$, we consider the subtrees $vX^*$ and $g(v)X^*$ of $X^*$, which are connected through the restriction of $g: vX^*\rightarrow g(v)X^*$. Since both $vX^*$ and $g(v)X^*$ are naturally isomorphic to the whole tree $X^*$, we can represent the above restriction map as an automorphism $g{\vert_v}: X^*\rightarrow X^*$, which is uniquely determined by the rule $g(vw)=g(v){g\vert_ v}(w)$.

\begin{definition}\label{Def: restriction map}
The element $g{\vert_ v}\in \text{Aut}(X^*)$ is the restriction of $g$ in $v$. It satisfies two obvious properties:
\begin{enumerate}
\item $g\vert_{v_1v_2}=(g{\vert_ {v_1}}){\vert_ {v_2}}$ for all $v_1,v_2\in X^*$.
\item $(g_1g_2){\vert_ v}={g_1}{\vert_ {g_2(v)}}{g_2}{\vert_ v}$ for all $g_1, g_2\in G$ and $v\in X^*$.
\end{enumerate}
This means that the restriction map is a $1$-cocycle ${-}_{\vert -}: G\times X^* \rightarrow G$.
\end{definition}

Now, we are ready to introduce the main definition:

\begin{definition}\label{Def: self-similar group}
Given a finite set $X$, a \underline{faithful} action of a (discrete) group $G$ on $X^\omega$ by isometries is self-similar if for all $g\in G$ and all $x\in X$ there exists $g{\vert_x}\in G$ such that, for all $w \in X^\omega$, we have $g(xw)=g(x)g{\vert_x}(w)$. Notice that this is equivalent to the definition of a self-similar action of $G$ on $X^*$ by (level-preserving) automorphisms; since the action is faithful, the element $g{\vert_x}$ is unique. Moreover, the restriction rules described in Definition \ref{Def: restriction map} allow us to inductively extend the definition to all $v, w\in X^*$. 
\end{definition}

Definition \ref{Def: self-similar group} means that $(G,X)$ is self-similar if $G$ is a subgroup of $\text{Aut}(X^*)$ invariant by the restriction $1$-cocycle. We can formalize this fact by using the permutational wreath product construction: $G\leq \text{Aut}(X^*)$ is self-similar if $G\leq S(X)\wr G$ (for more details, see e.g. \cite[Sections 1.4-1.5]{N-book}).

Let us illustrate this definition with a classical example:

\begin{example}\label{Exam:Grigorchuk group}
The Grigorchuk group is a group $G$ acting on $X^*$ for the alphabet $X=\{ 0, 1\}$, generated by four automorphisms $a,b,c,d\in \text{Aut}(X^*)$, encoded by the automaton %(\ref{f:nucleus.grigorchuk}):
%\begin{figure}[htbp]
%\begin{figure}[t]
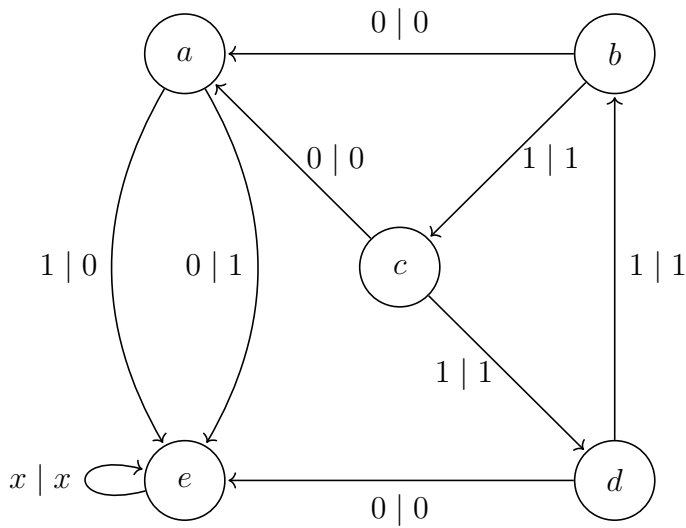
\begin{center}
\begin{tikzpicture}[->,shorten >=1pt,%
auto,node distance=4cm,semithick,
inner sep=5pt,bend angle=30]
\tikzset{every loop/.style={min distance=10mm,looseness=10}}
%\tikzset{every state/.style={minimum size=30pt}}.
\node[state] (C) {$c$};
\node[state] (A) [above left of=C] {$a$};
\node[state] (B) [above right of=C] {$b$};
\node[state] (D) [below right of=C]{$d$};
\node[state] (E) [below left of=C]{$e$};
%\tikzstyle{every node}=[font=\footnotesize]
\path  (E) edge [loop left] node [left] {$x\mid x$} (E)
        (A) edge [bend right] node [left] {$1\mid 0$} (E)
        (A) edge [bend left] node [left] {$0\mid 1$} (E)
        (B) edge [right] node [above] {$0\mid 0$} (A)
        (B) edge [right] node [right] {$1\mid 1$} (C)
        (C) edge [right]  node [right] {$0\mid 0$} (A)
        (C) edge [right]  node [left] {$1\mid 1$} (D)
        (D) edge [right] node [right] {$1\mid 1$} (B)
        (D) edge [left] node [below] {$0\mid 0$} (E);
\end{tikzpicture}
\end{center}
The group also can be defined recursively by the corresponding action and restriction of these generators on $X$:
\begin{enumerate}
\item $a(0)=1$, $a{\vert_0}=e$ \& $a(1)=0$, $a{\vert_1}=e$.
\item  $b(0)=0$, $b{\vert_0}=a$ \& $b(1)=1$, $b{\vert_1}=c$.
\item  $c(0)=0$, $c{\vert_0}=a$ \& $c(1)=1$, $c{\vert_1}=d$.
\item  $d(0)=0$, $d{\vert_0}=e$ \& $d(1)=1$, $d{\vert_1}=b$.
\end{enumerate}
%\caption{State diagram for the automaton of the Grigorchuk group~\label{f:nucleus.grigorchuk}}
%\end{figure}
%Finally, the permutational wreath product description of $G$ corresponds to take $S(X)=\{\id, \sigma\}$, and then state the wreath recursion:
%\begin{enumerate}
%\item $a=\sigma$.
%\item $b=(a,c)$.
%\item $c=(a,d)$.
%\item $d=(e, b)$.
%\end{enumerate}
\end{example}

Nekrashevych \cite{N1, N2} defined a $C^*$-algebra associated to a self-similar group (also known as Nekrashevych algebra), as follows:

\begin{definition}\label{Def: Nekrashevych algebras}
Let $(G,X)$ be a self-similar group. We define $\C (G,X)$ to be the universal complex $\ast$-algebra generated by
\[\{u_g \mid g\in G\}\cup \{s_x \mid x\in X\} \]
satisfying the following relations:
\begin{enumerate}
\item $s_x^*s_y=\delta_{x,y}$ for all $x,y\in X$.
\item $\sum\limits_{x\in X}s_xs_x^*=1$.
\item The map $g\mapsto u_g$ is a unitary representation of $G$.
\item $u_gs_x=s_{g(x)}u_{g_{\vert x}}$ for all $g\in G$ and all $x\in X$.
\end{enumerate}
We define the $C^*$-algebra $\mathcal{O}(G,X)$ to be the universal $C^*$-algebra given by the same set of generators and relations; it turns out to be the norm-completion of $\C (G,X)$ in a suitable norm.
\end{definition}
Given any field $K$, we can also define $K (G,X)$ to be the universal $K$ algebra with involution $\ast$ satisfying the same relations as above (see e.g. \cite{HPSS} for a general definition when the field $K$ is replaced by a unital commutative ring $R$). For this $K$-algebra, we use the notation $\Oo_{G,X}^K$.

When $G=\{1\}$, $\C (G,X)$ is the complex Leavitt algebra $L_{\vert X\vert}$, and $\mathcal{O}(G,X)$ is the Cuntz algebra $\mathcal{O}_{\vert X\vert}$ (see e.g. \cite{Raeburn}); in fact, it is easily seen that these graph algebras are subalgebras of $\C (G,X)$ and $\mathcal{O}_{\vert X\vert}$ (respectively) for any self-similar group $(G,X)$. 

Nekrashevych gave presentations of $\mathcal{O}_{\vert X\vert}$ both as Cuntz-Pimsner algebras \cite{N1} (which allowed stating the nuclearity of the algebra) and as groupoid $C^*$-algebras associated with a certain ample, second countable groupoid $\mathcal{G}_{(G,X)}$ \cite{N2}. Since the action of $G$ on $X^*$ is faithful, the groupoid $\mathcal{G}_{(G,X)}$ is effective, while the fact that the underlying graph of $(G,X)$ has a unique vertex implies that $\mathcal{G}_{(G,X)}$ is minimal. As we see later on, this forces the algebra $\mathcal{O}(G,X)$ to be simple in reasonable situations (e.g., when $\mathcal{G}_{(G,X)}$ is Hausdorff \cite{BCFS}), but not in general (see e.g. \cite[Section 5]{CEPSS}).

\subsection{Katsura algebras}\label{Subsect:Katsura}

Katsura \cite{Kat1} constructed a family of $C^*$-algebras of combinatorial nature that includes, up to isomorphism, all Kirchberg algebras in the UCT. We recall the definition and basic properties of Katsura algebras that will be needed in the sequel, borrowed from \cite{Kat1}.

\begin{definition}\label{Def:KatAlgData}
Let $N\in \N\cup \{\infty\}$, let $A\in M_N(\Z^+ )$ and $B\in M_N(\Z )$ be row-finite matrices. Define a set $\Omega_A $ by 
\[\Omega_A
:=\{ (i,j)\in \{ 1, 2, \dots ,N\}\times \{ 1, 2, \dots ,N\} \mid A_{i,j}\geq 1 \}.\] 
For each $i\in \{ 1, 2, \dots ,N\}$, define a set $\Omega_A (i)\subset \{ 1, 2, \dots ,N\}$ by 
\[\Omega_A (i):=\{ j \in \{ 1, 2, \dots ,N\}\mid (i,j)\in \OmA \}.\]
Notice that, by definition, $\Omega_A (i)$ is finite for all $i$.  Finally, fix the following relation:

\quad (0) $\Omega_A (i)\ne \emptyset $ for all $i$, and $B_{i,j}=0$ for $(i,j)\not \in \OmA $.
\end{definition}

With these data we can define Katsura algebras:

\begin{definition}\label {Def:KatAlgAlgebra}
Define $\OAB $ to be the universal $C^*$-algebra generated by
mutually orthogonal projections $\{ q_i\}_{i=1}^N$, partial unitaries $\{ u_i\}_{i=1}^N$ with $u_iu_i^*=u_i^*u_i=q_i$,
and partial isometries $\{ s_{i,j,n}\}_{(i,j)\in \OmA , n\in \Z }$ satisfying the relations:
\begin{enumerate}
\item $s_{i,j,n}^*s_{i,j,n}=q_j$ for all $(i,j)\in \OmA $ and $n\in \Z $.
\item $q_i=\sum \limits _{j\in \OmA (i)}\sum \limits _{n=1}^{A_{i,j}}s_{i,j,n}s_{i,j,n}^*$ for all $i$.
\item $s_{i,j,n}u_j=s_{i,j, n+A_{i,j}}$ and $u_is_{i,j,n}=s_{i,j, n+B_{i,j}}$ for all $(i,j)\in \OmA $ and $n\in \Z$.
\end{enumerate}
\end{definition}
Given any field $K$, we can also consider $\Oo_{A,B}^K$ to be the universal K-algebra with the same presentation as above (see also \cite{HPSS}).

\begin{remark}\label{Rem:KatsuraConditions}
The following facts holds: 
\begin{enumerate}
\item The $C^*$-algebra $\OAB $ is separable, nuclear, in the UCT class \cite [Proposition 2.9]{Kat1}. 
\item If the matrices $A,B$ satisfy the following
additional properties:
\begin{enumerate}
\item $A$ is irreducible, and
\item $\vert A_{ii}\vert\geq 2$ and $B_{i,i}=1$ for every $1\leq i\leq N$,
\end{enumerate}
then the $C^*$-algebra $\OAB $ is purely infinite simple, and hence a Kirchberg algebra \cite [Proposition 2.10]{Kat1}.
\item The $K$-groups of $\OAB $ are \cite [Proposition 2.6]{Kat1}:
\begin{enumerate}
\item  $K_0(\OAB )\cong \mbox {coker}(I-A)\oplus \mbox {ker}(I-B)$, and
\item $K_1(\OAB )\cong \mbox {coker}(I-B)\oplus \mbox {ker}(I-A)$.
\end{enumerate}
\item Every Kirchberg algebra can be represented, up to isomorphism, by an algebra $\OAB $ for matrices $A,B$ satisfying the conditions in \ref{Rem:KatsuraConditions}(2) \cite [Proposition 4.5]{Kat2}.
\end{enumerate}
\end{remark}

At first sight, it could seem that there is no relation between the definition of Ne\-kra\-she\-vych algebras and that of Katsura algebras. Let us show that, in fact, we can change Katsura's picture in such a way that both definitions almost match, at least when $N< \infty$.

As in \cite{Kat1}, let us assume that we have two $N\times N$ matrices $A$ and $B$ with integer entries, and such that $A_{i,j}\geq 0$, for all $i$ and $j$.  We may then consider the graph $E$ with vertex set
\[E^0 = \{1,2,\ldots ,N\},\]
and such that, for each pair of vertices $i,j \in E^0$, the set of edges from vertex $j$ to vertex $i$ is a set with $A_{i,j}$ elements, say
\[\{e_{i,j,n} : 0\leq n<A_{i,j}\}.\]
Moreover, assuming that $A$ has no identically zero rows, it is easy to see that $E$ has no sources.

Define an action $\sigma $ of $\Z$ on $E$, which is trivial on $E^0$, and which acts on edges as follows: given $m\in \Z$, and $e_{i,j,n} \in E^1$, let $(\hat{k},\hat{n})$ be the unique pair of integers such that
\[mB_{i,j}+n=\hat{k} A_{i,j} + \hat{n} \text{ and } 0\leq \hat{n}<A_{i,j}.\]
That is, $\hat{k}$ is the quotient and $\hat{n}$ is the remainder of the Euclidean division of $mB_{i,j}+n$ by $A_{i,j}$.  We then put
\[\sigma_m(e_{i,j,n})=e_{i,j,\hat{n}}.\]
In other words, $\sigma_m$ corresponds to the addition of $mB_{i,j}$ to the variable ``$n$'' of ``$e_{i,j,n}$'', taken modulo $A_{i,j}$. In turn, it is easy to check that the map $\varphi: \Z \times E^1 \rightarrow \Z$ defined by
\[\varphi (m, e_{i,j,n})= \hat{k}\]
is a $1$-cocycle in the sense of Definition \ref{Def: restriction map}.

Observe that if $A_{i,j}=0$, then there are no edges from $j$ to $i$, so the value $B_{i,j}$ is entirely irrelevant for the above construction. Therefore, it makes no difference to assume that
\[A_{i,j}=0 \Rightarrow B_{i,j}=0.\]

So, Katsura algebra $\OAB$ seems to be isomorphic to a sort of Nekrashevych algebra $\Oo_{(\Z , E)}$, defined via an action of $\Z$ on $E$ by the automorphism $\sigma$, ``twisted'' by the 1-cocycle $\varphi$, just by sending each $u_m$ to the $m ^{\text{th}}$ power of the unitary
\[u:=\sum _{i=1}^Nu_i\]
in $\OAB$, and sending $s_{e_{i,j,n}}$ to $s_{i,j,n}$; the formalization of this idea is the basis of the definition of a self-similar graph.

Notice that, when $N=1$, the graph $E$ for Katsura's algebras is the same as the one we used above in the description of Nekrashevych's example. However, the former is not a special case of the latter because, contrary to what is required in \cite {N2}, the group action might not be faithful.

\section{Self-similar graphs and their algebras}

Looking at the previous constructions, we can extract common elements to both families of algebras: a (finite) directed graph $E$, a discrete group $G$ acting on $E^0\sqcup E^1$ by graph automorphisms, and a length-preserving $1$-cocycle $\varphi: G\times (E^0\sqcup E^1)\rightarrow G$ allowing to define a sort of self-similar identity: for all $g\in G$ and all $x\in E^0\sqcup E^1$ we have
\[ g\cdot x:=g(x)\varphi(g, x).\]
Using this identity, we could then define a $C^*$-algebra with generators corresponding to a graph $C^*$-algebra $C^*(E)$ extended by a self-similar action of a (non necessarily faithful) unitary copy $u(G)$ of $G$ in such a way that, for all $g\in G$ and all $x\in E^0\sqcup E^1$ we have
\[u_g s_x=s_{g(x)} u_{\varphi(g,x)}.\]

Building a coherent algebra will require to take care of some technical details. Let us formalize the idea. We follow \cite{EP2} throughout this section.\vspace{.2truecm}

\subsection{Self-similar graphs}\label{SubSect:SSGclassic}
Let $E = (E^0, E^1, r , d )$ be a directed graph, where $E^0$ denotes the set of \emph{vertices}, $E^1$ is the set of \emph{edges}, $r$ is the \emph{range} map, and $d$ is the \emph{source}, or \emph{domain} map.

By definition, a \emph{source} in $E$ is a vertex $x\in E^0 $, for which $r^{-1}(x)=\emptyset $.  Thus, when we say that a graph has \emph{no sources}, we mean that $r^{-1}(x)\ne \emptyset $, for all $x \in E^0$.

By an \emph{automorphism} of $E$ we shall mean a bijective map
\[\sigma : E^0 \sqcup E^1 \to E^0 \sqcup E^1\]
such that $\sigma (E^i)\subseteq E^i$, for $i = 0,1$, and moreover such that $r \circ \sigma = \sigma \circ r$, and $d \circ \sigma = \sigma \circ d$, on $E^1$. It is evident that the collection of all automorphisms of $E$ forms a group under composition.

By an action of a group $G$ on a graph $E$ we mean a homomorphism of the group from $G$ to the group of all automorphisms of $E$.

If $X$ is any set, and if $\sigma $ is an action of a group $G$ on $X$, we shall say that a map
\[\varphi :G\times X \to G\]
is a \emph{one-cocycle} for $\sigma$, when
\[\varphi (gh, x) = \varphi \big(g,\sigma _h(x)\big)\varphi (h,x),\]
for all $g,h \in G$, and all $x \in X$.  In particular, $\varphi (1,x) = 1$ for every $x$.\vspace{.2truecm}

The \textbf{Standing Hypothesis} in the classical construction we define is that $G$ will be a discrete group, $E$ will be a finite graph with no sources, $\sigma$ will be an action of $G$ on $E$, $\varphi:G\times E\rightarrow G$ will be a $1$-cocycle for the restriction of $\sigma$ to $E^1$, which must satisfy
\begin{equation}\label{equa:StanHypSSG}
\sigma_{\varphi(g,e)}(x)=\sigma_g(x) \hspace{.2truecm} \forall g\in G, \hspace{.2truecm} \forall e\in E^1, \hspace{.2truecm} \forall x\in E^0.
\end{equation}

As we shall see later, the finiteness of $E$ could be removed from the hypotheses. Also, as remarked in \cite{BKQ}, $(\ref{equa:StanHypSSG})$ can be replaced by the weaker requirement
\begin{equation}\label{equa:StanHypSSGweak}
\sigma_{\varphi(g,e)}(s(e))=\sigma_g(s(e)) \hspace{.2truecm} \forall g\in G, \hspace{.2truecm} \forall e\in E^1.
\end{equation}

By \emph{path} in $E$ of \emph{length} $n\geq 1$ we mean any finite sequence of the form
\[\alpha = \alpha _1\alpha _2\ldots \alpha _n,\]
where $\alpha _i \in E^1$, and $d (\alpha _i) = r (\alpha _{i+1})$, for all $i$ (this is the usual convention when treating graphs from a categorical point of view, in which functions compose from right to left). The \emph{range} of $\alpha $ is defined by $r (\alpha ) = r (\alpha _1)$, while the \emph{source} of $\alpha $ is defined by $d (\alpha ) = d (\alpha _n)$. A vertex $x \in E^0$ will be considered a path of length zero, in which case we set $r (x ) = d (x) = x $.

For every integer $n\geq 0$ we denote by $E^n$ the set of all paths in $E$ of length $n$. Finally, we denote by $E^*$ the sets of all finite paths, and by $E^{\leq n}$ the set of all paths of length at most $n$, namely
\[ E^* = \bigcup\limits_{k\geq 0}E^k \text{ and } E^{\leq n} =  \bigcup\limits _{k=0}^nE^k.\]
It is a straightforward computation to extend $\sigma$ and $\varphi$ to finite paths. From now on, we adopt the shorthand notation
\[g\alpha = \sigma_g(\alpha ).\]

Then, we can write

\begin{proposition}[{c.f. \cite[Proposition 2.4 \& Equation 2.5]{EP2}}]\label{Prop:Equacoes}
For every $g, h\in G$, for every $x \in E^0$, and for every $\alpha $ and $\beta $ in $E^*$ such that $d (\alpha )=r (\beta )$, one has that
\begin{enumerate}
\item $(gh)\alpha = g(h\alpha )$,
\item $\varphi (gh, \alpha ) = \varphi (g,h\alpha )\varphi (h,\alpha ),$
\item $\varphi (g,x ) = g$,
\item $r (g\alpha )=g r (\alpha )$,
\item $d (g\alpha )=g d (\alpha )$,
\item $\varphi (g,\alpha )x =g x $,
\item $g(\alpha \beta ) = (g\alpha ) \varphi (g,\alpha )\beta $,
\item $\varphi (g, \alpha \beta )=\varphi (\varphi (g,\alpha ),\beta )$.
\end{enumerate}
\end{proposition}

Once we have the correct definition for self-similar graphs, we are ready to define a $C^*$-algebra associated with it.

\subsection{Algebras of self-similar graphs}
Now, we fix a graph $E$, an action of a group $G$ on $E$, and a $1$-cocycle $\varphi $ under (\ref{equa:StanHypSSG}). We proceed to build a $C^*$-algebra from these data.

\begin{definition}\label{Def:DefineOGE}
We define $\OGE $ to be the universal unital $C^*$-algebra generated by a set
\[\{p_x \mid x\in E^0\}\cup \{s_e \mid  e \in E^1\} \cup \{u_g \mid g \in G\},\]
subject to the following relations:
\begin{enumerate}
\item $\{p_x \mid x\in E^0\}\cup \{s_e \mid  e \in E^1\}$ is a Cuntz-Krieger $E$-family \cite{Raeburn},
\item the map $u:G\rightarrow \OGE $, defined by the rule $g\mapsto u_g$, is a unitary representation of $G$,
\item $u_gs_e=s _{ge }u_{\varphi (g,e )}$, for every $g \in G$, and $e \in E^1$,
\item $u_gp_x =p_{gx}u_g$, for every $g \in G$, and $x \in E^0$.
\end{enumerate}
\end{definition}

\noindent Given any field $K$, we can also define $\Oo_{G,E}^K$ to be the universal $K$ algebra with the same presentation as above (see \cite{HPSS} for a general version).\vspace{.2truecm}

This construction generalizes our inspiration examples, as well as some other well-known constructions, as follows:

\begin{example}\label{Exam:Nekrashevych}
Let $(G,X)$ be a self-similar group as in \cite[Definition 2.1]{N1}. We may then consider a graph $E$ having only one vertex and such that $E^1=X$.  If we define
\[\varphi (g,x )=g{\vert_ x} ,\]
then the triple $(G, E, \varphi)$ satisfies (\ref{equa:StanHypSSG}), and it is easy to show that $\OGE $ is isomorphic to the Nekrashevych algebra $\mathcal{O}_{(G,X)}$.
\end{example}

\begin{example}\label{Exam:Katsura}
Let us assume that we are given two $N\times N$ matrices $A$ and $B$ with integer entries, and such that $A_{i,j}\geq 0$, for all $i$ and $j$. Then, taking the graph $E$ whose incidence matrix is $A$, and defining the action of $\Z$ on $E$ and the 1-cocycle as at the end of Subsection \ref{Subsect:Katsura}, we can easily show that the Katsura algebra $\OAB$ is isomorphic to the self-similar graph $C^*$-algebra $\Oo_{(\Z , E)}$, just by sending each $u_m$ to the $m ^{\text{th}}$ power of the unitary
\[u:=\sum _{i=1}^Nu_i\]
in $\OAB$, and sending $s_{e_{i,j,n}}$ to $s_{i,j,n}$
\end{example}

\begin{example}\label{Exam:CrossedProduct}
Given any finite graph $E$, and any action $\sigma $ of a group $G$ on $E$, the map $\varphi : G\times E^1 \rightarrow G$ defined by
\[\varphi (g,a)=g, \text{ for } g\in G, a\in E^1\]
is a $1$-cocycle, and the triple $(G,E, \varphi )$ satisfies (\ref{equa:StanHypSSG}). By Definition \ref{Def:DefineOGE}(c), we have
\[u_gs_au_g^*=s_{ga}\]
for any $g$ in $G$, and every $a$ in $ E^1$. Therefore, it is easy to see that
$\OGE $ is isomorphic to the crossed product of the graph $C^*$-algebra $C^*(E)$ \cite {Raeburn} by $G$, relative to the natural action of $G$ on $C^*(E)$ induced by $\sigma $. In particular, if $\sigma $ is the trivial action, we have that $\OGE $ is the maximal tensor product of $C^*(E)$ by the full
group $C^*$-algebra of $G$.
\end{example}

\begin{example}\label{Exam:GraphAlgebra}
Given any finite graph without sources, and any action $\sigma $ of a group $G$ on $E$ fixing the vertices, consider the map $\varphi : G\times E^1 \rightarrow G$ defined by
\[\varphi (g,a)=1 \text{ for } g\in G, a\in E^1.\]
It is easy to see that $\varphi $ is a $1$-cocycle, and that the triple $(G,E, \varphi )$ satisfies (\ref{equa:StanHypSSG}). Since $E$ has no sources, we have,
for any $g$ in $ G$, that
\[u_g=\sum\limits_{x\in E^0}u_gp_x=\sum\limits_{x\in E^0}\sum \limits _{a\in r^{-1} (x)}u_gs_as_a^*=
\sum\limits_{x\in E^0}\sum\limits_{a\in r^{-1} (x)}s_{ga}s_a^*,\]
that therefore lies in the copy of $C^*(E)$ within $\OGE $.  Since the natural representation of
$C^*(E)$ in $\OGE $ is faithful, the conclusion is that $\OGE \cong C^*(E)$.
\end{example}

In order to have a useful picture of $\OGE$ that allows us to study its structure and properties, we need to fix a couple of basic facts.

\begin{lemma}[{\cite[Lemma 3.8]{EP2}}]\label{Lem:RuleforFinitePaths}
Given $\alpha \in E^*$, and $g \in G$, one has that
\[u_g s_\alpha =s_{g\alpha }u_{\varphi (g,\alpha )}.\]
\end{lemma}

\begin{proposition}[{\cite[Proposition 3.9]{EP2}}]\label{Prop:InvSemPicture}
Let
\[\mathcal{S} = \left\{s_\alpha u_g s _\beta ^* \mid \alpha ,\beta \in E^*, g \in G, d (\alpha )=g d (\beta )\right \} \cup \{0\}.\]
Then $\mathcal{S}$ is closed under multiplication and adjoints, and its closed linear span coincides with $\OGE $.
\end{proposition}

In fact, $\mathcal{S}$ is an inverse subsemigroup of $\OGE$ that generates the entire algebra. We will use an abstract version of $\mathcal{S}$ to produce a suitable ample groupoid $\mathcal{G}$ such that $\OGE$ turns out to be $\ast$-isomorphic to the full groupoid $C^*$-algebra $C^*(\mathcal{G})$. Once this is done, we will be able to use the results of \cite{EP2} to characterize key properties of the groupoid used to give a description of simplicity and pure infiniteness of $\OGE$.

\subsection{The inverse semigroup $\SGE$}\label{SubSect:InvSemigrp}

We are working with a self-similar graph $(G,E, \varphi)$ under (\ref{equa:StanHypSSG}). Then

\begin{definition}\label{Def:inverseSemigroup}
Over the set
\[\SGE =\left\{ (\alpha ,g,\beta ) \in E^*\times G\times E^*\mid d (\alpha )=gd (\beta )\right\}\cup \{0\},\]
consider a binary \emph{multiplication} operation defined by
\[(\alpha ,g,\beta ) (\gamma ,h,\delta ) = \left \{
\begin{array}{cc}
(\alpha g \varepsilon ,\hfil \varphi (g,\varepsilon ) h,\hfil \delta ), & \text {if } \gamma =\beta \varepsilon ,  \\
(\alpha ,\ g\varphi (h^{-1} ,\varepsilon )^{-1} ,\ \delta h^{-1} \varepsilon ), & \text {if } \beta =\gamma \varepsilon ,\\\
  0, & \text {otherwise,}
\end{array}
\right.\]
and a unary \emph{adjoint} operation defined by
\[(\alpha ,g,\beta )^*:= (\beta ,g^{-1} , \alpha ). \]
Furthermore, the subset of $\SGE $ formed by all elements $(\alpha ,g,\beta )$, with $g=1$, will be denoted by $\SE$.
\end{definition}

It is easy to see that $\SE $ is closed under the above operations, and that it is isomorphic to the inverse semigroup generated by the canonical partial isometries in the graph $C^*$-algebra $C^*(E)$. Then, we have

\begin{proposition}[{\cite[Proposition 4.3]{EP2}}]\label{Prop:SGE-InvSemgrp}
$\SGE $ is an inverse semigroup with zero.
\end{proposition}

From Proposition \ref{Prop:SGE-InvSemgrp} it is easy to state that the idempotent semilattice of $\SGE $, denoted by $\EGE $, is given by
 \[\EGE = \left\{(\alpha ,1,\alpha )\mid \alpha \in E^*\right\} \cup \{0\}.\]
Evidently $\EGE $ is also the idempotent semi-lattice of $\SE $. For simplicity, we adopt the short-hand notation
\[f_\alpha = (\alpha ,1,\alpha ) \text{ for } \alpha \in E^*.\]

The following is a standard fact in the theory of graph $C^*$-algebras.

\begin{proposition}[{\cite[Proposition 4.6]{EP2}}]\label{Prop:ProdIdemp}
If $\alpha ,\beta \in E^*$, then
\[
f_\alpha f_\beta = \left \{
\begin{array}{cl}
f_\alpha, & \text{if there exists } \gamma \text{ such that } \alpha =\beta \gamma , \\
f_\beta, & \text {if there exists } \gamma \text{ such that } \alpha \gamma =\beta , \\
0, & \text {otherwise.}
\end{array}
\right .\]
\end{proposition}

Recall that if $\alpha $ and $\beta $ are in $E^*$, we say that $\alpha \preceq \beta $, if $\alpha $ is a \emph{prefix} of $\beta $, i.e.~if there exists $\gamma \in E^*$, such that $\alpha \gamma = \beta $.  Therefore, it follows from Proposition \ref{Prop:ProdIdemp} that
\[f_\alpha \leq f_\beta \iff \beta \preceq \alpha.\]

Another easy consequence of Proposition \ref{Prop:ProdIdemp} is that, for any two elements $e,f \in \EGE$, one has that either $e\perp f$, or $e$ and $f$ are comparable.  It follows that
\[e\Cap f \Rightarrow e\leq f ,\text { or } f\leq e.\]

A basic property that will play a role later on is the following.

\begin{definition}\label{Def:PseudoFree}
Let $(G,E,\varphi)$ be a self-similar graph under (\ref{equa:StanHypSSG}). We say that $(G,E,\varphi)$ is \emph{pseudo free} if whenever $(g,e)\in G\times E^1$ is such that $ge=e$ and $\varphi(g,e)=1$ then $g=1$.
\end{definition}

It is easy to see that pseudo-freeness implies this apparently stronger property.

\begin{proposition}[{\cite[Proposition 5.6]{EP2}}]\label{Prop:VarStar}
Suppose that $(G,E,\varphi)$ is pseudo-free. Then, for all $g_1,g_2 \in G$, and $\alpha \in E^*$, one has that
\[g_1\alpha =g_2\alpha \text { \ and \ } \varphi (g_1,\alpha )=\varphi (g_2,\alpha ) \Rightarrow g_1=g_2.\]
\end{proposition}  
  
For some future applications, we need the following definition.

\begin{definition}\label{Def:StronglyFixedPath}
If $g\in G$ and $\alpha $ is a finite path satisfying $g\alpha = \alpha  \text{ and } \varphi (g,\alpha )=1$,
we say that $\alpha $ is \emph{strongly fixed} by $g$.  In addition, if no proper prefix of $\alpha $ is strongly fixed by $g$, we say that $\alpha $ is a \emph{minimal} strongly fixed path for $g$.
\end{definition}

Thus,

\begin{proposition}[{\cite[Propositions 5.3 \& 5.5]{EP2}}]\label{Prop:StrFixElts}
Given $g$ in $G$, let $M_g$ be the set of all minimal strongly fixed paths for $g$.  Then the set of all strongly fixed paths for $g$ is given by
\[\bigsqcup\limits_{\mu \in M_g}\{\mu \gamma \mid \gamma \in E^*,\ d(\mu ) = r(\gamma )\}.\]
In particular, if $(G,E,\varphi)$ is pseudo-free and $\alpha\in E^*$ is strongly fixed for some $g\in G$, then $g=1$.
\end{proposition} 

Hence,

\begin{proposition}[{\cite[Propositions 5.7\& 5.8]{EP2}}]\label{Prop:LemDominateIdempotent}
Let $\alpha $, $\beta $ and $\gamma $ be finite paths in $E$, and let $g\in G$ be such that $d(\alpha )=g d(\beta )$, so that $s:=(\alpha ,g,\beta )$ is a general nonzero element of $\SGE $ and $e:=(\gamma ,1,\gamma )$ is a general nonzero idempotent element of $\SGE $.  Then $e\leq s$, if and only
\begin{enumerate}
\item $\alpha =\beta $,
\item $\gamma =\alpha \tau $, for some finite path $\tau $,
\item $\tau $ is strongly fixed by $g$.
\end{enumerate}
Thus, $\SGE$ is an $E^*$-unitary semigroup if and only iff $(G,E,\varphi)$ is pseudo-free.
\end{proposition}

\subsection{Representing $C^*(E)$ and $G$ into $\OGE $}

Here, we study natural representations of the graph $C^*$-algebra $C^*(E)$ and of the group $G$ in $\OGE $.

Given that
\[\{p_x \mid x \in E^0\}\cup \{s_e \mid e \in E^1\}\]
is a Cuntz-Krieger $E$-family, the universal property of the graph $C^*$-algebra $C^*(E)$ \cite{Raeburn} provides for the existence of a $\ast$-homomorphism
\[\iota : C^*(E)\to \OGE ,\]
sending the canonical Cuntz-Krieger $E$-family of $C^*(E)$ to the corresponding one within $\OGE $. It is easy to see that

\begin{proposition}[{\cite[Proposition 11.1]{EP2}}]\label{Prop:PropGraphSeaStarisSubalg}
The $\ast$-homomorphism $\iota $ above is injective.
\end{proposition}

With respect to the injectivity of the representation of $G$ into $\OGE $, we have to work a bit more to obtain a result in the line of Proposition \ref{Prop:PropGraphSeaStarisSubalg}.

\begin{lemma}[{\cite[Lemma 11.2]{EP2}}]\label{Lem:LempiInjImplyUInjUinj}
Let $\pi : \SGE \rightarrow \OGE $ and $u:G\rightarrow \OGE $ be the natural maps. If $\pi $ is injective, then so is $u$.
\end{lemma}

Then, it can be shown that

\begin{proposition}[{\cite[Corollary 11.7]{EP2}}]\label{Prop:CorpiIsInj}
If $(G,E, \varphi )$ is pseudo-free, then $\pi : \SGE \rightarrow \OGE $ is injective. Hence, if $(G,E, \varphi )$ is pseudo-free, then $u:G\rightarrow \OGE $ is injective.
\end{proposition}

Proposition \ref{Prop:CorpiIsInj} provides the best situation possible, as the next example shows:

\begin{example}
Let $E$ be the graph with only one vertex and one edge, and let $G$ be any noncommutative group. Fix the trivial action of $G$ on $E$, and let $\varphi $ be the $1$-cocycle of $G$ defined by $\varphi (g, a)=1$ for every $g\in G, a\in  E^1$. Then, it is easy to see that $\OGE \cong C^*(E)\cong C(\mathbb {T})$, which is a commutative $C^*$-algebra, so it cannot contain any faithful copy of $G$.
\end{example}

\section{$\OGE$ as groupoid algebra}\label{Sect:GroupoidPicture}

One of the main tools for understanding the structure of $\OGE$ is to represent it as a (full) groupoid $C^*$-algebra of a second countable ample groupoid. To do this, we use the classical approach due to Exel \cite{Exel1}.

The first step is to show that $\OGE$ is universal for tight representations of $\SGE$; roughly speaking, given an inverse semigroup $S$ and a unital ring $R$, a $\ast$-representation $\phi:S\rightarrow R$ is said to be tight if whenever $\{f_1,\dots ,f_n\}$ is a cover of $e$ (that is, for any $0\ne g \leq e$ there exists $1\leq j\leq n$ such that $g\cdot f_j\ne0$), then $\phi(e)=\bigvee_{i=1}^n \phi(f_i)$. For a precise definition of the concept, see \cite[Definition 13.1]{Exel1} (also \cite{ExelCover}).

\begin{theorem}[{\cite[Proposition 6.2 \& Theorem 6.3]{EP2}}]\label{Theorem:TightRep}
Let $(G,E,\varphi)$ be a self-similar graph under (\ref{equa:StanHypSSG}). Then:
\begin{enumerate}
\item The map
\[\pi : \SGE \to \OGE ,\]
defined by $\pi (0)=0$, and
\[\pi (\alpha ,g,\beta ) = s_\alpha u_gs_\beta ^*,	\]
is a tight representation.
\item Let $A$ be a unital $C^*$-algebra and let $\rho :\SGE \to A$ be a tight representation.  Then there exists a
unique unital *-homomorphism $\psi :\OGE \to A$, such that the diagram
 
\special{em:linewidth 0.4pt} \unitlength 1mm \linethickness{0.4pt}
\begin{picture}(80.00,30.00)
\put(82.00,3.00){\makebox(0,0){$A$}}
\put(55.00,25.00){\makebox(0,0){$\SGE$}}
\put(85.00,25.00){\makebox(0,0){$\OGE$}}
\put(83.00,22.00){\vector(0,-1){15.00}}
\put(70.00,27.00){\makebox(0,0){$\pi$}}
\put(90.00,15.00){\makebox(0,0){$\psi$}}
\put(65.00,10.00){\makebox(0,0){$\rho$}}
\put(62.00,25.00){\vector(1,0){16.00}}
\put(60.00,22.00){\vector(1,-1){19.00}}
\end{picture}

commutes.
\end{enumerate}
\end{theorem}

Hence,

\begin{corollary}[{\cite[Corollary 6.4]{EP2}}]\label{Cor:UniversalTightAlgebra}
Let $(G,E,\varphi)$ be a self-similar graph under (\ref{equa:StanHypSSG}). Then:
\begin{enumerate}
\item $\OGE $ is isomorphic to the (full) $C^*$-algebra of the tight groupoid of the inverse semigroup $\SGE$, denoted $C^*(\mathcal{G}_{\text{tight}}(\SGE))$.
\item If $K$ is any field, then $\Oo_{G,E}^K$ is isomorphic to the Steinberg $K$-algebra of the tight groupoid of the inverse semigroup $\SGE$, denoted $A_ K(\mathcal{G}_{\text{tight}}(\SGE))$ (see e.g. \cite{HPSS} for further details).
\end{enumerate}
\end{corollary}

Thus, we need to define $\mathcal{G}_{\text{tight}}(\SGE)$, and describe its topology. The groupoid $\CG$ is the quotient (by the so-called germ relation) of the transformation groupoid of the action of $\SGE$ on the space of tight filters over the idempotent semi-lattice $\EGE $ of $\SGE $ (see e.g. \cite[Section 4]{Exel1}). Since we need a clear picture of this construction, let us be more precise in this description.

By an \emph{infinite path} in $E$ we shall mean any infinite sequence of the form
\[\xi =\xi _1\xi _2\ldots ,\]
where $\xi _i \in E^1$, and $d (\xi _i) = r (\xi _{i+1})$, for all $i$.  The set of all infinite paths will be denoted by $E^\infty $.  Given an infinite path
\[\xi =\xi _1\xi _2\ldots \in E^\infty ,\]
and an integer $n\geq 0$, denote by $\xi\vert_ n$ the finite path of length $n$ given by
\[\xi\vert_ n = \left \{
\begin{array}{cc}
\xi _1\xi _2\ldots \xi _n, & \hbox { if } n\geq 1, \\
r (\xi _1), & \hbox { if } n=0.
\end{array}
\right.
\]

\begin{proposition}[{\cite[Proposition 8.1]{EP2}}]\label{Prop:ActionOnInfWords}
There is a unique action
\[(g,\xi ) \in G\times E^\infty \mapsto g\xi \in E^\infty\]
of $G$ on $E^\infty $ such that
\[{(g\xi )}\vert_n = g (\xi\vert_ n),\]
for every $g \in G$, $\xi \in E^\infty $, and $n \in \N$.
\end{proposition}

Given an infinite path $\xi \in E^\infty $, we may look at the subset
\[\mathcal{F}_\xi = \{f_{\xi\vert_n}\mid n \in \N\} \subseteq \EGE ,\]
which turns out to be an ultrafilter over $\EGE $. Denoting the set of all ultrafilters over $\EGE $ by $\widehat{\EGE}_\infty $, one may also show \cite [Proposition 19.11]{Exel1} that the correspondence
\[\xi \in E^\infty \mapsto \mathcal{F}_\xi \in \widehat{\EGE}_\infty\]
is bijective, and we use it to identify $ E^\infty$ and $ \widehat{\EGE}_\infty$. Furthermore, this correspondence may be proven to be a homeomorphism when $E^\infty $ is equipped with the product topology.

Since $E$ is finite, $E^\infty $ is compact by Tychonov's Theorem, and consequently so is $\widehat{\EGE}_\infty $. Being the closure of $\widehat{\EGE}_\infty $ within $\widehat{\EGE}$ the space $\widehat{\EGE}_\text{tight} $ formed by the tight filters \cite[Theorem 12.9]{Exel1}, it coincides with $\widehat{\EGE}_\infty$.

Identifying $\widehat{\EGE}_\text{tight} $ with $E^\infty $, as above, we can transfer the canonical action of $\SGE $ from the former to the latter, resulting in the following: for each element $(\alpha ,g,\beta ) \in \SGE $, we associate the partial homeomorphism of $E^\infty $ whose domain is the \emph{cylinder}
\[ Z(\beta) := \{\eta \in E^\infty \mid \eta =\beta \xi , \hbox { for some } \xi \in E^\infty \},\]
and which sends each $\eta =\beta \xi \in Z(\beta)$ to $\alpha g\xi $, where the meaning of ``$g\xi $'' is as in
Proposition \ref{Prop:ActionOnInfWords}.

As before, we do not use any special symbol to indicate this action, using module notation instead:
\[(\alpha ,g,\beta )\eta = \alpha g\xi\, , \forall (\alpha ,g,\beta )\in \SGE\, , \forall \eta =\beta \xi \in Z(\beta).\]
Then, we have the transformation groupoid
\[\SGE \times E^{\infty} :=\{(\alpha,g,\beta ; \xi) \mid \xi\in Z(\beta) \} \]
with operation $(s, \xi) \cdot (t, s\xi)=(st, \xi)$ and $(s,\xi)^{-1}= (s^*, s\xi)$.

On this groupoid, we can define a germ relation, compatible with the groupoid operations: $(s; \xi)\sim (t; \omega)$ if and only if $\xi=\omega$ and there exists $e\in \xi$ with $e\leq s^*s, t^*t$ such that $se=te$. If we denote the equivalence class of $(s;\xi)$ by $[s; \xi]$, then the tight groupoid is
\[\CG:=  \left \{ [\alpha, g,\beta; \xi] \mid  (\alpha ,g,\beta ) \in \SGE ,\ \xi \in Z(\beta) \right \} .\]

Now, let us look at the topology of $\CG$. Recall from \cite [Proposition 4.14]{Exel1} that, if $\mathcal{S} $ is an inverse semigroup acting on a locally compact Hausdorff topological space $X$, then the corresponding groupoid of germs, say $\mathcal{G} $, is topologized by means of the
basis consisting of sets of the form
\[\Theta (s,U),\]
where $s \in \mathcal{S} $, and $U$ is an open subset of $X$, contained in the domain of the partial homeomorphism attached to $s$ by the given action.  Each $\Theta (s,U)$ is in turn defined by
\[\Theta (s,U) = \left\{[s,x]\in \mathcal{G} \mid x\in U\right\}.\]
Notice that, by definition, $\Theta (s,U)$ is an open bisection, and it is compact if so is $U$.

If we restrict the choice of the $U$'s above to a predefined basis of open sets of $X$, e.g.~the collection of all
cylinders in $E^\infty $ in the present case, we evidently get the same topology on the groupoid of germs.  Therefore, we see that a basis for its topology consists of sets of the form
\[\Theta (\alpha ,g,\beta ;\gamma ) := \left\{(\alpha, g,\beta ;\xi) \in \CG : \xi \in Z(\gamma) \right\}.\]
A simple computation shows that, if it is not empty, then it coincides with $\Theta (\alpha ,g,\beta ;\beta )$. Thus, we redefine
\[\Theta (\alpha ,g,\beta) := \left\{(\alpha, g,\beta ;\xi) \in \CG : \xi \in Z(\beta) \right\},\]
and we conclude

\begin{proposition}[{\cite[Proposition 9.4]{EP2}}]
The collection of all sets of the form $\Theta (\alpha ,g,\beta )$, where $(\alpha ,g,\beta )$ range in $\SGE $, is a basis of open compact bisections for the topology of $\CG$.
\end{proposition}

Now, we are ready to give characterizations of the topological properties enjoyed by $\CG$ in terms of the action of $\SGE$ on $E^\infty$ (see \cite{EP2} for a thorough account of the connections).

\subsection{Topological properties of $\CG$}

We are mainly interested in characterizing the topological properties of $\CG$ connected with determining the simplicity and pure infiniteness of $\OGE\cong C^*(\CG)$.

\subsubsection{\underline{Hausdorff property for $\CG$}}

The main result is a characterization, associated to the fact that an \'etale groupoid $\mathcal{G}$ is Hausdorff if and only if its unit space $\mathcal{G}^{(0)}$ is closed \cite[Proposition 3.10]{EP2}.

\begin{theorem}[{\cite[Theorem 12.2]{EP2}}]\label{Thm:HausdorffGroupoid}
Let $(G,E,\varphi)$ be a self-similar graph under (\ref{equa:StanHypSSG}). Then, the following statements are equivalent:
\begin{enumerate}
\item $\CG $ is Hausdorff.
\item for every $g$ in $G$, there are at most finitely many minimal strongly fixed paths for $g$,
\end{enumerate}
\end{theorem}

As a natural consequence, we have

\begin{corollary}[{\cite[Proposition 12.1]{EP2}}]
Let $(G,E,\varphi)$ be a self-similar graph under (\ref{equa:StanHypSSG}). If it is pseudo-free, then $\CG$ is a Hausdorff groupoid.
\end{corollary}

\subsubsection{\underline{Minimality for $\CG$}}

We study conditions under which $\CG$ is minimal. To this end, we need to describe two relations among the vertices in $E^0$ which are relevant for the question at hand.  First of all, let us say
that
\[x \rightharpoonup y\]
if there exists a path $\alpha $ in $E^*$ such that $d (\alpha )=x$ and $r (\alpha )=y$.  Notice that this relation is reflexive and transitive, but this is neither symmetric nor antisymmetric, and hence it is neither an equivalence relation nor an order relation. 

The other relation we have in mind is simply the orbit relation, defined by
\[x\sim y\]
if there exists $g$ in $G$ such that $gx=y$.  Unlike ``$\rightharpoonup $'', it is well known that ``$\sim $'' is an equivalence relation.

We may then consider the smallest transitive relation extending both ``$\rightharpoonup $'' and ``$\sim $'', by saying that vertices $x$ and $y$ are related when one may find a sequence of vertices $x_0,x_1,\ldots ,x_{2n}$ such that
\[x =x_0 \rightharpoonup x_1 \sim x_2 \rightharpoonup x_3 \sim \ldots  \sim x_{2n-2} \rightharpoonup x_{2n-1} \sim x_{2n} = y.\hspace{1truecm} (\dagger)\]

Then, we have

\begin{proposition}[{\cite[Proposition 13.2]{EP2}}]\label{Prop:TwoRelations}
Let $x$ and $y$ be vertices in $E^0$.  Then the following statements are equivalent;
\begin{enumerate}
\item there exists a vertex $u$ such that $x \rightharpoonup u \sim y$,
\item there exists a vertex $v$ such that $x \sim v \rightharpoonup y$.
\end{enumerate}
\end{proposition}

\begin{definition}
Given $x$ and $y$ in $E^0$, we say that
\[x\gg y\]
if the equivalent conditions of Proposition \ref{Prop:TwoRelations} are satisfied.
\end{definition}

Observe that ``$\gg $'' coincides with the relation $(\dagger)$, thanks to Proposition \ref{Prop:TwoRelations}, and hence it is clearly transitive. It is also evident that ``$\gg $'' is reflexive but, again, it is neither symmetric nor antisymmetric. But it is possible to turn it into a partial order by
identifying elements whenever antisymmetry fails: two vertices $x$ and $y$ in $E^0$ will be called equivalent, in symbols
\[x\approx y\]
whenever $x\gg y$ and $y\gg x$.  Writing $[x]$ for the equivalence class of each $x$ in $E^0$, the set of all
equivalent classes, namely
\[\frac{E^0}{\approx} = \big \{[x]: x\in E^0\big \}\]
becomes a partially ordered set via the well-defined order relation.
\[[x] \geq  [y] \iff x \gg y.\]

\begin{definition}
Let $(G,E,\varphi)$ be a self-similar graph under (\ref{equa:StanHypSSG}). We say that:
\begin{enumerate}
\item $E$ is \emph{$G$-transitive} if, for any two vertices $x$ and $y$ in $E^0$, one has that $x\gg y$,
\item $E$ is \emph{weakly $G$-transitive} if, given any infinite path $\xi $, and any vertex $x$ in $E^0$, there is some vertex $v$ along $\xi $ such that $v\gg x$.
\end{enumerate}
\end{definition}

The notion of $G$-transitivity generalizes the well-known notion of transitivity in graphs.  When it holds, $E^0$ has a single equivalence class. On the other hand, weak $G$-transitivity is inspired by the notion of cofinality for graphs.

\begin{theorem}[{\cite[Theorem 13.6]{EP2}}]\label{Thm:CharacMinimal}
Given $(G,E,\varphi )$ under (\ref{equa:StanHypSSG}), one has that the following statements are equivalent:
\begin{enumerate}
\item the action of $\SGE $ on $E^\infty $ is irreducible,
\item $\CG $ is minimal,
\item $E$ is weakly $G$-transitive.
\end{enumerate}
 If, in addition, we have that $E$ does not have sinks, then the above conditions are also equivalent to:
\hspace{.2truecm}\rm{(4)} $E$ is $G$-transitive.
\end{theorem}

\subsubsection{\underline{Effectiveness for $\CG$}}

We discuss conditions under which $\CG$ is an effective groupoid, a condition that is intimately related to the action of $\SGE$ on $E^\infty $ being topologically free.  The reader is referred to
\cite[Section 4]{EP2} for the definition of the notion of topologically free actions of inverse semigroups. 

\begin{definition}$\mbox{ }$ 
\begin{enumerate}
\item A \emph{circuit}\footnote{Circuits are also called loops or cycles in the graph $C^*$-algebra literature.  Our preference for \emph{circuits} comes from the fact that it is the terminology of choice among graph theorists and also because in the established graph theory terminology, the word \emph{loop} refers to a single edge whose source and range coincide.} is a finite path $\gamma \in E^*$ of nonzero length such that $d (\gamma )= r (\gamma )$.
\item A \emph{$G$-circuit} is a pair $(g,\gamma )$, where $g\in G$, and $\gamma \in E^*$ is a finite path of nonzero length such that $d(\gamma )= g r (\gamma )$.
\end{enumerate}
\end{definition}

Thus, a $G$-circuit needs a little help from the group to close it up.  Notice that a (usual) circuit $\gamma $ may be concatenated infinitely many times, producing an infinite path
\[\xi =\gamma \gamma \gamma \ldots\]
Moreover, if
\[s = \left(\gamma ,1,d (\gamma )\right),\]
then, regarding the standard action of $\SGE $ on $E^\infty $, it is easy to see that $s\xi =\xi$, that is to say, $\xi $ is a fixed point for $s$.  It is also possible to create fixed points from $G$-circuits as
follows:

\begin{proposition}[{\cite[Proposition 14.2]{EP2}}]\label{Prop:GcircuitGivesFixed}
Given a $G$-circuit $(g,\gamma )$, define a sequence $\{\gamma ^n\}_{n\geq 1}$ of finite paths, and a sequence $\{g_n\}_{n\geq 1}$ of group elements, recursively by $\gamma ^1=\gamma $, \ $g_1=g$, and
\[\left\{
\begin{array}{ccl}
\gamma ^{n+1} &=& g_n\gamma ^n, \\
g_{n+1} &=& \varphi (g_n,\gamma ^n),
\end{array}
\right.\]
for all $n\geq 1$. Then
\begin{enumerate}
\item $d (\gamma ^n)=r (\gamma ^{n+1})$, for all $n\geq 1$,
\item the concatenation \ $\xi =\gamma ^1\gamma ^2\gamma ^3\ldots $ \ is a well defined infinite path,
\item for every finite path $\beta $ such that $d (\beta )=r (\gamma )$, one has that $s:= (\beta \gamma ,g,\beta )$ lies in $\SGE $, and $\beta \xi $
is a fixed point for $s$.
\end{enumerate}
\end{proposition}

The above method does not give us all fixed points of every single element $s$ in $\SGE $, but in certain cases it does:

\begin{proposition}[{\cite[Proposition 14.3]{EP2}}]\label{Prop:DescribeSomeFixedPoints}
Given $s := (\alpha ,g,\beta )$ in $\SGE $, suppose that $|\alpha |>|\beta |$.  Then, regarding the action of $\SGE $ on $E^\infty $, one has that:
\begin{enumerate}
\item $s$ admits at most one fixed point,
\item if $s$ admits a fixed point $\zeta $, then there is a $G$-circuit $(g,\gamma )$ such that $\alpha =\beta \gamma $, and $\zeta $ coincides with the fixed point $\beta \xi $ mentioned in Proposition \ref{Prop:GcircuitGivesFixed}(3), constructed from $(g,\gamma )$.
\end{enumerate}
\end{proposition}

We will eventually be interested in determining the conditions under which the action of $\SGE $ on $E^\infty $ is topologically free, so the fixed points that will really interest us are the interior ones. Under the conditions of Proposition \ref{Prop:DescribeSomeFixedPoints}, when there is at most one fixed point, the existence of interior fixed points hinges on whether or not the unique fixed point is isolated in $E^\infty $.  we now introduce certain concepts designed to study isolated fixed points. Recall that our graph $E$ has no sources, which means that $r^{-1}(x)\ne \emptyset$ for every vertex $x$.

\begin{definition}$\mbox{ }$
\begin{enumerate}
\item We shall say that a vertex $x$ in $E^0$ is a \emph{simple vertex} if $r^{-1}(x)$ is a singleton.
\item Given a path $\gamma =\gamma _1\gamma _2\ldots \gamma _n$ in $E^*$, where each $\gamma _i$ is in $E^1$, we say that $\gamma $ \emph{has no entry} if $d (\gamma _i)$ is a simple vertex for every $i=1,\ldots ,n$.
\item If the above condition fails, we say that $\gamma $ \emph{has an entry}.
\end{enumerate}
\end{definition}

\begin{proposition}[{\cite[Proposition 14.4]{EP2}}]\label{Prop:IsolatedIffNoEntry}
Under the conditions of Proposition \ref{Prop:DescribeSomeFixedPoints}(2), let $(g,\gamma )$ be the $G$-circuit and let $\zeta $ be the fixed
point for $s$ mentioned there.  Then the following statements are equivalent:
\begin{enumerate}
\item $\zeta $ is an isolated point in $E^\infty $,
\item $\gamma $ has no entry.
\end{enumerate}
\end{proposition}

Since we are interested in topologically free actions, we would like to avoid isolated fixed points, and hence we will be interested in situations where every $G$-circuit has an entry.  However, given that we are working with finite graphs only, the action of $G$ on $E$ is not relevant in this respect.  In precise terms, what we mean is that

\begin{proposition}[{\cite[Proposition 14.6]{EP2}}]\label{Prop:TudoIgual}
Given $(G,E,\varphi )$ under (\ref{equa:StanHypSSG}), the following statements are equivalent:
\begin{enumerate}
\item every $G$-circuit has an entry,
\item every circuit  has an entry.
\end{enumerate}
\end{proposition}

Since we need to use the finiteness of $E$ in a very strong way above, to extend this theory to infinite graphs we might have to distinguish between conditions $(1)$ and $(2)$ in Proposition \ref{Prop:TudoIgual}. We should point out that a graph in which every circuit has an entry is usually said to satisfy condition (L).

Now, we have work to do to treat the remaining case $|\alpha |=|\beta |$.

\begin{proposition}[{\cite[Proposition 14.7]{EP2}}]\label{Prop:FixForG}
Let $s:=(\alpha ,g,\beta )\in \SGE $, with $|\alpha |=|\beta |$, and suppose that $s$ admits a fixed point.  Then:
\begin{enumerate}
\item $\alpha =\beta $,
\item the fixed points of $s$ in $E^\infty $ are precisely the elements of the form $\zeta =\beta \xi $, where $\xi $ is an infinite path such that $r (\xi ) =d (\beta )$, and $g\xi =\xi $.
\end{enumerate}
\end{proposition}

The conclusion of Proposition \ref{Prop:FixForG} is that when $|\alpha |=|\beta |$, understanding the fixed points for $s$ requires understanding the fixed points for the action of $g$ on $E^\infty $.  One may easily describe such fixed points in terms of the action of $G$ on $E$ and the cocycle $\varphi $, but apparently there is no smart way to control each and every one of them.  However, since our main interest is in studying topological freeness, we need only focus on large (meaning open) sets of fixed points:

\begin{proposition}[{\cite[Proposition 14.8]{EP2}}]\label{Prop:GIntFixPoints}
Suppose that $s:=(\alpha ,g,\alpha )$ lies in $\SGE $, and that $\zeta $ is an interior fixed point for $s$. Then, there is a finite path $\gamma $, such that:
\begin{enumerate}
\item $g\gamma =\gamma $,
\item $d (\alpha ) = r (\gamma )$,
\item $\zeta \in Z(\alpha \gamma)$,
\item the group element $h: = \varphi (g,\gamma )$ pointwise fixes $Z(d (\gamma ))$ (i.e., every point in $Z(d (\gamma ))$ is fixed by $h$).
\end{enumerate}
Conversely, if $\gamma $ is any finite path satisfying (1), (2) and (4), then every $\zeta \in Z(\alpha \gamma)$ is a (necessarily interior) fixed point for $s$.
\end{proposition}

Looking for conditions under which the standard action of $\SGE $ on $E^\infty $ is to\-po\-lo\-gi\-cal\-ly free, one should probably worry about group elements fixing whole cylinders, as in Proposition \ref{Prop:GIntFixPoints}(4). The following notion is designed to pinpoint situations in which whole cylinders of the form $Z(x)$ are in fact fixed.

\begin{definition}\label{Def:DefineSlack}
Given $g\in G$, and $x\in E^0$, we shall say that $g$ is \emph{slack} at $x$, if there is a non-negative integer $n$ such that all finite paths $\gamma $ with $r (\gamma )=x$, and $|\gamma |\geq n$, are strongly fixed by $g$, as defined in Definition \ref{Def:StronglyFixedPath}.
\end{definition}

Clearly, if $\gamma $ is strongly fixed by $g$, then $g$ fixes any finite path extending $\gamma $, and hence also all infinite paths in $Z(\gamma)$. If $g$ is slack at $x$, and if $n$ is as in Definition \ref{Def:DefineSlack}, notice that
\[Z(x )= \bigcup\limits_{r(\gamma )=x, |\gamma |=n} Z(\gamma) ,\]
and since each $\gamma $ that occurs above is strongly fixed by $g$, we have that $g$ pointswise fixes $Z(\gamma)$, and hence also the whole cylinder $Z(x)$.

Notice that a path of length zero, namely a vertex $x$, is never strongly fixed by a nontrivial group element $g$, because
\[ \varphi (g,x)=g \neq 1.\]
The concept of slackness above should therefore be seen as the best replacement for the notion of being strongly fixed in the case of a vertex.

We are now ready for a main result:

\begin{theorem}[{\cite[Theorem 14.10]{EP2}}]\label{Thm:MainTopFreeSGE}
Given $(G,E,\varphi )$ under (\ref{equa:StanHypSSG}), the action of $\SGE $ on $E^\infty $ is topologically free if and only if the following two conditions hold:
\begin{enumerate}
\item every $G$-circuit has an entry,
\item given a vertex $x$, and a group element $g$ fixing every infinite path in $Z(x)$, then necessarily $g$ is slack at $x$.
\end{enumerate}
\end{theorem}

\begin{remark}
If for any $g\in G\setminus \{1\}$ and for any $x\in E^0$ there exists $\eta \in Z(x)$ such that $g\eta\ne \eta$, then Theorem \ref{Thm:MainTopFreeSGE}(2) holds trivially. 
\end{remark}

\begin{corollary}[{\cite[Corollary 14.13]{EP2}}]\label{Cor:TopFreeUnderEStarUni}
If $(G,E,\varphi )$ is a self-similar graph under (\ref{equa:StanHypSSG}) and pseudo free, then the action of $\SGE $ on $E^\infty $ is topologically free if and only if the following two conditions hold:
\begin{enumerate}
\item every $G$-circuit has an entry (which is the same as saying that every circuit has an entry by Proposition \ref{Prop:TudoIgual}),
\item for every $g$ in $G$, with $g\neq 1$, and for every $x$ in $E^0$, there is at least one $\zeta $ in $Z(x)$ such that $g\zeta \neq \zeta $.
\end{enumerate}
\end{corollary}

An important case for the theory of self-similar groups is when $G$ acts faithfully (meaning that if $g\xi =\xi $, for all $\xi $ in $E^\infty $, then $g=1$) on $E^\infty $, and $E$ is a graph with a single vertex.

\begin{corollary}[{\cite[Corollary 14.14]{EP2}}]\label{Cor:TopFreeUnderSSGroup}
If $(G,E,\varphi )$ is a self-similar graph under (\ref{equa:StanHypSSG}), suppose moreover that:
\begin{enumerate}
\item $E$ has a single vertex, and at least two edges,
\item $G$ acts faithfully on $E^\infty $.
\end{enumerate}
Then the action of $\SGE $ on $E^\infty $ is topologically free.
\end{corollary}

Having understood topological freeness, an immediate consequence of \cite[Theorem 4.7]{EP2} is the following:

\begin{corollary}[{\cite[Corollary 14.15]{EP2}}]\label{Cor:TopFreeSameEssPrin}
If $(G,E,\varphi )$ is a self-similar graph under (\ref{equa:StanHypSSG}), then $\CG$ is effective if and only if Theorem \ref{Thm:MainTopFreeSGE}(1-2) holds.
\end{corollary}

\subsubsection{\underline{Local contractivity for $\SGE$}}

In \cite[Section 6]{EP2} local contractivity is studied for groupoids and for actions of inverse semigroups.  We will now use these results to characterize local contractivity for the tight groupoid associated to an inverse semigroup $\mathcal{S}$.

\begin{theorem}[{\cite[Theorem 15.1]{EP2}}]\label{Thm:LocContrVsTopFree}
If $(G,E,\varphi )$ is a self-similar graph under (\ref{equa:StanHypSSG}), one has that the following statements are equivalent:
\begin{enumerate}
\item $\SGE$ is a locally contracting inverse semigroup,
\item the action $\theta :\SGE\curvearrowright E^\infty $ is locally contracting,
\item $\CG$ is a locally contracting groupoid,
\item every circuit in $E$ has an entry.
\end{enumerate}
\end{theorem}

It is worth noticing that many results of \cite{EP2} used in the proof of Theorem \ref{Thm:LocContrVsTopFree}, such as \cite[Proposition 6.3]{EP2}, \cite[Theorem 6.5]{EP2} and \cite[Proposition 6.7]{EP2}, comparing local contractivity for groupoids, inverse semigroups, and actions, are either one-way implications only, or the converse depends on special conditions, which are largely fulfilled when the graph is \emph{finite}.  

However, it is well known that the above condition on circuits is not sufficient for local contractivity for the groupoid associated to \emph{infinite} graphs. In fact, in the case of infinite graphs, a characterization of local contractivity for the corresponding groupoid will certainly not follow from the fact that every circuit has an entry, since this is false for infinite graphs, as mentioned above.

Moreover, the condition on the existence of entries for circuits completely ignores the group $G$, but, again, a generalization to infinite graphs will probably depend on the action. We will revisit that in Section \ref{Sect:CountGraph}

\subsection{Simplicity and pure infiniteness for $\OGE $}

Now, we use the results in the previous subsections to characterize when $\OGE $ is simple and purely infinite. To this end, we first recall the known results for (algebraic) groupoid algebras.

In the case of Steinberg algebras, the basic result is the following: 

\begin{theorem}{\cite[Theorem 3.5]{St2}}
\label{Thm:SteinbergSimple}
Let $\mathcal{G}$ be an ample groupoid such that $\mathcal{G}^{(0)}$ is Hausdorff. If $A_K (\mathcal{G})$ is simple, then $\mathcal{G}$ is effective and minimal. The converse holds if $\mathcal{G}$ is Hausdorff.
\end{theorem}

When Hausdorff condition fails, the best result is 

\begin{theorem}[{\cite[Theorem 3.14]{CEPSS}}]\label{Thm:simple}
Let $\mathcal{G}$ be a second-countable, ample groupoid such that $\mathcal{G}^{(0)}$ is Hausdorff.
Then $A_K (\mathcal{G})$ is simple if and only if the following three conditions are satisfied: \begin{enumerate}
\item $\mathcal{G}$ is minimal,
\item $\mathcal{G}$ is effective, and
\item for every nonzero $f\in A_K (\mathcal{G})$, $\text{supp}(f)$\footnote{Here, given $f\in A_K (\mathcal{G})$, we define $\text{supp}(f):=\{x\in \mathcal{G} \mid f(x)\ne 0\}$.} has nonempty interior.
\end{enumerate}
\end{theorem}

As we can see with some examples related to self-similar graphs, condition $(3)$ in Theorem \ref{Thm:simple} can fail.

In the case of groupoid $C^*$-algebras, the condition is a bit more involved. Here, we no longer restrict our attention to ample groupoids; instead we consider arbitrary (second-countable, locally-compact) \'etale groupoids (with Hausdorff unit space).  Here we mainly deal with the ``open'' support,
(unconventionally) defined as
\[\text{supp}(f) = {\{x \in X\mid f(x) \neq 0\}}.\]
We denote the set of continuous
functions with compact support by
\[C_c(X) :=\{ f:X \to \C \mid f \text{ is continuous and $\exists$ compact $K$ such
that $f(x) = 0 \ \forall \  x \notin K$}\}.\]
We write $\mathcal{C}(\Grpd)$ for Connes' algebra of
functions $f : \Grpd \to \C$ linearly spanned by the spaces $C_c(U)$ for open
bisections $U$ contained in $\Grpd$.   We view a function in $C_c(U)$ as a function on $\Grpd$ by
defining it to be $0$ outside of $U$. In some papers and texts, this algebra $\mathcal{C}(\Grpd)$ is simply
denoted $C_c(\Grpd)$, but we avoid this because its elements are in general not continuous.

For $u \in \Grpd^{(0)}$, we write $L_u$ for the regular representation $L_u : \mathcal{C}(\Grpd) \to
\mathcal{B}(\ell^2(\Grpd_u))$ satisfying
\[L_u(f)\delta_\gamma = \sum_{\alpha \in \Grpd_{r(\gamma)}}
f(\alpha)\delta_{\alpha\gamma} \text{ for } f \in \mathcal{C}(\Grpd).\]
By definition, $C^*_r(\Grpd)$ is the
completion of the image of $\mathcal{C}(\Grpd)$ under $\bigoplus_{u \in \Grpd 0} L_u$.

Recall that if $\Grpd$ is an \'etale groupoid and $a \in C^*_r(\Grpd)$, then we can define a
function $j(a) : \Grpd \to \C$ by $j(a)(\gamma) = \big(L_{s(\gamma)}(a) \delta_{s(\gamma)}
\mid \delta_\gamma\big)$. It is routine to check that $j(f) = f$ for $f \in \mathcal{C}(\Grpd)$.
Since $\bigoplus_u L_u$ is faithful, the map $a \mapsto j(a)$
is injective. Write $B(\Grpd)$ for the vector space of all bounded functions $f : \Grpd \to \C$,
and regard $B(\Grpd)$ as a normed vector space under $\|\cdot\|_\infty$. For $a \in C^*_r(\Grpd)$
and $\gamma \in \Grpd$, we have
\[
|j(a)(\gamma)| = \big|\big(L_{s(\gamma)}(a)\delta_{s(\gamma)} \mid \delta_\gamma\big)\big| \le \|L_{s(\gamma)}(a)\| \le \|a\|,
\]
so $j$ is an injective norm-decreasing linear map from $C^*_r(\Grpd)$ to $B(\Grpd)$.

Then, the right characterization of simplicity is

\begin{theorem}[{\cite[Theorem 4.10]{CEPSS}}]\label{Thm:c*simple}
Let $\Grpd$ be a second-countable, locally compact, \'etale groupoid such that $\Grpd^{(0)}$ is Hausdorff.
\begin{enumerate}
\item If $C^*(\Grpd)$ is simple, then $C^*(\Grpd) = C^*_r(\Grpd)$, $\Grpd$ is effective and for every nonzero $a \in C^*_r(\Grpd)$, $\text{supp}(j(a))$ has nonempty interior.
\item If $C^*_r(\Grpd)$ is simple, then $\Grpd$ is minimal.
\item If $\Grpd$ is minimal and effective and for every nonzero $a \in C^*_r(\Grpd)$, $\text{supp}(j(a))$ has nonempty interior, then $C^*_r(\Grpd)$ is simple.
\end{enumerate}
In particular, if $C^*(\Grpd) = C^*_r(\Grpd)$, then $C^*(\Grpd)$ is simple if and only if $\Grpd$ is
minimal, effective and for every nonzero $a \in C^*_r(\Grpd)$, $\text{supp}(j(a))$ has nonempty interior.
\end{theorem}

Given that condition $C^*(\Grpd) = C^*_r(\Grpd)$ is central in determining the simplicity of the algebra, we need to establish sufficient conditions to ensure this equality. Obviously, the nuclearity of the algebra is equivalent to $C^*(\Grpd) = C^*_r(\Grpd)$, and a sufficient, but not necessary, condition to have nuclearity is to require $\Grpd$ to be amenable \cite{A}. We will delay the proof of this fact to a subsequent section, but it is enough to require the group $G$ in a self-similar graph $(G, E, \varphi)$ under (\ref{equa:StanHypSSG}) to be amenable to guarantee that $\OGE$ is nuclear \cite[Corollary 10.16]{EP2}.

The other natural condition is to guarantee that $\CG$ is Hausdorff. Under these requirements, we have a first result:

\begin{theorem}[{\cite[Theorem 16.1]{EP2}}]\label{Thm:ThmSimple}
If $(G,E, \varphi )$ is a self-similar graph under (\ref{equa:StanHypSSG}), $G$ is amenable, and that for every $g\in G$ there are at most finitely many minimal strongly fixed paths for $g$, then $\OGE $ is simple if and only if the following conditions are satisfied:
\begin{enumerate}
\item $E$ is weakly-$G$-transitive.
\item Every $G$-circuit has an entry.
\item Given a vertex $x$, and a group element $g$ fixing $Z(x)$ pointwise, then necessarily $g$ is a slack at $x$.
\end{enumerate}
\end{theorem}
The same result holds for the algebraic $K$-algebra $\Oo_{G,E}^K$ over any field $K$, the only difference being that, in this case, the amenability of $G$ does not play any role.

The proof is clear because first condition in Theorem \ref{Thm:ThmSimple} is equivalent to $\CG$ being Hausdorff, Theorem \ref{Thm:ThmSimple}(1) is equivalent to $\CG$ being minimal, and Theorem \ref{Thm:ThmSimple}(2-3) are equivalent to $\CG$ being effective, so Theorem \ref{Thm:c*simple} applies.

Moreover,

\begin{theorem}[{\cite[Theorem 16.2]{EP2}}]\label{Thm:Thmpi}
Let $(G,E, \varphi )$ be a self-similar graph under (\ref{equa:StanHypSSG}), and let $G$ be an amenable group. If $\CG$ is effective, then every hereditary subalgebra of $\OGE $ contains an infinite projection.
\end{theorem}

Recall that a unital ring $R$ is purely infinite simple if it is a simple, non-division ring, and for any $0\ne a\in R$ there exist $x,y\in R$ such that $xay=1$. Thus,

\begin{corollary}[{\cite[Corollary 16.3]{EP2}}]\label{Cor:CorolPinfSimple}
If $(G,E, \varphi )$ is a self-similar graph under (\ref{equa:StanHypSSG}), the group $G$ is amenable and
$\CG$ is Hausdorff, then, whenever $\OGE$ is simple, it is necessarily also purely infinite (simple).
\end{corollary}

One fundamental point is that given $(G,E, \varphi )$ is a self-similar graph under (\ref{equa:StanHypSSG}), even if $G$ is amenable, this does not imply that $\CG$ is Hausdorff, whence Theorem \ref{Thm:ThmSimple} no longer applies.

\begin{example}[{\cite[Subsection 5.4]{CEPSS}}]\label{Exam:KatsuraNoHausdorff}
we construct a Katsura algebra $\Oo_{A,B}$ such that its tight groupoid fails to be Hausdorff, but the algebra still enjoys being simple (the details are described in \cite[Subsection 5.4]{CEPSS}).

Consider the matrices
\begin{equation}\label{eq:matrices}
A = \left[\begin{matrix}
2&1&0\\
1&2&1\\
1&1&2
\end{matrix}\right]
\hspace{.5cm} \text{and} \hspace{.5cm}
B = \left[\begin{matrix}
1&2&0\\
2&1&2\\
0&2&1
\end{matrix}\right].
\end{equation}
Let $E_A = (E_A^0, E_A^1, r, s)$ be the directed graph whose incidence matrix is $A$. The matrix $B$ determines an action $\Z\curvearrowright E_A$ and a cocycle $\varphi: \Z\times E^1_A\to \Z$, as described after Remark \ref{Rem:KatsuraConditions}. This action and cocycle are described for $1\in\Z$ as follows:
\begin{align*}
1\cdot e_{ii}^0 = e_{ii}^1, \hspace{1cm}& \varphi(1, e_{ii}^0) = 0 \quad  \text{ for } i = 1, 2, 3,\\
1\cdot e_{ii}^1 = e_{ii}^0, \hspace{1cm}& \varphi(1, e_{ii}^1) = 1 \quad \text{ for } i = 1,2,3,\\
1\cdot e_{12} = e_{12}, \hspace{1cm}& \varphi(1, e_{12}) = 2, \\
1\cdot e_{21} = e_{21}, \hspace{1cm}& \varphi(1, e_{21}) = 2,\\
1\cdot e_{32} = e_{32}, \hspace{1cm}& \varphi(1, e_{32}) = 2, \\
1\cdot e_{23} = e_{23}, \hspace{1cm}& \varphi(1, e_{23}) = 2, \\
1\cdot e_{13} = e_{13}, \hspace{1cm}& \varphi(1, e_{13}) = 0.
\end{align*}
Let $(\Z, E_A,\varphi)$ be the Katsura triple associated to the matrices \eqref{eq:matrices}, and let $\Grpd_{\Z, E_A,\varphi}$ be the associated groupoid. Then $\Grpd_{\Z, E_A,\varphi}$ is minimal, effective and non-Hausdorff \cite[Lemmas 5.12 \& 5.14]{CEPSS}. But $\mathcal{S}_{\Z, E_A}$ satisfies an additional condition \cite[Lemma 5.15]{CEPSS} that guarantees the simplicity of $\OAB$ \cite[Theorem 5.16]{CEPSS}. The same result holds for the algebraic Katsura algebra $\Oo_{A,B}^K$ over any field $K$.
\end{example}

In general, it is possible to construct self-similar groups $(G,X)$ such that, even if $\Grpd_{(G,X)}$ is minimal and effective, $\Oo (G,X)$ fails to be simple. Or even more strange, $\Oo (G,X)$ is simple, but $K(G,X)$ fails to be simple for fields $K$ having nonzero characteristic.

\begin{example}[{\cite[Subsection 5.6]{CEPSS}}]\label{Exam:GrigorchukNonSimple}
Let $X = \{0,1\}$, and let $a, b, c, d$ be the length-preserving bijections of $X^*$ determined by the formulas
\[
\begin{array}{lll}
a\cdot(0w) = 1w ,&& c\cdot (0w) = 0(a\cdot w),\\
a\cdot(1w) = 0w, && c\cdot (1w) = 1(d\cdot w),\\
b\cdot(0w) = 0(a\cdot w), && d\cdot (0w) = 0w,\\
b\cdot(1w) = 1(c\cdot w), && d\cdot (1w) = 1(b\cdot w),
\end{array}
\]
for every $w\in X^*\cup X^\infty$. The group $G$ generated by the set of bijections $\{a,b,c,d\}$ is called the {\em Grigorchuk group}. We have that $(G, X)$ is a faithful self-similar action with restrictions given by
\[
\begin{array}{lll}
a\vert_0 = e,&& c\vert_0 = a,\\
a\vert_1 = e, && c\vert_1 = d,\\
b\vert_0 = a, && d\vert_0 = e,\\
b\vert_1 = c, && d\vert_1 = b.
\end{array}
\]
Each of the elements $b, c,$ and $d$ have infinitely many minimally fixed words, whence $\Grpd_{G,X}$ is a minimal, effective, non Hausdorff ample groupoid. According to \cite[Theorem 5.22]{CEPSS}, the Steinberg $K$-algebra  $A_K (\Grpd_{G,X})$ for any field $K$ of characteristic different from $2$ is simple, and so is the $C^*$-algebra $\Oo_{G,X}$.

But if $\text{char}(K)=2$, then $A_K (\Grpd_{G,X})$ is shown that it is not simple, by explicitly constructing a nonzero singular function $f\in A_K (\Grpd_{G,X})$ \cite[Corollary 5.26]{CEPSS}.
\end{example}

\begin{example}[c.f. {\cite{SS2}}]\label{Exam:GrigorchukErschler}
Let $X = \{0,1\}$, and let $a, b, c, d$ be the length-preserving bijections of $X^*$ determined by the formulas
\[
\begin{array}{lll}
a\cdot(0w) = 1w, && c\cdot (0w) = 0w,\\
a\cdot(1w) = 0w, && c\cdot (1w) = 1(d\cdot w),\\
b\cdot(0w) = 0(a\cdot w), && d\cdot (0w) = 0(a\cdot w)w,\\
b\cdot(1w) = 1(b\cdot w), && d\cdot (1w) = 1(c\cdot w),
\end{array}
\]
for every $w\in X^*\cup X^\infty$. The group $G$ generated by the set of bijections $\{a,b,c,d\}$ is called the {\em Grigorchuk-Erschler group}. We also have that $(G, X)$ is a faithful self-similar action with restrictions given by
\[
\begin{array}{lll}
a\vert_0 = e,&& c\vert_0 = e,\\
a\vert_1 = e, && c\vert_1 = d,\\
b\vert_0 = a, && d\vert_0 = a,\\
b\vert_1 = b, && d\vert_1 = c.
\end{array}
\]
Each of the elements $b, c,$ and $d$ have infinitely many minimally fixed words, whence $\Grpd_{G,X}$ is a minimal, effective, non Hausdorff ample groupoid. But the Steinberg $K$-algebra $A_K (\Grpd_{G,X})$ is not simple for any field $K$, and so is the $C^*$-algebra $\Oo_{G,X}$. 
\end{example}
These are clear examples that simplicity for these kind of algebras is far to be characterized using only topological properties of the underlying groupoid. Steinberg and Szak\'acs \cite{SS1, SS2} did a thorough study of this particular phenomenon, providing families of examples that enjoy properties similar to Examples \ref{Exam:GrigorchukNonSimple} and \ref{Exam:GrigorchukErschler}. 

\section{$\OGE$ as Cuntz-Pimsner algebra}

Inspired by Nekrashevych's paper \cite{N1}, we give a description of $\OGE $ as a Cuntz-Pimsner algebra \cite{Pimsner}.  With this we will be able to prove that $\OGE $ is nuclear and that $\CG$ is amenable when $G$ is an amenable group.  As before, we work with self-similar graphs $(G,E,\varphi)$ under (\ref{equa:StanHypSSG}).\vspace{.2truecm}

We begin by introducing the algebra of coefficients over which the relevant Hilbert bimodule, also known as a
correspondence, will later be constructed.

Since the action of $G$ on $E$ preserves length, we see that the set of vertices of $E$ is $G$-invariant, so we get an action of $G$ on $E^0$ by restriction.  By dualization, $G$ acts on the algebra $C(E^0)$ of complex valued functions on $E^0$ (notice that, since $E^0$ is a finite set, $C(E^0)$ is nothing but ${\C}^{|E^0|}$). We may therefore form the crossed-product $C^*$-algebra
\[A = C(E^0) \rtimes G.\]
Since $C(E^0)$ is a unital algebra, there is a canonical unitary representation of $G$ in the crossed product, which we denote by $\{v_g\}_{g\in G}$. On the other hand, $C(E^0)$ is also canonically isomorphic to a subalgebra of $A$ and we therefore identify these two algebras without further warning.

For each $x $ in $E^0$, we denote the characteristic function of the singleton $\{x \}$ by $q_x $, so
that 
\[\{q_x \mid x \in E^0\}\] 
is the canonical basis of $C(E^0)$, and thus $A$ coincides with the closed linear span of the set
\[\left\{q_x v_g\mid  x \in E^0, g\in G\right\}.\]
For later reference, notice that the covariance condition in the crossed product reads
\[v_g q_x = q_{gx } v_g \hspace{.2truecm} \forall x \in E^0\hspace{.2truecm} \forall g\in G.\]

Our next step is to construct a correspondence over $A$.  In preparation for this, we denote by $A^e $ the right ideal of $A$ generated by $q_{d (e )} $, for each $e\in E^1$.  In technical terms,
\[ A^e = q_{d (e)} A. \]
With the obvious right $A$-module structure, and the inner product defined by
\[ \langle y,z\rangle = y^*z\hspace{.2truecm} \forall y,z\in A^e,\]
one has that $A^e$ is a right Hilbert $A$-module. Notice that this is not necessarily a full Hilbert module since
$\langle A^e ,A^e \rangle $ is the two-sided ideal of $A$ generated by $q_{d(e)} $, which might be a proper ideal in some cases.

Since $A$ is spanned by the elements of the form $q_x v_g$, $A^e$ is spanned by the elements of the form
$q_{d(e)} q_x v_g$, but, since the $q$'s are mutually orthogonal, this is either zero or equal to $q_{d(e)} v_g$. Therefore, we see that
 \[ A^e = \overline {\hbox {span}} \{q_{d(e)} v_g \mid g\in G\}.\]

Introducing the right Hilbert $A$-module, which will later be given the structure of a correspondence over $A$, we define
\[ M = \bigoplus _{e \in E^1} A^e.\]

Observe that if $x$ is a vertex which is the source of many edges, say
\[d^{-1}(x) = \{e_1,e_2,\ldots ,e_n\},\]
then
\[A^{e_i} = q_{d (e_i)} A = q_xA, \]
for all $i$, so that $q_x A$ appears many times as a direct summand of $M$.  However, these copies of $q_x A$ should be suitably distinguished, according to which edge $e_i$ is being considered. On the other hand, notice that if $d^{-1}(x) = \emptyset $, then $q_x A$ does not appear among the summands of $M$.

Addressing the fullness of $M$, observe that
\[ \langle M,M\rangle = \sum\limits_{x \in E^0,d^{-1}(x) \neq \emptyset }A q_x A. \]
So, when $E$ has no \emph{sinks}, that is, when $d^{-1}(x)$ is nonempty for every $x$, it follows that $M$ is
full.

Given $e\in E^1$, the element $q_{d(e)}$, when viewed as an element of $A^e \subseteq M$, will play a very special role in what follows, so we give it a special notation, namely
\[t_e := q _{d(e)}.\]
There is a small risk of confusion here in the sense that, if $e_1 ,e_2 \in E^1$ are such that
\[x := d (e_1)=d (e_2), \]
then we assign $q_x $ to both $t_{e_1}$ and $t_{e_2}$.  However, the coordinate in which $q_x$ appears in $t_{e_i}$ is determined by the corresponding $e_i$. So, if $e_1\neq e_2$, then $t_{e_1}\neq t_{e_2}$.

In order to completely dispel any confusion, here is the technical definition:
\[t _e= (m_f )_{f \in E^1},\]
where
\[m_f = \left\{
\begin{array}{cc}
 q_{d(e)}, &\text {if } f =e , \\
 0, &\text{otherwise.}
\end{array}
\right.
\]

We should notice that $t_e q_{d(e)}= t_e$, and that any element $y\in M$ can be written uniquely as $y=\sum\limits_{e \in E^1} t_e y_e$, where each $y_e \in A^e $.

As the next step in constructing a correspondence over $A$, we have to define a certain $\ast$-homomorphism from $A$ to the algebra $\mathcal{L} (M)$ of adjointable linear operators on $M$. Since $A$ is a crossed product algebra, this will be accomplished once we produce a covariant representation $(\psi ,V )$ of the $C^*$-dynamical system $\left(C(E^0),G\right)$ on $M$.  We begin with the group representation $V$.

\begin{definition}
For each $g\in G$, let $V_g$ be the linear operator on $M$ given by
\[V_g\left(\sum\limits_{e \in E^1} t_e y_e \right) = \sum\limits_{e \in E^1} t _{g e } v_{\varphi (g,e )}y_e ,\]
whenever $y_e \in A^e $, for each $e\in E^1$.
\end{definition}

By the uniqueness of the representation of $y$ above, it is clear that $V_g$ is well defined.

\begin{proposition}[{\cite[Proposition 10.7]{EP2}}]
Each $V_g$ is a unitary operator on $M$. Moreover, the correspondence $g\mapsto V_g$ is a unitary representation of $G$.
\end{proposition}

In order to complete our covariant pair, we must now construct a $\ast$-homomorphism from $C(E^0)$ to $\mathcal{L} (M)$. With this in mind, we give the following:

\begin{definition}
For every $x\in  E^0$, let
\[M_x = \bigoplus _{e \in r^{-1}(x)} A^e ,\]
which we view as a complemented submodule of $M$. In addition, we let $Q_x$ be the orthogonal projection from $M$ to $M_x$, so that
\[
Q_x (t_e y) = [r(e)=x] t_e y\hspace{.2truecm} \forall e \in E^1 \hspace{.2truecm} \forall y\in A.\]
\end{definition}

Observe that the $Q_x$ are pairwise orthogonal projections and that $ \sum\limits _{x \in E^0} Q_x = 1$.

\begin{definition}
Let $\psi :C(E^0) \to \mathcal{L} (M)$ be the unique unital $\ast$-homomorphism such that
\[\psi (q_x) = Q_x\hspace{.2truecm} \forall x \in E^0.\]
\end{definition}

From our working hypothesis that $E$ has no sources we see that, for every $x\in E^0$, there is some $e \in E^1$ such that $r(e)=x $. So,
\[Q_x(t_e)=t_e ,\]
whence $Q_x \neq 0$.  Consequently $\psi $ is injective.

\begin{proposition}[{\cite[Proposition 10.10]{EP2}}]
The pair $(\psi, V)$ is a covariant representation of the $C^*$-dynamical system $\left(C(E^0),G\right)$ in $\mathcal{L} (M)$.
\end{proposition}

Then, there exists a $\ast$-ho\-mo\-mor\-phism
\[\Psi : C(E^0) \rtimes G \to \mathcal{L} (M),\]
such that
\[\Psi (q_x )=Q_x \hspace{.2truecm} \forall x \in E^0, \text{ and }\Psi (v_g)=V_g \hspace{.2truecm} \forall g\in G.\]

Equipped with the left-$A$-module structure provided by $\Psi $, we then have that $M$ is a correspondence over $A$.

Our next goal is to prove that $\OGE $ is naturally isomorphic to the Cuntz-Pimsner algebra associated to the correspondence $M$, which we denote by $\Oo _M$.  As a first step, we identify a certain Cuntz-Krieger $E$-family.

\begin{proposition}[{\cite[Proposition 10.12]{EP2}}]\label{Prop:RelInOM}
The following relations hold within $\Oo_M$:
\begin{enumerate}
\item For every $x \in E^0$, one has that $\sum_{e \in r^{-1}(x)}t_et_e^* = q_x $.
\item $\sum _{e \in E^1}t_et_e^* = 1$.
\item The set $ \{q_x \mid x \in E^0\}\cup \{t_e \mid e \in E^1\}$ is a Cuntz-Krieger $E$-family.
\end{enumerate}
\end{proposition}

Thus,

\begin{proposition}[{\cite[Proposition 10.13]{EP2}}]\label{Prop:OneMapForCuntzPimsner}
There exists a unique surjective $\ast$-ho\-mo\-mor\-phism
\[\Lambda :\OGE \to \Oo_M\]
such that $\Lambda (p_x ) = q_x $, $\Lambda (s_e ) = t_e $, and $\Lambda (u_g) = v_g$.
\end{proposition}

Let us now prove that $\Lambda $ is bijective by providing an inverse to it.  Since $A$ is the crossed product $C^*$-algebra $C(E^0) \rtimes G$, the universal property of $\OGE$ guarantees the existence of a $\ast$-homomorphism
\[\theta _A: A \to \OGE ,\]
sending $q_x$ to $p_x $, and $v _g$ to $u_g$. For each $e\in E^1$, consider the linear mapping
\[\theta _M : M \to \OGE ,\]
given, for every $m = (m_e)_{e \in E^1} \in M$, by
\[\theta _M(m)= \sum _{e\in E^1}s_e \theta _A(m_e ) \in \OGE .\]
Notice that $\theta _M (t_e )= s_e $, for all $s \in E^1$, because
\[\theta _M (t_e )=s_e \theta _A(q_{d(e)}) =s_e p_{d(e)} =s_e.\]

\begin{lemma}[{\cite[Lemma 10.14]{EP2}}]
The pair $(\theta _A,\theta _M)$ is a representation of the correspondence $M$ in the sense of \cite[Theorem 3.4]{Pimsner}, meaning that for all $y\in A$ and all $\xi ,\xi' \in M$,
\begin{enumerate}
\item $\theta _M(\xi )\theta _A(y) = \theta _M(\xi y),$
\item $\theta _A(y)\theta _M(\xi ) = \theta _M(y\xi ),$
\item $\theta _M(\xi )^*\theta _M(\xi ') = \theta _A(\langle \xi ,\xi '\rangle ).$
\end{enumerate}
\end{lemma}

It is well known \cite [Theorem 3.4]{Pimsner} that the Toeplitz algebra for the correspondence $M$, usually denoted $\mathcal{T}_M$, is universal for representations of $M$, so there exists a $\ast$-ho\-mo\-mor\-phism
 \[\Theta _0: \mathcal{T}_M\to \OGE ,\]
coinciding with $\theta _A$ on $A$ and with $\theta _M$ on $M$.

\begin{theorem}[{\cite[Theorem 10.15]{EP2}}]\label{Thm:CPPicture}
The map $\Theta _0$, defined above, factors through $\Oo_M$, providing a $\ast$-iso\-mor\-phism
\[\Theta :\Oo_M\to \OGE ,\]
such that $\Theta (q_x ) = p_x $, $\Theta (t_e)=s_e $, and $\Theta (v _g) = u_g$, for all $x \in E^0$, $e \in E^1$, and $g\in G$.
\end{theorem}

When $G$ is amenable, $C(E^0)  \rtimes G$ is nuclear.  Then, from Theorem \ref{Thm:CPPicture}, the fact that Toeplitz-Pimsner algebras over nuclear coefficient algebras are nuclear \cite[Theorem 4.6.25]{BO}, and so are quotients of nuclear algebras \cite [Theorem 9.4.4]{BO}, we obtain:

\begin{corollary}[{\cite[Corollary 10.16]{EP2}}]\label{Cor:Amenabuilidade}
 If $G$ is amenable then $\OGE $ is nuclear.
\end{corollary}

\begin{remark} 
Since $E^0$ is finite, the nuclearity of $C(E^0)  \rtimes G$ is equivalent to the amenability of $G$.  However, if the present construction is generalized to infinite graphs, one could produce examples of non-amenable groups acting amenably on $E^0$, in which case $C(E^0) \rtimes G$ would be nuclear. The proof of
Corollary \ref{Cor:Amenabuilidade} could then be adapted to prove that $\OGE $ is nuclear.
\end{remark}

\begin{corollary}[{\cite[Corollary 10.18]{EP2}}]\label{Cor:Newamenabuilidade}
If $G$ is amenable, then $\CG$ is an amenable groupoid.
\end{corollary}

\section{Self-similar graphs of countable graphs}\label{Sect:CountGraph}

In this section, we extend the construction of a $C^*$-algebra associated to a self-similar graph to the case of arbitrary countable graphs. We reduce the problem to the row-finite case with no sources by using a desingularization process. Finally, we characterize the properties associated with simplicity in this case.

The results of this section, contained in \cite{EPS}, are essential to extend the scope of the results obtained in \cite{EP2} for Katsura algebras over finite matrices to the general case; this guarantees, by \cite{Kat1}, that every Kirchberg algebra in the UCT is the full groupoid $C^*$-algebra of a second countable amenable ample groupoid.

\subsection{Extending self-similar graphs to countable graphs}

First, we look at the problem of extending our class to the case of countably infinite graphs with no restrictions (i.e. sources, sinks, and infinite receivers are admitted).

Notice that, if $E^0$ is countably infinite, then the algebra $C^*(E)$ is no longer unital. Thus, if we intend to get a unitary representation of the group $G$ associated to the algebra, we need to consider unitary representations of $G$ in the multiplier algebra $\mathcal{M}(\OGE)$. This approach does not give us an intrinsic definition of the object, but it is very helpful to fix it precisely. In our search for such a definition, we follow the model of Katsura algebras \cite{Kat1}, where the unitary associated to an element of $\Z$ is written in terms of partial unitaries associated to the projections $p_x$ for $x\in E^0$.\vspace{.2truecm}

\begin{noname}\label{basicdata2}
{\rm The basic data for our construction is a triple $(G,E,\varphi)$ composed of:
\begin{enumerate}
\item A countable directed graph $E=(E^0, E^1, r, s)$.
\item A discrete group $G$ acting on $E$ by graph automorphisms.
\item A 1-cocycle $\varphi: G\times E^1\rightarrow G$ satisfying the property
\[
\varphi(g,a)\cdot x=g\cdot x \text{ for every } g\in G, a\in E^1, x\in E^0.
\]
\end{enumerate}
We will refer to that triple as a self-similar graph under (\ref{equa:StanHypSSG}).}
\end{noname}
As remarked before, this condition can be weakened to
\[
\varphi(g,a)\cdot s(a)=g\cdot s(a) \text{ for every } g\in G, a\in E^1,
\]
(see e.g. \cite[Appendix A]{LRRW18}). Notice that all the results in Subsection \ref{SubSect:SSGclassic} are true for $(G,E, \varphi)$ as in (\ref{basicdata2}). So, we assume that we have such a triple $(G,E, \varphi)$.

In order to build our algebra, we again resort to a set of generators 
\[\{p_x : x\in E^0\}\cup\{s_a : a\in E^1\}\]
as in the case of classical self-similar graphs, generating a copy of $C^*(E)$ into our algebra. However, once we allow $E^0$ to be infinite, the lack of units prevents us from having unitary elements in our algebra, and hence the inclusion of the missing $u_g$ becomes problematic.

The strategy we follow to circumvent this problem is to look at the elements $\{u_g \colon g\in G\}$ into the multiplier algebra ${\mathcal{M}(\OGE)}$ of our algebra $\OGE$, and suppose that the assignment $g\mapsto u_g$ defines a unitary group representation $u:G\rightarrow \mathcal{M}(\OGE)$. Then, in the multiplier algebra, we shall have relations of the form
\[ u_gs_a=s_{ga}u_{\varphi(g,a)} \]
for every $g\in G, a\in E^*$. Since $\OGE$ is an ideal of ${\mathcal{M}(\OGE)}$, the product $u_{g, x}:=u_gp_x$ is an element of $\OGE$. In addition, the fact that the $p_x$ add up to the identity of ${\mathcal{M}(\OGE)}$ in the strict topology ensures that $\sum\limits_{x\in E^0}u_{g, x}=u_g$ strictly.  So, we add to our generators a collection of symbols
\[\{u_{g, x}: x\in E^0,  g\in G\},\]
which would play the role of proxies for the unitary elements we are prevented to consider.

Using the relations in ${\mathcal{M}(\OGE)}$ (inspired by the relations in the original definition), we conclude that each $u_{g,x}$ is a partial isometry satisfying:
\begin{enumerate}
\item $u_{g,x}u_{g,x}^*=p_{g\cdot x}$.
\item $u_{g,x}^*u_{g,x}=p_{x}$.
\item $u_{g, s(a)}s_a=s_{g\cdot s(a)}u_{\varphi (g,a), s(a)}$.
\item $u_{g,x}p_x=p_{g\cdot x}u_{g,x}$.
\end{enumerate}

If we take the above identities as relations of our definition of $\OGE$, we see that $\sum\limits_{x\in E^0}u_{g,x}$ indeed converges to an element $u_g\in \mathcal{M}(\OGE)$, and it is easy to see that it is unitary. Since we are interested in getting a unitary representation of $G$ in $\mathcal{M}(\OGE)$, we need to be sure that $u_gu_h=u_{gh}$ for every $g,h\in G$. A
simple computation shows that this occurs exactly when 
\[u_{gh, h^{-1}\cdot x}=u_{g,x}u_{h, h^{-1}\cdot x}\]
for every
$g,h\in G$ and every $x\in E^0$. In view of all these facts, we give the following definition:

\begin{definition}\label{Def:OGE_infty}
{\rm Given a triple $(G,E, \varphi)$ as in  (\ref{basicdata2}), we define $\OGE$ to be the universal $C^*$-algebra as follows:
\begin{enumerate}
\item \underline{Generators}:
\[\{p_x: x\in E^0\}\cup\{s_a: a\in E^1\} \cup \{u_{g,x}: g\in G,\, x\in E^0\}.\]
\item \underline{Relations}:
\begin{enumerate}
\item $\{p_x: x\in E^0\}\cup\{s_a: a\in E^1\}$ is a Cuntz-Krieger $E$-family in the sense of \cite{Raeburn}.
\item $u_{g,x}$ is a partial isometry with:
\begin{enumerate}
\item $u_{g,x}u_{g,x}^*=p_{g\cdot x}$,
\item $u_{g,x}^*u_{g,x}=p_{x}$,
\end{enumerate}
for every $g\in G, x\in E^0$.
\item $u_{gh, h^{-1}\cdot x}=u_{g,x}u_{h, h^{-1}\cdot x}$ for every $g,h\in G$ and $x\in E^0$.
\item $u_{g, s(a)}s_a=s_{g\cdot a}u_{\varphi(g,a), s(a)}$ for every $g\in G, a\in E^1$.
\item $u_{g,x}p_x=p_{g\cdot x}u_{g, x}$ for every $g\in G, x\in E^0$.
\end{enumerate}
\end{enumerate}
}
\end{definition}

\begin{remark}\label{Rem:OGE_infty}
{\rm
When $E^0$ is finite, Definition \ref{Def:OGE_infty} coincides with Definition \ref{Def:DefineOGE}. Mo\-re\-over, we have:
\begin{enumerate}
\item $u_g:=\sum\limits_{x\in E^0}u_{g,x}$ is a unitary in $\mathcal{M}(\OGE)$.
\item The map $u:G\rightarrow \mathcal{M}(\OGE)$ defined by the rule $g\mapsto u_g$ is a unitary $\ast$-representation of $G$.
\item $u_{g,x}=u_gp_x$ for every $g\in G$ and $x\in E^0$.
\item $u_gp_x=p_{gx}u_g$ for every $g\in G$ and $x\in E^0$.
\item $u_gs_a=s_{ga}u_{\varphi(g,a)}$ for every $g\in G$ and $a\in E^1$.
\end{enumerate}
}
\end{remark}

Following the method adopted in the classic case, we wish to associate an
abstract inverse semigroup to these data.  Even though the present
context of infinite graphs is significantly more general than before, nothing prevents us from interpreting Definition \ref{Def:inverseSemigroup} verbatim here:

\begin{definition}\label{Def:inverseSemigroup_infty}
{\rm Given a self-similar graph $(G,E,\varphi)$ as in (\ref{basicdata2}), we define a $\ast$-semigroup ${\SGE}$ as follows: as a set,
\[{\SGE}=\{ (\alpha,g,\beta): \alpha, \beta\in E^*,\ g\in G,\ s(\alpha)=g\cdot s(\beta)\}\cup \{ 0\}.\]
The multiplication operation is defined as in \cite[Definition
4.1]{EP2}, and the semilattice of idempotents turns out to be
\[\EGE=\{(\alpha , 1, \alpha) : \alpha \in E^*\}\cup \{ 0\}.\]
}
\end{definition}

With small adaptations, the previous results from the case of finite graphs hold here as well. Moreover, there is a semigroup representation
\[
\begin{array}{cccc}
 \pi :& \SGE & \longrightarrow   & \OGE  \\
 & (\alpha, g, \beta) & \mapsto  &   s_{\alpha}u_{g,s(\alpha)}s_{\beta}^*
\end{array}.
\]

The essential point is the following result:

\begin{theorem}[{\cite[Theorem 2.5]{EPS}}]\label{Thm:GroupoidRowFiniteNoSources}
Let $(G,E, \varphi)$ be a self-similar graph as in (\ref{basicdata2}) such that $E$ is a row-finite graph with no sources. Then:
\begin{enumerate}
\item The semigroup homomorphism $ \pi :\SGE \rightarrow \OGE $ is a universal tight representation of $\SGE$.
\item $\OGE\cong C^*_{\text{tight}}(\SGE)\cong C^*(\CG)$.
\end{enumerate}
\end{theorem}

\vspace{.2truecm}

In order to extend the results in \cite{EP2} to this context, we reduce ourselves to a situation where the graph $E$ is row-finite without sources. We show that this is possible via a ``desingularization'' process, inspired by the one developed in \cite{DrinTom} for graph $C^*$-algebras.

%%%%%%%%%%%%%%%%%%%%%%%%%%%%%%%%%%%%%%%%%%%%%%%%%%%%%%%%%%%

\subsection{Desingularizing self-similar graphs}\label{sec:desing}

Suppose that we have a self-similar graph $(G,E, \varphi)$ as in (\ref{basicdata2}), and let $F$ denote the desingularized graph of $E$ obtained in \cite{DrinTom}. In this section, we show that we can define an action $G\curvearrowright F$ and a $1$-cocycle $\widehat{\varphi}:G\times F\rightarrow G$ extending the ones in the original triple, such that $\OGE$ and $\mathcal{O}_{G,F}$ are strong Morita equivalent $C^*$-algebras.

\begin{remark}\label{rem:action_sourcesreceivers}
{\rm
Let $(G,E, \varphi)$ be a triple as in (\ref{basicdata2}), and let $x\in E^0$ be a vertex. Then:
\begin{enumerate}
\item If $x$ is a source, then so is $gx$ for every $g\in G$.
\item If $x$ is an infinite receiver, then so is $gx$ for every $g\in G$.
\end{enumerate}
In view of this, when defining the desingularization, we need to keep track of the fact that, for any vertex in the orbit of a singular vertex, we must define the tails added to it as orbit-connected parts of the graph.
}
\end{remark}

%%%%%%%%%%%%%%%%%%%%%%%%%

\subsubsection{Desingularizing a source}\label{ssec:source}

Now, we explain how to construct a triple $(G,F, \varphi)$ that desingularizes a source in a triple $(G,E, \varphi)$. \vspace{.2truecm}

Let $x\in E^0$ be a source, and let $\widehat{G}$ be a set of representatives of the orbits of $x$ under the action $G\curvearrowright E^0$. We add an infinite tail into $x$ and to all vertices in the orbit of $x$. Specifically, define
\[
F^1 = E^1 \bigsqcup \{e_{i,g}\}_{g\in \widehat{G}, i\geq 1}
\]
\[
F^0 = E^0 \bigsqcup \{v_{i,g}\}_{g\in \widehat{G}, i\geq 1}
\]
We extend the range map, the source map, and the action of $G$ to $F$ as follows:
\begin{enumerate}
	\item $s(e_{i,g})=r(e_{i+1,g}) = v_{i,g}$ for every $i\geq 1$, and $r(e_{1,g})=gx$, while $ge_{i, 1_G}=e_{i,g}$  for every $i\geq 1$.
\item For any $h\in G$ with $hx = x$, we let $he_{i, 1_G}=e_{i,1_G}$ for every $i\geq 1$.

\end{enumerate}
At this point, we have applied the Drinen-Tomforde desingularization construction in such a way that it is coherent with the action $G\curvearrowright E$.

The next step is to extend the $1$-cocycle. To do this, for any $h\in G$, $g\in\widehat{G}$ and any $i\geq 1$ we define $\widehat{\varphi} (h, e_{i,g}):=h$. Once this is done, what we have obtained is:
\begin{enumerate}
\item[(3)] A graph $F$ that extends $E$, constructed using the Drinen-Tomforde desingularization process.
\item[(4)] An action $G\curvearrowright F$ that extends the original action $G\curvearrowright E$.
\item[(5)] A $1$-cocycle $\widehat{\varphi}:G\times F\rightarrow G$ extending ${\varphi}:G\times E\rightarrow G$.
\end{enumerate}

Then, the self-similar graph $(G,F, \widehat{\varphi})$ is the desingularization of $(G,E, \varphi)$ at the source $x$.

%%%%%%%%%%%%%%%%%%%%%%%%%

\subsubsection{Desingularizing an infinite receiver}\label{ssec:emiter}

Now, we explain how to construct a triple $(G,F, \varphi)$ that desingularizes an infinite receiver in a triple $(G,E, \varphi)$. \vspace{.2truecm}

Let $x\in E^0$ be an infinite receiver, let $\widehat{G}$ be a set of representatives of the orbits of $x$ under the action $G\curvearrowright E^0$, and list $r^{-1}(x)=\{a_i\}_{i\geq 1}$. We remove all the $a_i$ and replace them with a tail into $x$ together with a connecting edge for each $a_i$, and then repeat this for each vertex in the orbit of $x$. Specifically, let
\[
F^1 = \left(E^1 \setminus \{a_i: i\geq 1\}\right) \bigsqcup \{e_{i,g}\}_{g\in \widehat{G}, i\geq 1}  \bigsqcup \{f_{i,g}\}_{g\in \widehat{G}, i\geq 1}
\]
\[
F^0 = E^0 \bigsqcup \{v_{i,g}\}_{g\in \widehat{G}, i\geq 1}
\]
We define $\alpha_{i,g}=e_{1,g}e_{2,g}\cdots e_{i-1,g}f_{i,g}$, and extend the range map, source map, and action of $G$ to $F$ as follows:
\begin{enumerate}
	\item For each $g\in \widehat{G}$ and $i\geq 1$, $s(e_{i,g})=r(e_{i+1,g}) = v_{i,g}$, $r(f_{i,g})=gr(e_{i,1_G})$ and $s(f_{i,g})=gs(a_i)$ for every $i\geq 1$, $r(e_{1,g})=gx$, while $ge_{i,1_G}=e_{i,g}$ and $gf_{i,1_G}=f_{i,g}$  for every $i\geq 1$.
	\item For any $h\in G$ such that $hx = x$, define $he_{i,1_G} = e_{i,1_G}$, and $hf_{i,1_G}=f_{i,1_G}$.
\end{enumerate}
Again, we have applied the Drinen-Tomforde desingularization construction in such a way that it is coherent with the action $G\curvearrowright E$.

The next step is to extend the $1$-cocycle. To do this, we define for any $g\in\widehat{G} ,h\in G$ and any $i\geq 1$:
\begin{enumerate}
\item[(3)] $\widehat{\varphi} (h, e_{i,g}):=h$.
\item[(4)] $\widehat{\varphi} (h, f_{i,g})=\varphi(h, a_j)$.
\end{enumerate}
In particular, $\widehat{\varphi} (h, \alpha_{i,g})=\varphi(h, a_j)$. Once this is done, what we have obtained is:
\begin{enumerate}
\item[(5)] A graph $F$ that extends $E$, constructed using the Drinen-Tomforde desingularization process.
\item[(6)] An action $G\curvearrowright F$ that extends the original action $G\curvearrowright E$.
\item[(7)] A $1$-cocycle $\widehat{\varphi}:G\times F\rightarrow G$ extending ${\varphi}:G\times E\rightarrow G$.
\end{enumerate}

Then, the self-similar graph $(G,F, \widehat{\varphi})$ is the desingularization of $(G,E, \varphi)$ on the infinite receiver $x$.

%%%%%%%%%%%%%%%%%%%%%%%%%

\subsubsection{The desingularization result}\label{ssec:DesingThm}

Now, we check that the results of Drinen and Tomforde about their desingularization process for $C^*(E)$ extend to this context.

A simple inspection shows that \cite[Lemmas 2.9 \& 2.10]{DrinTom} extend to our context. Now, we arrange the proof of \cite[Theorem 2.11]{DrinTom} in order to obtain the desired Morita equivalence between $\OGE$ and $\mathcal{O}_{G,F}$. We follow the notation of the proof of \cite[Theorem 2.11]{DrinTom}. Let $E$ be a graph with a singular vertex $v_0$, let
\[\{t_e,q_v \mid e\in F^1, v\in F^0\}\]
be the canonical set of generators for $C^*(F)$, and let
\[\{s_e,p_v \mid e\in E^1, v\in E^0\}\]
be the Cuntz-Krieger $E$-family constructed \textbf{inside} $C^*(F)$ in \cite[Lemma 2.9]{DrinTom}. Recall that $\{s_e,p_v\}$ is defined as follows:
\begin{enumerate}
\item For every $v\in E^0$, $p_v:=q_v$.
\item For every $e\in E^1$ such that $r(e)\in E^0_{\text{rg}}$, $s_e:=t_e$.
\item For every $e\in E^1$ such that $r(e)\not\in E^0_{\text{rg}}$, we have that $e=a_j$ for some $j\geq 1$, and thus $s_e:=t_{\alpha_j}$.
\end{enumerate}

Define $B:=C^*(\{s_e,p_v\})$ and $p:=\sum\limits_{v\in E^0}q_v\in \mathcal{M}(C^*(F))$. Then, \cite[Theorem 2.11]{DrinTom} shows that
\[C^*(E)\cong B\cong pC^*(F)p\]
and that $p\in \mathcal{M}(C^*(F))$ is a full projection.

Now, observe that:
\begin{enumerate}
\item For every $g\in G$ and $x\in E^0$ we have that $u_gp_x=u_gq_x=q_{gx}u_g=p_{gx}u_g$.
\item If $r(e)\in E^0_{\text{rg}}$, then $u_gs_e=u_gt_e=t_{ge}u_{\widehat{\varphi}(g,e)}=s_{ge}u_{{\varphi}(g,e)}$, since $\varphi$ and $\widehat{\varphi}$ match on $G\times E$.
\item If $r(e)\not\in E^0_{\text{rg}}$, then
\[u_gs_e=u_gt_{\alpha_j}=u_gt_{e_1}t_{e_2}\cdots t_{e_{j-1}}t_{f_j}=t_{ge_1}t_{ge_2}\cdots t_{ge_{j-1}}u_gt_{f_j}=\]
\[t_{ge_1}t_{ge_2}\cdots t_{ge_{j-1}}t_{gf_j}u_{\widehat{\varphi}(g, f_j)}=t_{g\alpha_j}u_{\widehat{\varphi}(g, \alpha_j)}=s_{ge}u_{\varphi(g,e)}\]
by definition of $\widehat{\varphi}$ in this case.
\end{enumerate}

Thus, the $C^*$-algebra isomorphism
\[
\begin{array}{cccc}
\Phi: & C^*(E) & \rightarrow   & C^*(\{s_e, q_x\})  \\
 & P_x & \mapsto  &  q_x \\
 & T_e & \mapsto  &   s_e
\end{array}
\]
satisfies
\begin{enumerate}
\item For every $g\in G$ and $x\in E^0$, $u_g\Phi(P_x)=\Phi(P_{gx})u_g$.
\item For every $g\in G$ and $e\in E^1$, $u_g\Phi(T_e)=\Phi(T_{ge})u_{\widehat{\varphi}(g,e)}$.
\end{enumerate}

Hence, by the universal property of $\OGE$, $\Phi$ extends to a $C^*$-algebra isomorphism
\[\Phi: \OGE \rightarrow  C^*(\{u_g, s_e, q_x\})\subseteq \mathcal{O}_{G,F}.\]
Moreover, following \cite[Theorem 2.11]{DrinTom}, we have a $C^*$-isomorphism
\[
\begin{array}{cccc}
\Psi: & \Phi(C^*(E)) & \rightarrow   & pC^*(F)p  \\
 & s_{\alpha}s_{\beta}^* & \mapsto  &  ps_{\alpha}s_{\beta}^*p \\
\end{array}.
\]
Notice that in $pC^*(F)p$, for all $g\in G$ and $x\in E^0$ we have $pu_gp_xp=pp_{gx}u_gp$, and since $u_gp=pu_g$ by the definition of $p$, this is equal to $pp_{gx}pu_g$. Since a homomorphism from $\OGE$ to $p\mathcal{O}_{G,F}p$ extending $\Psi$ should send $u_gp_x\mapsto pu_gpxp$ and $p_{gx}u_g \mapsto pp_{gx}pu_g$, such an extension will be compatible with the defining relations of $\OGE$. Similarly for every $g\in G$ and $a\in E^1$ we have that $u_gs_a$ and $s_{ga}u_{\widehat{\varphi}(g,a)}$ will map to the same element in $p\mathcal{O}_{G,F}p$. So, again by the universal property of $\OGE$, the isomorphism $\Psi$ extends to a $C^*$-algebra isomorphism $\OGE\cong p\mathcal{O}_{G,F}p$. Also, as in \cite[Theorem 2.11]{DrinTom}, $p\in \mathcal{M}(C^*(F))\subset \mathcal{M}(\mathcal{O}_{G,F})$ is a full projection.

Summarizing, we have:

\begin{theorem}[{\cite[Theorem 3.2]{EPS}}]\label{Thm:Desin}
Let $(G,E, \varphi)$ be a self-similar graph, and let $(G,F, \widehat{\varphi})$ be its desingularization. Then, there exists a full projection $p\in \mathcal{M}(\mathcal{O}_{G,F})$ such that $\OGE\cong p\mathcal{O}_{G,F}p$. In particular, $\OGE$ and $\mathcal{O}_{G,F}$ are strongly Morita equivalent.
\end{theorem}

Hence, up to Morita equivalence, we can assume that $E$ is a countable,
row-finite graph without sources, and thus $\OGE$ is a full groupoid
$C^*$-algebra by Theorem \ref{Thm:GroupoidRowFiniteNoSources}.

\begin{remark}\label{Rem:NoGroupoidPreserving} 
This highly roundabout way to achieve this result does not imply that the groupoid in point is the tight groupoid of $\SGE$. However, we can keep track of the construction to argue that the corresponding tight groupoid is Morita equivalent to $\CG$.
\end{remark}

%%%%%%%%%%%%%%%%%%%%%%%%%%%%%%%%%%%%%%%%%%%%%%%%%%%%%%%%%%%

\subsection{Characterizing properties of $\OGE$}\label{sec:properties}

Now, we are ready to extend the characterizations of various properties
of $\CG$ (and thus the simplicity of $\OGE$) obtained in Section \ref{Sect:GroupoidPicture} to
the case of triples $(G,E,\varphi)$ with $E$ a countable arbitrary graph;
in this sense, recall that properties such as Conditions (L) and (K), or
cofinality, are preserved through the desingularization process, as
shown in \cite{DrinTom}.

Since several arguments in \cite{EP2} use the fact that $\Et\cong
E^{\infty}$ when $E$ is a finite graph without sources, we need to prove
this fact in the case where $E$ is infinite. Due to Theorem
\ref{Thm:Desin} and the previous remark, we can assume without loss of
generality that $E$ is row-finite without sources. So, we keep that
fact in force for the remainder of the section.

%%%%%%%%%%%%%%%%%%%%%%%%%

\subsubsection{A technical issue}\label{ssec:Technical}

Suppose that $(G,E,\varphi)$ is as above, where $E$ is a row-finite
graph without sources. Then it is easy to see that $\Eu \cong
E^{\infty}$. But since $E$ is infinite, the space of filters
$\mathcal{E}_0$ is locally compact but not compact, and therefore we cannot
guarantee that $\Eu=\Et$; the only thing we know is that
$\Eu$ is a dense subspace of $\Et$. \vspace{.2truecm}

By carefully revisiting the characterizations of ultrafilters and tight filters \cite[Section 12]{Exel1}, we prove that

\begin{proposition}[{\cite[Proposition 4.1]{EPS}}]\label{Prop:UltraIsTight}
If $(G,E,\varphi)$ is a self-similar graph with $E$ row-finite graph without sources, then $\Et= \Eu$.
\end{proposition}

%%%%%%%%%%%%%%%%%%%%%%%%%

\subsubsection{The properties}\label{ssec:Properties}

Finally, we are ready to obtain the desired characterizations. Notice that, since the final aim is to characterize the simplicity of $\OGE$, and this property is Morita invariant, using Theorem \ref{Thm:Desin} we can reduce the problem to the case where $E$ is a row-finite graph without sources. Moreover, under this restriction, Proposition \ref{Prop:UltraIsTight} holds. So, we can use the arguments in \cite{EP2} involving actions $\SGE \curvearrowright E^{\infty}$ in this context. Thus, we can look at the results in Section \ref{Sect:GroupoidPicture} and fix the hypotheses so that they work in this context.

First, with respect to Hausdorffness of $\CG$, we have the following result.

\begin{theorem}[{\cite[Theorem 4.2]{EPS}}]\label{Thm:Hauss}
Let $(G,E,\varphi)$ a self-similar graph with $E$ a row-finite graph without sources. Then, the following statements are equivalent:
\begin{enumerate}
\item For every $g\in G$ and for every $x\in E^0$ there exists a finite number of minimal strongly fixed paths for $g$ with range $x$.
\item $\CG$ is Hausdorff.
\end{enumerate}
\end{theorem}

Next, we characterize minimality as follows:

\begin{theorem}[{\cite[Theorem 4.3]{EPS}}]\label{Thm:Minimal}
Let $(G,E,\varphi)$ a self-similar graph with $E$ being a row-finite graph without sources. Then, the following statements are equivalent:
\begin{enumerate}
\item The standard action $G\curvearrowright E^{\infty}$ is irreducible.
\item $\CG$ is minimal.
\item $E$ is weakly $G$-transitive.
\end{enumerate}
\end{theorem}

Finally, we characterize when the groupoid is effective, as follows.

\begin{theorem}[{\cite[Theorem 4.4]{EPS}}]\label{Thm:EssPrin}
Let $(G,E,\varphi)$ a self-similar graph with $E$ being a row-finite graph without sources. Then, the standard action $G\curvearrowright E^{\infty}$ is topologically free (equivalently, $\CG$ is effective) if and only if the following two conditions hold:
\begin{enumerate}
\item Every $G$-circuit has an entry.
\item Given a vertex $x\in E^0$ and a group element $g\in G$, if $g$ fixes $Z(x)$ pointwise then $g$ is slack at $x$.
\end{enumerate}
\end{theorem}

Hence, we conclude the following characterization of simplicity.

\begin{theorem}[{\cite[Theorem 4.5]{EPS}}]\label{Thm:Simple}
Let $(G,E,\varphi)$ be a self-similar graph with $E$ countable graph, and let $(G,F ,\widehat{\varphi})$ be its desingularization. If $G$ is amenable and $\mathcal{G}_{\text{tight}}(\Semi_{G,F})$ is Hausdorff, then $\OGE$ is simple if and only if the following conditions are satisfied:
\begin{enumerate}
\item $F$ is weakly $G$-transitive.
\item Every $G$-circuit in $F$ has an entry.
\item Given a vertex $x\in F^0$ and a group element $g\in G$, if $g$ fixes $Z_F(x)$ pointwise then $g$ is slack at $x$.
\end{enumerate}
\end{theorem}

We now turn to the question of when $\OGE$ is purely infinite. Unfortunately, when $E$ is an infinite graph, the equivalence in \cite[Theorem 15.1]{EP2} fails, because the implication $(iv)\Rightarrow (i)$ does not work, even for the underlying graph $C^*$-algebra $C^*(E)$. Fortunately, there is a condition, generalizing  \cite[Theorem 15.1(iv)]{EP2}, which allows to show an analog result, although the implication $(4)\Rightarrow (1)$ below (\cite[Theorem 16.1, $(4)\Rightarrow (1)$]{EP2}) remains unproved.

\begin{theorem}[{\cite[Theorem 4.6]{EPS}}]\label{Thm:PurInf}
Let $(G,E,\varphi)$ a self-similar graph with $E$ countable, row-finite graph with no sinks. Then, given the following statements:
\begin{enumerate}
\item For every $x\in E^0$ there exists $\alpha_x\in E^ *$ with $r(\alpha_x)=x$ and there exists a $G$-circuit $(g, \gamma)$ with entries such that $s(\alpha_x)=r(\gamma)$.
\item $\SGE$ is a locally contracting inverse semigroup.
\item The standard action $\theta: G\curvearrowright E^{\infty}$ is locally contracting.
\item $\CG$ is a locally contracting groupoid,
\end{enumerate}
we have that $(1)\Rightarrow (2) \Rightarrow (3) \Rightarrow (4)$.
\end{theorem}

In particular, if $\CG$ is minimal, effective, and $E$ contains at least one $G$-circuit, then $\CG$ is locally contracting. Hence,

\begin{corollary}[{\cite[Corollary 4.7]{EPS}}]\label{Cor:PurInf}
Let $(G,E,\varphi)$ a self-similar graph with $E$ countable graph, and let $(G,F ,\widehat{\varphi})$ be its desingularization. If $G$ is amenable, $\mathcal{G}_{\text{tight}}(\Semi_{G,F})$ is Hausdorff, and $(G,F,\widehat{\varphi})$ contains at least one $G$-circuit \footnote{That fact was remarked by H. Larki}, then whenever $\OGE$ is simple, it is necessarily also purely infinite (simple).
\end{corollary}

\section{Some generalizations of self-similar graphs}

In this section, we go through various generalizations of self-similar graphs, with the aim of connecting the original construction with Zappa-Sz\'ep products of categories with groupoids.

One, quite unexplored, generalization of the concept is due to B\'edos, Kwaliszewski, and Quigg \cite{BKQ}, where they extend the original definition to self-similar actions of \emph{locally compact} groups acting on \emph{topological} graphs. To define the corresponding algebras, they use a Cuntz-Pimsner algebra model. They also shall be acknowledged for coining the term \emph{Exel-Pardo algebras} for referring to self-similar graph algebras.\vspace{.2truecm} 

Other generalizations are the following:

\subsection{Self-similar $k$-graphs}

Graph $C^*$-algebras have been a successful source of knowledge for operator algebraists, as allowed to easily link properties of algebras with the combinatorial properties of ``simpler'' objects \cite{Raeburn}. But they also have limitations on the scope of the families of algebras that fall under their umbrella. So, quickly appeared new generalizations of directed graphs covering the gaps: ultragraphs, labelled graphs, etc.

One of the most successful constructions has been higher rank graphs (or $k$-graphs, for $k\in \N$). The difference between this construction and the previous ones is that to define it, Kumjian and Pask deliberately avoided using a graph language, and preferred a categorical approach, which will turn out to be useful in the sequel.

\begin{definition}[{\cite[Definitions 1.1]{KP}, \cite[Definitions 2.1]{RSY}}]\label{Def:k-graph}
Given $k\in \N$, a graph of rank $k$ (or $k$-graph) $(\Lambda, d)$ consists of a
countable category $\Lambda = (\text{Obj}(\Lambda),\text{Hom}(\Lambda), r, s)$ together with a functor $d : \Lambda \rightarrow \N^k$, called the degree map, which satisfies the \emph{unique factorization property} (UFP): for every $\lambda\in \Lambda$ and $m, n \in \N^k$ with $d(\lambda) = m + n$, there are unique elements $\mu, \eta \in \Lambda$ such that $\lambda=\mu\eta$,
$d(\mu) = m$ and $d(\eta) = n$. Elements $\lambda \in \Lambda$ are called paths. For $m \in \N^k$ and $v \in \text{Obj}(\Lambda)$, we define $\Lambda^m := \{\lambda \in \Lambda \mid d(\lambda) = m\}$ and $\Lambda^m(v) := \{\lambda \in \Lambda^m \mid r(\lambda) = v\}$. A morphism between two $k$-graphs $(\Lambda_1, d_1)$ and $(\Lambda_2, d_2)$ is a functor $f : \Lambda_1 \rightarrow \Lambda_2$ that respects the degree maps. $(\Lambda, d)$ is row-finite if for each $v \in \text{Obj}(\Lambda)$ and $m \in \N^k$, the set $\Lambda^m(v)$ is
finite; $(\Lambda, d)$ has no sources if $\Lambda^m(v)\ne \emptyset$ for all $v \in \text{Obj}(\Lambda)$ and $m \in \N^k$.
The unique factorization property says that there is a unique path of degree $0$ at each vertex, and hence allows us to identify $\text{Obj}(\Lambda)$ with $\Lambda^0$.
\end{definition}

Then, we can construct a $C^*$-algebra associated to this new object. 

\begin{definition}[c.f. {\cite[Definition 2.4]{LY}}]\label{Def:k-GraphAlgebra}
Let $\Lambda$ be a $k$-graph. Then, a family of partial isometries $\{S_\mu\}_\mu\in \Lambda$ in a $C^*$-algebra $B$ is called a Cuntz-Krieger $\Lambda$-family if:
\begin{enumerate}
\item $\{S_\mu\}_\mu\in {\Lambda^0}$ is a family of mutually orthogonal projections;
\item $S_{\mu\nu}=S_\mu S\nu$ if $s(\mu)=r(\nu)$;
\item $S_\mu^*S_\mu=S_{s(\mu)}$ for all $\mu\in \Lambda$; and
\item $S_v=\sum_{\mu\in v\Lambda^n}S_\mu S_\mu^*$ for all $v\in \Lambda^0$, $n\in \N^k$.
\end{enumerate}
The $C^*$-algebra $\Oo_\Lambda$ generated by a universal Cuntz-Krieger $\Lambda$-family $\{s_\mu\}_\mu\in \Lambda$ is called the $k$-graph $C^*$-algebra of $\Lambda$.
\end{definition}
\vspace{.2truecm}

Li and Yang \cite{LY} extended the definition of classical self-similar graphs from actions of groups on graphs to actions of groups on \emph{$k$-graphs}.

\begin{definition}[{\cite[Definition 3.1]{LY}}]\label{Def:Auto_k-G}
Let $\Cat$ be a $k$-graph. A bijection $\varphi: \Cat \to \Cat$ is an automorphism of $\Cat$ if:
\begin{enumerate}
\item $\varphi(\Cat^n)\subseteq \Cat^n$ for all $n\in \N^k$;
\item $s\circ \varphi =\varphi \circ s$ and $r\circ \varphi =\varphi \circ r$.
\end{enumerate}
If $G$ is a (discrete countable) group, an action of $G$ on $\Cat$ is a homomorphism from $G$ to $\text{Aut}(\Cat)$. For $g\in G$ and $\mu\in \Cat$, we write $g\cdot \mu:=\varphi (g)(\mu)$.
\end{definition}

Using that notion, they define their generalization.

\begin{definition}[{\cite[Definition 3.2]{LY}}]\label{Def:SSk-G}
Let $\Cat$ be a $k$-graph, and let $G$ be a group acting on $\Cat$. The action is self-similar if there exists a restriction map $G\time \Cat \to G$ defined by $(g,\mu)\mapsto g\vert_\mu$ such that:
\begin{enumerate}
\item $g\cdot (\mu\nu)=(g\cdot \mu)(g\vert_\mu \cdot \nu)$ for all $g\in G$, $\mu, \nu\in \Cat$ with $s(\mu)=r(\nu)$;
\item $g\vert_v=g$ for all $g\in G$, $v\in \Cat^0$;
\item $g\vert_{\mu\nu}={g\vert_\mu}\vert_\nu$ for all $g\in G$, $\mu, \nu\in \Cat$ with $s(\mu)=r(\nu)$;
\item $e\vert_\mu=e$ for $e\in G$ the neutral element, and all $\mu\in \Cat$;
\item $(gh)\vert_\mu=g\vert_{h\cdot \mu}h\vert_\mu$ for all $g,h\in G$, $\mu\in \Cat$. 
\end{enumerate}
Then, $\Cat$ is called a self-similar $k$-graph over $G$, and $G$ is called a self-similar group on $\Cat$.
\end{definition}

Furthermore, they define the universal $C^*$-algebra associated to it \emph{when $\Cat$ as a finite set of objects}:

\begin{definition}[{\cite[Definition 3.2]{LY}}]\label{Def:algebraSSk-G}
Let $\Cat$ with $\vert \Cat^0\vert < \infty$ be a self-similar $k$-graph over a group $G$. We define the self-similar $k$-graph $C^*$-algebra $\Oo_{G, \Cat}$ to be the universal unital $C^*$-algebra generated by a family of unitaries $\{u_g\}_{g\in G}$ and a Cuntz-Krieger family $\{s_\mu\}_{\mu\in \Cat}$ satisfying:
\begin{enumerate}
\item $u_{gh}=u_gu_h$ for all $g,h\in G$;
\item$u_gs_\mu=s_{g\cdot \mu}u_{g\vert_\mu}$ for all $g\in G$, $\mu\in \Cat$.
\end{enumerate}
\end{definition}

They showed that $\Oo_{G, \Cat}$ can be realized as a Cuntz-Pimsner algebra associated to a suitable product system \cite[Theorem 4.6]{LY}. Under the additional condition of $(G,\Cat)$ being pseudo-free, they construct a Deacoun-Renault groupoid $\Grpd_{G,\Cat}$ \cite[Definition 5.2]{LY}, and showed that, if \emph{$G$ is amenable}, then $\Oo_{G, \Cat}$ is isomorphic to the full groupoid $C^*$-algebra associated to $\Grpd_{G,\Cat}$ \cite[Theorem 5.10]{LY}. In this case, they characterize when the algebra $\Oo_{G, \Cat}$ is simple and purely infinite \cite[Section 6]{LY}.\vspace{.2truecm}

An interesting generalization was given by Larki \cite{Larki}. His approach is based on the association of a suitable inverse semigroup $\mathcal{S}_{G,\Cat}$ with the self-similar $k$-graph $(G,\Cat)$ (which generalizes the semigroup $\mathcal{S}_\Cat$ introduced in \cite{FMY}). This allows us to avoid restrictions about finiteness conditions on $\Cat$, or pseudo-freeness of $(G,\Cat)$. Since the notation in \cite{Larki} slightly changes from that in \cite{LY}, let us fix it explicitly:

\begin{definition}\label{defn2.3}
Let $\Cat$ be a finitely aligned $k$-graph and $G$ be a discrete group. A triple $(G,\Cat,\varphi)$ is called a \emph{self-similar $k$-graph} whenever the following properties hold:
\begin{enumerate}
  \item $G$ acts on $\Cat$ by a group homomorphism $g\mapsto \psi_g$. we write $g\cdot \mu$ for $\psi_g(\mu)$ to ease the notation.
  \item $\varphi:G\times \Cat\rightarrow G$ is a 1-cocycle for the above action $G$ on $\Cat$ such that for every $g\in G$, $\mu,\nu\in \Cat$ and $v\in \Cat^0$ we have:
  \begin{itemize}
    \item[(a)] $\varphi(gh,\mu)=\varphi(g,h \cdot \mu)\varphi(h,\mu)$,
    \item[(b)] $g\cdot (\mu\nu)=(g\cdot \mu)(\varphi(g,\mu) \cdot \nu)$,
    \item[(c)] $\varphi(g,\mu\nu)=\varphi(\varphi(g,\mu),\nu)$, and
    \item[(d)] $\varphi(g,v)=g$.
  \end{itemize}
\end{enumerate}
\end{definition}

If there is no ambiguity, we often write a self-similar $k$-graph $(G,\Cat,\varphi)$ as $(G,\Cat)$. \vspace{.2truecm}

\noindent {\bf Standing assumption:} Throughout this subsection, we work with self-similar $k$-graphs $(G,\Cat)$ over finitely aligned $k$-graphs $\Cat$.

\begin{definition}\label{defn2.4}
Let $(G,\Cat)$ be a self-similar $k$-graph. A \emph{$(G,\Cat)$-family} is a set
\[\{s_\mu:\mu\in \Cat\}\cup\{u_{v,g}:v\in\Cat^0,g\in G\}\]
in a $C^*$-algebra satisfying the following relations:
\begin{enumerate}
  \item $\{s_\mu:\mu\in\Cat\}$ is a Cuntz-Krieger $\Cat$-family,
  \item $u_{v,e_G}=s_v$ for all $v\in\Cat^0$,
  \item $u_{v,g}^*=u_{g^{-1} \cdot v ,g^{-1}}$ for all $v\in\Cat^0$ and $g\in G$,
  \item $u_{v,g}s_\mu=\delta_{v,g\cdot r(\mu)}s_{g\cdot \mu}u_{g\cdot s(\mu),\varphi(g,\mu)}$ for all $v\in\Cat^0$, $\mu\in\Cat$, and $g\in G$,
  \item $u_{v,g}u_{w,h}=\delta_{v,g\cdot w}u_{v,gh}$ for all $v,w\in\Cat^0$ and $g,h\in G$.
\end{enumerate}
The $C^*$-algebra $\Oo_{G\,\Cat}$ associated to $(G,\Cat)$ is the universal $C^*$-algebra generated by a $(G,\Cat)$-family $\{s_\mu,u_{v,g}\}$.
\end{definition}

Then, he defined his inverse semigroup: given a self-similar $k$-graph $(G,E, \varphi)$, we fix the set of triples $(\mu, g, \nu)\in \Cat \times G\times \Cat$ such that $s(\mu)=g\cdot s(\nu)$. 
\vspace{.2truecm}

\noindent {\bf Property $(\perp)$}: Two triples $(\mu,g,\nu)$ and $(\xi,h,\eta)$ are called {\it orthogonal}, denoted by $(\mu,g,\nu)\perp (\xi,h,\eta)$, if $\Cat^{\text{min}}(\mu,\xi)=\Cat^{\text{min}}(\nu,\eta)=\emptyset$.

\begin{definition}\label{defn3.1}
Let $\mathcal{S}_{G,\Cat}$ be the collection of all finite sets $F$ containing pairwise orthogonal triples $(\mu,g,\nu)$ satisfying $s(\mu)=g \cdot s(\nu)$. We equip $\mathcal{S}_{G,\Cat}$ with the multiplication
\begin{equation}\label{eq3.1}
EF:=\bigcup_{\substack{(\mu,g,\nu)\in E \\ (\xi,h,\eta)\in F}} \left\{ (\mu(g\cdot\alpha), \varphi(g,\alpha)\varphi(h,h^{-1}\cdot\beta), \eta(h^{-1}\cdot\beta)):(\alpha,\beta)\in\Cat^{\text{min}}(\nu,\xi) \right\}
\end{equation}
and the inverse
\[F^*:=\left\{(\nu,g^{-1},\mu):(\mu,g,\nu)\in F\right\}\]
for all $E,F\in \mathcal{S}_{G,\Cat}$.
\end{definition}

Then, he showed that $\mathcal{S}_{G,\Cat}$ is an inverse semigroup \cite[Proposition 3.6]{Larki}, and that the map
\[
\begin{array}{cccl}
\pi \colon & \mathcal{S}_{G,\Cat}  & \rightarrow   & \Oo_{G, \Cat}  \\
 & F & \mapsto  & \sum\limits_{(\mu, g, \nu) \in F}s_\mu u_{s(\mu), g}s_nu^* 
\end{array}
\]
is a universal tight representation, so that
\[\Oo_{G,\Cat}\cong C^*(\Grpd_{\text{tight}}(\mathcal{S}_{G,\Cat}))\]
(see \cite[Theorem 4.7]{Larki}).

Once this is at hand, he used the classical strategy \cite{EP1} to characterize the properties of the groupoid $\Grpd_{\text{tight}}(\mathcal{S}_{G,\Cat})$, and thus the simplicity of $\Oo_{G, \Cat}$, bypassing the technical restrictions needed in \cite{LY}.

\subsection{Self-similar actions of groupoids on graphs}

Laca, Ramagge, Raeburn, and Whittaker \cite{LRRW18} extended the definition of self-similar graphs to include \emph{faithful} actions of groupoids on the path space of \emph{finite} graphs. Since a groupoid $G$ has a large number of ``neutral'' elements (the units), we need to establish a suitable definition of the action of $G$ on a finite graph $E$. 

Given a graph $E$, and any $v\in E^0$, the subgraph $vE^*=\{\mu\in E^*\colon r(\mu)=v\} $
is a rooted tree with root $v\in E^0$, and $T_E=\bigsqcup_{v\in E^0} vE^*$ is a disjoint union of trees, or \emph{forest}. A \emph{groupoid} is a small category with inverses. Thus, a groupoid consists of a set $G^0$ of objects, a set $G$ of morphisms, two functions $c,d:G\to G^0$, and
a partially defined product $(g,h)\mapsto gh$ from $G^2:=\{(g,h): d(g)=c(h)\}$ to $G$ 
such that $(G,G^0,c,d)$ is a category and such that each $g\in G$ has an inverse $g^{-1}$. 

Since the groupoid operation is \emph{partial}, to define the action of the groupoid $G$ on $E$, we need to deal with partial automorphisms of $E^*$. These can be codified as partial automorphisms of the forest $T_E$, defined on each subtree and restricting to isomorphisms on its domain and range.

\begin{definition}[{\cite[Definition 3.1]{LRRW18}}]\label{defn: partial iso}
Suppose $E=(E^0,E^1,r,s)$ is a directed graph. A \emph{partial isomorphism} of $T_E$ consists of two vertices $v,w \in E^0$ and a bijection $g:vE^* \to wE^*$ such that
\begin{equation}
\label{defnpi}
g|_{vE^k}\text{ is a bijection onto }wE^k \text{ for }k \in \N, \text{ and } 
g(\mu e) \in g(\mu)E^1 \text{ for all } \mu e \in vE^*. 
\end{equation}
For $v\in E^0$, we write $\id_{v}:vE^* \to vE^*$ for the partial isomorphism given by
$\id_{v}(\mu)=\mu$ for all $\mu \in vE^*$. 
We denote the set of all partial isomorphisms of $T_E$ by $\text{PIso}(E^*)$.
We define domain and codomain maps $d,c\colon \text{PIso}(E^*) \to E^0$ so that $g:d(g)E^* \to c(g)E^*$; 
thus in~\eqref{defnpi} we have $d(g):=v$ and $c(g):=w$.
\end{definition}

The key point is

\begin{proposition}[{\cite[Proposition 3.2]{LRRW18}}]\label{prop:alltreeisos}
Suppose that $E=(E^0,E^1,r,s)$ is a directed graph with associated forest $T_E$. Then $(\text{PIso}(E^*),E^0,c,d)$ is a groupoid in which: the product is given by composition of functions, the identity isomorphism at $v\in E^0$ is $\id_v:vE^*\to vE^*$, and the inverse of $g\in \text{PIso}(E^*)$ is the inverse of the function $g:d(g)E^*\to c(g)E^*$.
\end{proposition}

Given a directed graph $E$ and a groupoid $G$ with unit space $E^0$, an \emph{action} of $G$ on the path space $E^*$ is  a (unit-preserving) groupoid homomorphism $\phi:G\to \text{PIso}(E^*)$; the action is \emph{faithful} if $\phi$ is one-to-one. 
If the homomorphism is fixed, we can write $g\cdot \mu$ for $\phi_g(\mu)$. 

\begin{definition}[{\cite[Definition 3.3]{LRRW18}}]
\label{faithful self-similar groupoid action}
Suppose $E=(E^0,E^1,r,s)$ is a directed graph and $G$ is a groupoid with unit space $E^0$ which acts \underline{faithfully} on $T_E$. The action is \emph{self-similar} if  for every $g\in G$ and $e\in d(g)  E^1$, there exists $h\in G$ such that
\begin{equation}
\label{selfsimilar groupoid defn}
g\cdot(e\mu)=(g\cdot e)(h\cdot \mu) \text{ for all } \mu\in s(e)E^*.
\end{equation}
Since the action is faithful, there is then exactly one such $h\in G$, and we write $g|_e:=h$.
\end{definition}

Definition~\eqref{faithful self-similar groupoid action} has some immediate and important consequences.

\begin{lemma}[{\cite[Lemma 3.4]{LRRW18}}]
\label{source and range observations}
Suppose $E=(E^0,E^1,r,s)$ is a directed graph and $G$ is a groupoid with unit space $E^0$ acting self-similarly on $T_E$. Then for $g,h\in G$ with $d(h)=c(g)$ and $e\in d(g) E^1$, we have
\begin{enumerate}
\item 
\label{source and range in SSG}
$d(g|_e)=s(e)$ and $c(g|_e)=s(g\cdot e)$,   
\item 
\label{equivariance}
$r(g\cdot e) = g\cdot r(e)$ and $s(g\cdot e) = g|_e\cdot s(e)$,
\item\label{resid}
if $g=\id_{r(e)}$, then $g|_e=\id_{s(e)}$, and
\item\label{multonedges}
$(hg)|_{e}=(h|_{g\cdot e})(g|_e)$.
\end{enumerate}
\end{lemma}

Notice that Lemma~\ref{source and range observations}\eqref{equivariance} implies that 
the source map may not be equivariant: $s(g\cdot e)\not=g\cdot s(e)$ in general. Indeed, $g\cdot s(e)$ will often not make sense: $g$ maps $d(g)E^*$ onto $c(g)E^*$, and $s(e)$ is not in $d(g)E^*$ unless $s(e)=d(g)$. This non-equivariance of the source map is one of the main points of difference with the original definition of self-similar action, and also differs from that in \cite{BKQ}.

With this new construction at hand, they construct the Toeplitz and Cuntz-Krieger algebras associated to this self-similar graph \cite[Section 4]{LRRW18}, and recover the Katsura algebras \cite[Example 7.7]{LRRW18} as a particular case. Moreover, it is shown \cite[Appendix A]{LRRW18} that if a classical self-similar graph $(G,E,\varphi)$ satisfies a certain faithfulness condition, then it gives rise to faithful self-similar groupoid actions; to be coherent with \cite{LRRW18} notation, let me write $(G,E, \sigma, \varphi)$ for $(G,E,\varphi)$ with action $\sigma: G \curvearrowright E$. 

\begin{proposition}[{\cite[Proposition A.2]{LRRW18}}]\label{EP trans group}
Suppose $(G,E,\sigma,\varphi)$ is a self-similar action of a group on a graph and let $G \times E^0:=\{(g,v)\}$ be the transformation groupoid with $d(g,v)=v$, $c(g,v)=\sigma_g^0(v)$, and unit space $E^0$, and consider the maps
\begin{align}
\label{EP action map}
(g,v) \cdot \mu &:= \sigma_g^*(\mu) \quad \text{ for $\mu \in vE^*$ and} \\
\label{EP restriction map}
(g,v)|_{\mu} &:= (\varphi^*(g,\mu),s(\mu)) \quad \text{ for $\mu \in vE^*$}.
\end{align}
Then \eqref{EP action map} defines a groupoid action of $G \times E^0$ on $T_E$ by partial isomorphisms. If the groupoid action is faithful, then $(G \times E^0,E)$ is a self-similar groupoid action.

Moreover, if $E$ has no sources, then the algebra $\Oo(K \times E^0,E)$ is isomorphic to the self-similar graph algebra $\OGE$.
\end{proposition}

\subsection{Self-similar graphs with a twist}

One topic in groupoid $C^*$-algebras theory which attracts a lot of interest is that of \emph{twisted groupoid $C^*$-algebras}, and their algebraic siblings, \emph{twisted Steinberg algebras} (see \cite{steintwist, resteintwist}, and references herewith).

Corti\~nas \cite{Willie} defined a twisted version of self-similar graphs using an algebraic appro\-ach, showed that it relates to the classical definition associated to twisted groupoid extensions, and studied simplicity, pure infiniteness, as well as $kk$-theoretical and homological results. Let us give a view of his work.\vspace{.2truecm}

Fix a self-similar graph $(G,E,\phi)$ (he termed it \emph{Exel-Pardo tuple}), and write $\cU(\ell)$ for the group of invertible elements of our ground ring $\ell$.  A \emph{twisted Exel-Pardo tuple} is an Exel-Pardo tuple $(G,E,\phi)$ together with a $1$-cocyle $c:G\times E^1\to \cU(\ell)$. Remark that 
\[
\phi_c:G\times E^1\to \cU(\ell)G\subset \cU(\ell[G]), \ \phi_c(g,e)=c(g,e)\phi(g,e)
\]
is a $1$-cocycle with values in the multiplicative group of the group algebra $\ell[G]$. We write $(G,E,\phi_c)$ for the twisted EP-tuple $(G,E,\phi,c)$. 

\begin{lemma}[{\cite[Lemma 2.3.1]{Willie}}]\label{lem:twistextend} Let $(G,E,\phi_c)$ be a twisted Exel-Pardo tuple. Then $c$ extends uniquely to a $1$-cocycle $G\times\cP(E)\to \cU(\ell)$ satisfying
\begin{equation}\label{eq:conditwist}
c(g,v)=1,\,\, c(g,\alpha\beta)=c(g,\alpha)c(\phi(g,\alpha),\beta)
\end{equation}
for all concatenable $\alpha,\beta\in\cP(E)$. 
\end{lemma}

\begin{remark}\label{rem:phic}
Lemma \ref{lem:twistextend} imply that, for extension of $\phi_c$ to a map $\ell[G]\times\cP(E)\to \ell[G]$, $\ell$-linear on the first variable, we have
\[
\phi_c(g,\alpha\beta)=\phi_c(\phi_c(g,\alpha),\beta).
\]
\end{remark}

Let $\cS$ be an inverse semigroup with inverse operation $s\mapsto s^*$, pointed by a zero element $\emptyset$. A (normalized) \emph{$2$-cocycle} with values in the pointed inverse semigroup $\cU(\ell)_\bu=\cU(\ell)\cup\{0\}$ is a map 
\[\omega:\cS\times\cS\to \cU(\ell)_\bu\] 
such that for all $s,t,u\in\cS$ and $x\in\cS\setminus\{\emptyset\}$, 
\begin{gather}\label{eq:2cocy1}
\omega(st,u)\omega(s,t)=\omega(s,tu)\omega(t,u),\\ 
\omega(\emptyset,s)=\omega(s,\emptyset)=0,\ \ \omega(x,x^*)=\omega(xx^*,x)=\omega(x,x^*x)=1.\label{eq:2cocy2}%\nonumber
\end{gather}
The \emph{twisted semigroup algebra} $\ell[\cS,\omega]$ is the $\ell$-module $\ell[\cS]=(\bigoplus_{s\in S}\ell s)/\ell \emptyset$ equipped with the $\ell$-linear multiplication $\cdot_\omega$ induced by
\begin{equation}\label{semitwist}
s\cdot_\omega t=st\omega(s,t).    
\end{equation}
A straightforward calculation shows that $\ell[\cS,\omega]$ is an associative algebra and that for every $s\in \cS$,
\begin{equation*}
s\cdot_\omega s^*\cdot_\omega s=s.
\end{equation*}

Let $(G,E,\phi_c)$ be a twisted Exel-Pardo tuple. On the inverse semigroup $\cS(G,E,\phi)$ we can define a map 
\begin{gather}\label{map:omega}
\omega:\cS(G,E,\phi)\times \cS(G,E,\phi)\to \cU(\ell)_\bu\\
\omega(0,s)=\omega(s,0)=0\ \ \forall s\in \cS(G,E,\phi),\nonumber\\
\omega((\alpha,g,\beta),(\gamma,h,\theta))=\left\{
\begin{array}{ll}
c(h,h^{-1}(\beta_1))& \text{ if } \beta=\gamma\beta_1,\\ 
c(g,\gamma_1)& \text{ if } \gamma=\beta\gamma_1,\\ 
0& \text{ otherwise.}
\end{array}\right. 
\nonumber
\end{gather}

\begin{lemma}[{\cite[Lemma 2.4.6]{Willie}}]\label{lem:omega}
The map $\omega$ of \eqref{map:omega} is a $2$-cocycle. 
\end{lemma}

Now, we proceed to define twisted Exel-Pardo algebras. First, we define the \emph{Cohn algebra} of the twisted Exel-Pardo tuple $(G,E,\phi_c)$ to be the quotient $C(G,E,\phi_c)$ of the free algebra on 
the set 
\[\{v,\ vg,\ e, \ eg,\ e^*,\ ge^*:\ \  v\in E^0,\ \ g\in G,\  \ e\in E^1\}, \] 
modulo the following relations:
\begin{gather}
v=v1,\ \ e=e1,\ \ 1e^*=e^*,\ eg=s(e)e(r(e)g),\ ge^*=(g(r(e))g)e^*s(e),  \label{eq:cohn0}\\
 e^*f=\delta_{e,f}r(e),\ (vg)wh=\delta_{v,g(w)}vgh, \label{eq:cohn1}\\
(vg)e=\delta_{v,g(s(e))}g(e)\phi_c(g,e),\ \ e^*vg=\delta_{v,s(e)}\phi_c(g,g^{-1}(e))(g^{-1}(e))^*.\label{eq:cohn2}
\end{gather} 
If $E^0$ is finite, then $C(G,E,\phi_c)$ is unital with unit $1=\sum_{v\in E^0}v$, and $g\mapsto \sum_{v\in E^0}vg$ is a group homomorphism $G\to \cU(C(G,E,\phi_c))$. From Proposition \ref{prop:cohnpres} it will follow that the latter is a monomorphism, and we identify $G$ with its image in $\cU(C(G,E,\phi_c))$. Thus, if $E^0$ is finite, then writing $\cdot$ for the product in $C(G,E,\phi_c)$, we have $vg=v\cdot g$ and $eg=e\cdot g$. However, if $E^0$ is infinite, we do not identify $G$ with any subset of $C(G,E,\phi_c)$.   

\begin{proposition}[{\cite[Proposition 3.1.4]{Willie}}]\label{prop:cohnpres}
Let $\omega$ be as in \eqref{map:omega}. There is an isomorphism of algebras
\[
\begin{array}{cccl}
\psi\colon & C(G,E,\phi_c) & \rightarrow  & \ell[\cS(G,E,\phi),\omega]  \\
 & vg & \mapsto  & (v,g,g^{-1}(v))  \\
 & eg & \mapsto  & (e,g,g^{-1}(r(e))) \\
 & ge^* &\mapsto   & (g(r(e)),g,e).
\end{array}
\]
\end{proposition}

Now, we take a suitable ideal of the Cohn algebra to construct the Exel-Pardo algebra. For each regular vertex $v$ of $E$ and each $g\in G$, we define
\begin{equation}\label{eq:qv}
q_vg:=vg-\sum_{s(e)=v}e\phi_c(g,g^{-1}(e))g^{-1}(e)^* \in C(G,E,\phi_c).
\end{equation}
As usual, we set $q_v=q_v1$; notice that $q_vg=q_v\cdot (vg)$. Recall that a vertex $v\in E^0$ is \emph{regular} if $0<|s^{-1}\{v\}|<\infty$; write $\reg(E)\subset E^0$ for the subset of regular vertices. We define a two-sided ideal
\begin{equation}\label{eq:kgec}
C(G,E,\phi_c)\vartriangleright \cK(G,E,\phi_c)=\langle q_v\mid v\in \reg(E)\rangle.
\end{equation}

\begin{proposition}[{\cite[Proposition 3.2.5]{Willie}}]\label{prop:kge} 
The $*$-algebra $\cK(G,E,\phi_c)$ is independent of $\phi_c$ up to canonical algebra isomorphism.
\end{proposition}

Let $(G,E,\phi_c)$ be a twisted EP-tuple. We define the \emph{twisted} Exel-Pardo algebra of $(G,E,\phi_c)$ to be the quotient algebra.
\[
L(G,E,\phi_c)=C(G,E,\phi_c)/\cK(G,E,\phi_c).
\]

Then,
\begin{proposition}[{\cite[Proposition 3.4.2]{Willie}}]\label{prop:lesubsetlge}
Let $(G,E,\phi_c)$ be any twisted EP tuple. Then the canonical homomorphism 
\[L(E)\to L(G,E,\phi_c)\] 
is injective. 
\end{proposition}

We relate this algebra with the twisted Steinberg algebra introduced in \cite{steintwist}. Let $(\cS,\emptyset)$ be a pointed inverse semigroup, and let $\omega:\cS\times\cS\to \cU(\ell)_\bu$ be a $2$-cocycle as above. Equip the smash product set
\begin{align*}
\cU(\ell)_\bu\land \cS=&\ \ \negthickspace \cU(\ell)_\bu\times \cS/\cU(\ell)_\bu\times\{\emptyset\}\cup \{0\}\times \cS\\
                            =&\ \ \negthickspace \cU(\ell)\times \cS/\cU(\ell)\times\{\emptyset\}.
\end{align*}
with the product
\[
(u\land s)(v\land t)=uv\omega(s,t)\land st.
\]
The result is a pointed inverse semigroup $\tilde{\cS}$, with inverse 
\[
(u\land s)^*=u^{-1}\land s^*.
\]
We write $0$ for the class of the zero element in $\tilde{\cS}$. The coordinate projection $\cU(\ell)\times \cS\to \cS$ induces a surjective semigroup homomorphism $\pi:\tilde{S}\to S$. Observe that $\pi^{-1}(\{\emptyset\})=\{0\}$ and that if $p\in\cS$ is a nonzero idempotent, we have a group isomorphism
\[
\pi^{-1}(\{p\})=\cU(\ell)\land p\cong \cU(\ell).
\]
Next, let $X$ be a locally compact, Hausdorff space, and suppose that there is an action of $\cS$ on $X$. Then $\tilde{\cS}$ also acts on $X$ through $\theta\circ\pi$, and so we can consider the groupoids of germs 
\[
\Grpd=\Grpd(\cS,X), \ \ \tilde{\Grpd}=\Grpd(\tilde{\cS},X).
\]
Equip $\cU(\ell)$ with the discrete topology and regard it as a groupoid over the one-point space. 
The map 
\begin{equation}\label{map:trivbund}
\tilde{\Grpd}\to \cU(\ell)\times \Grpd, \ \ [u\land s,x]\mapsto (u,[s,x])
\end{equation}
is a homeomorphism that maps a basic open set $[u\land s,U]$ to the basic open $\{u\}\times[s,U]$, and the sequence
\begin{equation}\label{seq:distwist}
\cU(\ell)\times\Grpd^{0}\to \tilde{\Grpd}\to \Grpd
\end{equation}
is a \emph{discrete twist} over the (possibly non-Hausdorff) \'etale groupoid $\Grpd$, in the sense of \cite[Section 2.3]{resteintwist}. Under the homeomorphism \eqref{map:trivbund}, the multiplication of $\tilde{\Grpd}$ corresponds to
\[
(u,[s,t(x)])(v,[t,x])=(uv\omega(s,t),[st,x]).
\]
The map
\begin{equation}\label{map:tildomega}
\tilde{\omega}:\Grpd^{(2)}\to\cU(\ell), \ \ \tilde{\omega}([s,t(x)],[t,x])=\omega(s,t)
\end{equation}
is a continuous, normalized $2$-cocycle in the sense of \cite[page 4]{steintwist}. 

If $\Grpd$ is an ample groupoid, its \emph{Steinberg} algebra $\cA(\Grpd)$ \cite[Proposition 4.3]{St} is the $\ell$-module spanned by all characteristic functions of compact slices, equipped with the convolution product
\begin{equation}\label{convo}
f\star g([s,x])=\sum_{t_1t_2=s}f[t_1,t_2(x)]g[t_2,x].
\end{equation}

Following \cite[Definition 5.2]{St}, we say that the action of $\cS$ on $X$ is \emph{ample} if $X$ is Boolean and $\dom(s)$ is a compact clopen subset for all $s\in\cS$. In this case, the characteristic function $\chi_{[s,\dom(s)]}\in\cA(\Grpd)$ for every $s\in\cS$, and the $\ell$-module map
\begin{equation}\label{map:ellstoag}
\begin{array}{ccl}
\ell[\cS] & \to &\cA(\Grpd)\\ 
s & \mapsto & \chi_{[s,\dom(s)]}
\end{array}
\end{equation}
is a homomorphism of algebras. One may also equip the $\ell$-module $\cA(\Grpd)$ with the \emph{twisted convolution product}
\begin{align}\label{convotwist}
f\star_{\omega} g([s,x])=&\sum_{t_1t_2=s}\tilde{\omega}([t_1,t_2(x)], [t_2,x])f[t_1,t_2(x)]g[t_2,x]\\
=&\sum_{t_1t_2=s}\omega(t_1,t_2)f[t_1,t_2(x)]g[t_2,x].\nonumber
\end{align}
The result is an algebra $\cA(\Grpd,\tilde{\omega})$. By \cite[Proposition 4.3]{St}, this is the same as the twisted Steinberg algebra defined in \cite[Proposition 3.2]{steintwist} under the assumption that $\Grpd$ is Hausdorff. Moreover, by \cite[Corollary 4.25]{steintwist}, the latter also agrees with the Steinberg algebra of the discrete twist \eqref{seq:distwist}. 
Comparison of formulas \eqref{convo}, \eqref{convotwist} and \eqref{semitwist} tells us that, whenever the $\ell$-linear map \eqref{map:ellstoag} is an algebra homomorphism $\ell[\cS]\to \cA(\Grpd)$, it is also an algebra homomorphism $\ell[\cS,\omega]\to \cA(\Grpd,\tilde{\omega})$. 

\begin{lemma}[{\cite[Lemma 4.1.7]{Willie}}]\label{lem:ellStostein}
Let $\cS$ be an inverse semigroup, let $\omega:\cS\times\cS\to \cU(\ell)$ be a $2$-cocycle, $\cS \curvearrowright X$ be an ample action, and let $\Grpd$ be the groupoid of germs. Let $\tilde{\omega}:\Grpd^{(2)}\to\cU(\ell)$ be as \eqref{map:tildomega}. Assume that the algebra homomorphism \eqref{map:ellstoag} is surjective with kernel $\cK$. Then $\cK\triqui\ell[\cS,\omega]$ is an ideal, and $\ell[\cS,\omega]/\cK\cong \cA(\Grpd,\tilde{\omega})$.
\end{lemma}

Now, we will specialize to the case of the inverse semigroup $\cS(G,E,\phi)$: let $E$ be a graph. Set
\[
\fX(E)=\{\alpha\mid \text{ infinite path in } E\}\cup\{\alpha\in\cP(E)\mid r(\alpha)\in\sing(E)\}.
\]
For $\alpha\in\cP(E)$, let
\[
\fX(E)\supset \cZ(\alpha)=\{x\in \fX(E) \mid \ x=\alpha y,\  y\in\fX(E)\}=\alpha\fX(E).
\]
For each finite set $F\subset \alpha\cP(E)$, let 
\[
\cZ(\alpha\setminus F)=\cZ(\alpha)\cap\left(\bigcup_{\beta\in F}\cZ(\alpha\beta)\right)^c.
\]
The sets $\cZ(\alpha\setminus F)$ are a basis of compact open sets for a locally compact, Hausdorff topology on $\fX(E)$ \cite[Theorem 2.1]{web}. We regard the latter as a topological space equipped with this topology; it is a Boolean space.  
Now assume that an Exel-Pardo tuple $(G,E,\phi)$ is given. Then there is an action of $\cS(G,E,\phi)$ on $\fX(E)$ such that $\dom(\alpha,g,\beta)=\cZ(\beta)$ and 
\[
(\alpha,g,\beta)(\beta x)=\alpha g(x)\in\cZ(\alpha)
\]
for all $x\in\beta\fX(E)$. Let $\Grpd(G,E,\phi)$ be the groupoid of germs associated to this action; it is an ample groupoid. If, moreover, $c:G\times E^1\to\cU(\ell)$ is a $1$-cocycle, and $\omega:\cS(G,E,\phi)\times\cS(G,E,\phi)\to\cU(\ell)_\bu$ is the semigroup $2$-cocycle of \eqref{map:omega}, then by \eqref{map:tildomega} we also have a twist 
\begin{equation}\label{map:tildomega2}
\tilde{\omega}:\Grpd(G,E,\phi)\times\Grpd(G,E,\phi)\to \cU(\ell).
\end{equation}

\begin{proposition}[{\cite[Proposition 4.2.2]{Willie}}]\label{prop:lgec=stein} Let $(G,E,\phi_c)$ be a twisted EP tuple and let \eqref{map:tildomega2} be the associated groupoid cocycle. Let $L(G,E,\phi_c)$ and $\cA(\Grpd(G,E,\phi),\tilde{\omega})$ be the twisted Exel-Pardo and Steinberg algebras. Then there is an algebra isomorphism 
\[
\begin{array}{cccc}
\Phi \colon & L(G,E,\phi_c) & \rightarrow  &  \cA(\Grpd(G,E,\phi),\tilde{\omega}) \\
 &\alpha g\beta^*  & \mapsto  & \chi_{[(\alpha,g,\beta),\cZ(\beta)]}.
\end{array}
\]
\end{proposition}

Some results about the simplicity of these algebras are given; for example:

\begin{theorem}[{\cite[Theorem 4.3.10]{Willie}}]\label{thm:spistein}
Let $(G,E,\phi_c)$ be a twisted EP-tuple with $G$ countable and $E$ countable and regular and let $\ell\supset\Q$ be a field.
If $C^*(G,E,\phi)$ is simple, then $L(G,E,\phi_c)$ is simple. Moreover, if $L(E)$ is simple purely infinite, then so is $L(G,E,\phi_c)$.
\end{theorem}

Also, a short view of twisted Katsura algebras is given. Let $E$ be a row-finite graph, $A=A_E\in \N_0^{(\reg(E)\times E^0)}$ its reduced incidence matrix, and $B\in \Z^{(\reg(E)\times E^0)}$ such that 
\begin{equation}\label{katcond0}
A_{v,w}=0\Rightarrow B_{v,w}=0.
\end{equation}
Under the obvious identification $\Z/\Z A_{v,w}\iso \N_{v,w}:=\{0,\dots,A_{v,w}-1\}$, translation by $B_{v,w}$ defines a bijection 
\[
\tau_{v,w}:\N_{v,w}\to \N_{v,w}, \ \ n\mapsto \overline{B_{v,w}+n}
\]
and therefore a $\Z$-action on $\N_{v,w}$; here $\bar{m}$ is the remainder of $m$ under division by $A_{v,w}$. A $1$-cocyle $\psi_{v,w}:\Z\times \N_{v,w}\to\Z$ for this action is determined by
\[
\psi_{v,w}(1,n)=\frac{B_{v,w}+n-\tau_{v,w}(n)}{A_{v,w}}.
\]
For every pair of vertices $(v,w)$ with $A_{v,w}\ne 0$, choose a bijection 
\[
\n_{v,w}:vE^1w\iso \{0,\dots, A_{v,w}-1\},
\]
and let
\[
\n=\coprod_{v,w}\n_{v,w}:E^1=\coprod_{v,w}vE^1w\to \N_0. 
\]
Consider the $\Z$-action on $E$ that fixes the vertices and is induced by the bijection $\sigma:E^1\to E^1$ that is determined by
\[
\n(\sigma(e))=\tau_{s(e),r(e)}(\n(e)).
\] 
Write $\Z=\langle t\rangle$ multiplicatively, and let $\phi:\Z\times E^1\to \Z$ be the $1$-cocycle determined by
\[
\phi(t,e)=t^{\psi_{s(e),r(e)}(1,\n(e))}.
\]
By \cite[Proposition 1.13]{HPSS}, the algebra $L(\Z,E,\phi)$ associated with the EP-tuple $(\Z,E,\phi)$ is the \emph{Katsura algebra} $\Oo_{A,B}=\Oo_{A,B}(\ell)$, itself the algebraic analogue of the $C^*$-algebra introduced by Katsura in \cite{Kat1}. Next, we consider a twisted version. For this purpose, we start with a row-finite matrix
\[
C\in \cU(\ell)^{(\reg(E)\times E^0)}
\]
such that 
\begin{equation}\label{katcond01}
A_{v,w}=0\Rightarrow C_{v,w}=1.
\end{equation}
 
Consider the $1$-cocycle $c:\Z\times E^1\to \cU(\ell)$ defined by
\begin{equation}\label{map:katwist}
c(t,e)=\left\{
\begin{array}{cl}
(-1)^{(A_{s(e),r(e)}-1)B_{s(e),r(e)}}C_{s(e),r(e)} & \text{ if } \n(e)=0,\\ 
1 &\text{ else}.
\end{array}
\right.    
\end{equation}

Write $\Oo_{A,B}^C=L(\Z,E, \phi_c)$ for the twisted Exel-Pardo algebra associated to the twisted EP-tuple $(\Z,E,\phi_c)$. 

We say that $(A,B)$ is \emph{KSPI} if $E$ is countable, \eqref{katcond0} holds, and if in addition we have:
\begin{itemize}
\item  For any pair $(v,w)\in E^0\times E^0$ there exists $\alpha\in\cP(E)$ such that $s(\alpha)=v$ and $r(\alpha)=w$. 
\item  For every vertex $v$, there are at least two distinct loops based at $v$.
\item $B_{v,v}=1$ for all $v\in E^0$. 
\end{itemize}

Katsura proved in \cite[Proposition 2.10]{Kat1} that if $(A,B)$ is a Katsura pair, then the $C^*$-algebra completion $C^*_{A,B}=\overline{\Oo_{A,B}(\C)}$ is simple and purely infinite. 

\begin{theorem}[{\cite[Theorem 5.4]{Willie}}]\label{thm:katpis}
Let $\ell\supset\Q$ be a subfield. Let $E$ be a countable regular graph, $A=A_E$ its incidence matrix, and $B\in \Z^{(E^0\times E^0)}$ such that $(A,B)$ is a KSPI pair. Let $C\in (\cU(\ell))^{(E^0\times E^0)}$ be arbitrary. Then the twisted Katsura algebra $\Oo_{A,B}^C$  is simple purely infinite. 
\end{theorem}

For a general field $\ell$, we have the following more restrictive simplicity criterion.

\begin{proposition}[{\cite[Proposition 5.5]{Willie}}]\label{prop:katsimp}
Let $E$, $A$ and $B$ be as in Theorem \ref{thm:katpis}, $\ell$ a field and $C\in\cU(\ell)^{(E^0\times E^0)}$ satisfying \eqref{katcond01}. In addition, assume that, whenever $A_{v,w}\ne 0=B_{v,w}$, for each $l\ge 1$ there exist finitely many paths $\alpha=e_1\dots e_n\in w\cP(E)$ such that $l\prod_{i=1}^{n-1}B_{s(e_i),r(e_i)}/A_{s(e_i),r(e_i)}\in\Z$. Then $\Oo_{A,B}^C$ is simple.
\end{proposition}

Now, we are ready to present a major generalization of self-similarity for graphs involving almost all the previous ideas in this section.

\section{Left cancellative small categories and their algebras}\label{Sect:LCSC}

Spielberg \cite{S1} extended the Kumjian-Pask idea of dealing with $k$-graphs \emph{as categories} to associate a $C^*$-algebra with a more general categorical object. To this end, he introduced \emph{categories of paths} --i.e. cancellative small categories with no (nontrivial) inverses-- as a generalization of higher rank graphs. The next definitions and results can be found in \cite[Sections 2 \& 3]{S1} and \cite[Section 2]{S2}

\subsection{Basic definitions} Given a small category $\Cat$, we denote by $\Obj$ its objects, and we identify $\Obj$ with identity morphisms, so that $\Obj\subseteq \Cat$. Given $\alpha\in\Cat$, we denote by $s(\alpha):=\dom(\alpha)\in \Obj$ and $r(\alpha):=\ran(\alpha)\in \Obj$. The right invertible elements of $\Cat$ are 
\[\Cat^{-1}:=\{\alpha\in\Cat\mid \exists \beta\in\Cat\text{ such that }\alpha\beta=s(\beta)\}\,.\]
\begin{definition}
	Given a small category $\Cat$, and let $\alpha,\beta,\gamma\in\Cat$:
	\begin{enumerate}
		\item $\Cat$ is \emph{left cancellative} if $\alpha\beta=\alpha\gamma$ then $\beta=\gamma$,
		\item $\Cat$ is \emph{right cancellative} if $\beta\alpha=\gamma\alpha$ then $\beta=\gamma$,
		\item $\Cat$ \emph{has no inverses} if $\alpha\beta=s(\beta)$ then $\alpha=\beta=s(\beta)$.
	\end{enumerate}

A \emph{category of paths} is a small category that is right and left cancellative and has no inverses. 
\end{definition}

Examples of left cancellative small categories (LCSC for short) are: 
\begin{enumerate}
\item Higher-rank graphs, and arbitrary subcategories of them.
\item The positive cone in a discrete ordered group (and, in particular, quasi-
lattice ordered groups).
\item $P$-graphs.
\end{enumerate}

Notice that if $\Cat$ is either left or right cancellative, then the only idempotents in $\Cat$ are $\Obj$. 

\begin{definition}[{c.f. \cite[Definition 1.2]{OP1}}]\label{defi1_1_1_2}
	Let $\Cat$ be a small category. Given $\alpha,\beta\in\Cat$, we say that $\beta$ \emph{extends $\alpha$} (equivalently \emph{$\alpha$ is an initial segments of $\beta$}) if there exists $\gamma\in\Lambda$ such that $\beta=\alpha\gamma$. We denote by $[\beta]=\{\alpha\in\Cat:\alpha\text{ is an initial segment of }\beta\}$. We write $\alpha\leq \beta$ if $\alpha\in[\beta]$.
\end{definition}

\begin{lemma}[{c.f. \cite[Lemma 1.3]{OP1}}]\label{lemma1_1_1_3}
Let $\Cat$ be a small category. Then
\begin{enumerate}
	\item the relation $\leq$ is reflexive and transitive,
	\item if $\Cat$ is left cancellative with no inverses, then $\leq$ is a partial order. 
\end{enumerate}
\end{lemma}

\begin{lemma}[{\cite[Lemma 2.3]{S2}}]\label{lemma1_1_1_4} 
	Let $\Cat$ be a LCSC, and let $\alpha,\beta\in\Cat$. Then, $\alpha\leq \beta$ and $\beta\leq \alpha$ if and only if $\beta\in\alpha\Cat^{-1}=\{\alpha\gamma:\gamma\in \Cat^{-1}\text{ with }r(\gamma)=s(\alpha)\}$.
\end{lemma}

We denote by $\alpha\approx\beta $ if $\beta \in\alpha\Cat^{-1}$. This is an equivalence relation. 
\begin{lemma}[{\cite[Lemma 2.5(ii)]{S2}}]\label{lemma1_1_1_5} 
	Let $\Cat$ be a LCSC, and let $\alpha,\beta\in\Cat$. Then the following statements are equivalent:\begin{enumerate}
		\item $\alpha\approx \beta$,
		\item $\alpha\Cat=\beta\Cat$,
		\item $[\alpha]=[\beta]$.
	\end{enumerate}
\end{lemma} 

\begin{notation}
	Let $\Cat$ be a LCSC. Given $\alpha,\beta\in\Cat$, we say :
	\begin{enumerate}
		\item $\alpha\Cap \beta$ if and only if $\alpha\Cat\cap\beta \Cat\neq\emptyset$,
		\item $\alpha\perp \beta$ if and only if $\alpha\Cat\cap\beta \Cat=\emptyset$.
	\end{enumerate}
\end{notation}

\begin{definition}
	Let $\Cat$ be a LCSC, and let $F\subset \Cat$. The elements of $\bigcap_{\gamma\in F}\gamma\Cat$ are \emph{the common extensions of $F$}. A common extension $\varepsilon$ of $F$  is \emph{minimal} if for any common extension $\gamma$ with $\varepsilon\in \gamma\Cat$ we have that $\gamma\in\varepsilon\Cat^{-1}=\{\varepsilon\mu \mid \mu\in \Cat^{-1}\text{ with }r(\mu)=s(\varepsilon)\}$.
\end{definition}

When $\Cat$ has no inverses, given $F\subseteq \Cat$ and given any minimal common extension $\varepsilon$ of $F$, if $\gamma$ is a common extension of $F$ with $\varepsilon\in \gamma\Cat$ then $\gamma= \varepsilon$. We denote by $\alpha\vee \beta$ the set of minimal extensions of $\alpha$ and $\beta$. Notice that if $\alpha\vee \beta\neq \emptyset$, then $\alpha\Cap \beta$, but the converse fails in general. 

\begin{definition}
	A LCSC $\Lambda$ is \emph{finitely aligned} if for every $\alpha,\beta\in\Lambda$ there exists a finite subset $\Gamma\subset \Lambda$ such that $\alpha\Lambda\cap\beta\Lambda=\bigcup_{\gamma\in \Gamma}\gamma\Lambda$.
\end{definition}

When $\Lambda$ is a finitely aligned LCSC, we can always assume that $\alpha\vee\beta=\Gamma$, where $\Gamma$ is a finite set of minimal common extensions of $\alpha$ and $\beta$.

\begin{definition}\label{definition2_3_7}
Let $\Lambda$ be a LCSC and $\alpha\in \Lambda$. A subset $F\subset r(\alpha)\Lambda$ is \emph{exhaustive with respect to $\alpha$} if for every $\gamma\in \alpha\Lambda$ there exists a $\beta\in F$ with $\beta\Cap\gamma$. We denote by $\mathsf{FE}(\alpha)$ the collection of finite sets of $r(\alpha)\Lambda$ that are exhaustive with respect to $\alpha$.
\end{definition}

Now, given a LCSC $\Cat$, we are ready to define the (Spielberg) algebra associated to $\Cat$ by generators and relations.

\begin{definition}[{\cite[Theorem 6.3]{S1}, \cite[Theorem 9.4]{S2}}]\label{Def:SpielbergAlgebra}
Given $\Lambda$ a (countable) finitely aligned LCSC, we define:
\begin{enumerate}
\item $C^*(\Lambda)$ to be the universal $C^*$-algebra generated by a family $\{T_\alpha \mid \alpha \in \Lambda\}$ satisfying:
\begin{enumerate}
\item $T_\alpha^*T_\alpha=T_{s(\alpha)}$.
\item $T_\alpha T_\beta=T_{\alpha \beta}$ if $s(\alpha)=r(\beta)$.
\item $T_\alpha T_\alpha^* T_\beta T_\beta^*=\bigvee_{\gamma\in \alpha \vee \beta} T_\gamma T_\gamma^*$.
\item $T_v=\bigvee_{\alpha \in F} T_\alpha T_\alpha^*$ for every $v\in \Lambda^0$ and for all $F\subset v\Lambda$ finite exhaustive set.
\end{enumerate}
\item Given any unital commutative ring $R$, we define $R\Lambda$ to be the $R$-algebra generated by a family $\{T_\alpha \mid \alpha \in \Lambda\}$ satisfying:
\begin{enumerate}
\item $T_\alpha^*T_\alpha=T_{s(\alpha)}$.
\item $T_\alpha T_\beta=T_{\alpha \beta}$ if $s(\alpha)=r(\beta)$.
\item $T_\alpha T_\alpha^* T_\beta T_\beta^*=\bigvee_{\gamma\in \alpha \vee \beta} T_\gamma T_\gamma^*$.
\item $T_v=\bigvee_{\alpha \in F} T_\alpha T_\alpha^*$ for every $v\in \Lambda^0$ and for all $F\subset v\Lambda$ finite exhaustive set.
\end{enumerate}
\end{enumerate}
\end{definition}

Spielberg presented a natural groupoid $\Grpd \vert_{\partial\Cat}$ associated to $\Cat$ and showed that, under mild hypotheses, the Spielberg algebra of $\Cat$ is isomorphic to the (reduced) groupoid algebra of $\Grpd \vert_{\partial\Cat}$. Here, we give a different picture of this groupoid, as constructed in \cite{OP1, OP2}.

\subsection{The inverse semigroup of a LCSC}

We start by recalling some basic concepts that we will need in the sequel.

\begin{definition}
Let $\Semi$ be an inverse semigroup, and let $\Idem$be its semi-lattice of idempotents. We say that $s,t\in \Semi$ are compatible, denoted by $s\sim t$, if both $s^*t$ and $st^*$ belong to $\Idem$.
\end{definition} 

\begin{lemma}[{\cite[Lemma 1.4.16]{L}}]\label{lemma1_1_2_5} Let $\Semi$ be an inverse semigroup, and let $\Sigma\subseteq \Semi$. If $\bigvee_{\alpha\in \Sigma}\alpha\in\Semi$ (the least-upper bound with the above defined partial order), then the elements of $\Sigma$ are pairwise compatible.
\end{lemma}

We say that an inverse semigroup $\Semi$ is \emph{(finitely) distributive} if whenever $\Sigma$ is a (finite) subset of $\Semi$ and $s\in\Semi$, if $\bigvee_{\alpha\in \Sigma}\alpha\in\Semi$ then $\bigvee_{\alpha\in \Sigma} s\alpha\in\Semi$ and $s\left(\bigvee_{\alpha\in \Sigma} \alpha\right)=\bigvee_{\alpha\in \Sigma}s\alpha$.

we say that $\Semi$ is \emph{(finitely) complete} if for every (finite) subset $\Sigma\subseteq \Semi$ of pairwise compatible elements we have that $\bigvee_{\alpha\in \Sigma}\alpha\in\Semi$.

Notice that this property is not necessarily inherited by its inverse subsemigroups. Indeed, the point is that given $\Sigma\subseteq \Semi$ a set of pairwise compatible elements, and $s\in\Semi$:
\begin{enumerate}
	\item Not necessarily $\bigvee_{\alpha\in \Sigma}\alpha \in\Semi$,
	\item even if $\bigvee_{\alpha\in \Sigma} \alpha\in\Semi$, it can happen that $\bigvee_{\alpha\in \Sigma}s\alpha\notin \Semi$.
\end{enumerate}\vspace{.3truecm}

Now, given a LCSC $\Lambda$, we define an inverse semigroup associated to $\Lambda$. We denote by $\cI(\Cat)$ the inverse semigroup of partial bijections of $\Cat$.

\begin{definition}
	Let $\Lambda$ be a LCSC. For any $\alpha\in\Lambda$, we define two maps:\begin{enumerate}
		\item  $\osh^\alpha:\alpha\Lambda\to s(\alpha)\Lambda$ given by $\alpha\beta \mapsto \beta\,,$
		\item  $\tau^\alpha:s(\alpha)\Lambda\to \alpha\Lambda$ given by $\beta \mapsto \alpha\beta\,.$
	\end{enumerate}
\end{definition}
Clearly $\osh^\alpha$ is injective, and since $\Lambda$ is left cancellative so is $\tau^\alpha$, whence both maps belong to $\cI(\Cat)$. Moreover, 
\[\osh^\alpha=\osh^\alpha\tau^\alpha\osh^\alpha\qquad\text{and}\qquad \tau^\alpha=\tau^\alpha\osh^\alpha\tau^\alpha\,,\]
for every $\alpha\in\Lambda$.

\begin{lemma}[{\cite[Lemma 2.3]{OP1}}]\label{lemma1_2_3} Given a LCSC $\Lambda$, we define the semigroup
	\[\Semi_\Lambda:=\left\langle \osh^\alpha,\tau^\alpha:\alpha\in \Lambda\right\rangle \subseteq \cI(\Cat).\]
Then $\Semi_\Lambda$ is an inverse semigroup.
\end{lemma}

Now, we connect $\Semi_\Lambda$ to the semigroups appearing in \cite{DM,S1,S2}.

\begin{definition}
	Let $\Lambda$ be a small category. A \emph{zigzag} is an even tuple of the form 
	\[\xi=(\alpha_1,\beta_1,\alpha_2,\beta_2,\ldots,\alpha_n,\beta_n)\]
	with $\alpha_i,\beta_i\in\Lambda$, $r(\alpha_i)=r(\beta_i)$ for every $1\leq i\leq n$ and $s(\alpha_{i+1})=s(\beta_{i})$ of every $1\leq i<n$. We denote by $\Zig$ the set of zigzags of $\Lambda$. Given $\xi\in\Zig$, we define $s(\xi)=s(\beta_n)$, $r(\xi)=s(\alpha_1)$ and $\bar{\xi}=(\beta_n,\alpha_n,\ldots,\beta_1,\alpha_1)$.
\end{definition}

Every $\xi\in \Zig$ defines a zigzag map $\varphi_\xi\in \cI(\Lambda)$ by
\[\varphi_\xi=\osh^{\alpha_1}\tau^\beta_1\cdots \osh^{\alpha_n}\tau^{\beta_n}\,.\]
We denote $\ZigM=\{\varphi_\xi:\xi\in\Zig \}$.  

\begin{remark}\label{rema1_3_2} $\mbox{ }$

	\begin{enumerate}
		\item For every $\alpha\in\Lambda$ we can define $\xi_\alpha:=(r(\alpha),\alpha)$. Notice that $\varphi_{\xi_\alpha}=\tau^\alpha$ and $\varphi_{\bar{\xi}_\alpha}=\osh^\alpha$.
		\item $\ZigM$ is closed by concatenation, and $\varphi_{\xi_1}\circ\varphi_{\xi_2}=\varphi_{\xi_1\xi_2}$.
		\item For every $\xi\in \Zig$, then $\varphi_{\bar{\xi}}=\varphi_\xi^{-1}$. 
		\end{enumerate}
\end{remark}

Thus,
\begin{lemma}[{\cite[Section 7.2]{BKQS}}]\label{zig_zag}
	If $\Lambda$ is a LCSC, then $\ZigM$ is an inverse semigroup. Moreover, $\ZigM=\Semi_\Lambda$.
\end{lemma}

To understand the structure of $\Semi_\Lambda$, we need to restrict ourselves to the case of finite aligned LCSC. First, we introduce a definition.

\begin{definition}\label{Def:Cond U}
Let $\Lambda$ be a finitely aligned LCSC, and let $s\in \Semi_\Lambda$. We say that a presentation $s=\bigvee_{i=1}^n\tau^{\alpha_i}\sigma^{\beta_i}$ is \emph{irredundant} if for all $1\leq i\ne j\leq n$ we have $\alpha_i\not\in [\alpha_j]$ and $\beta_i\not\in [\beta_j]$.
\end{definition}

\begin{lemma}[{\cite[Lemma 3.3 \& Theorem 6.3]{S1}}]\label{lemma1_2_4}
If $\Lambda$ is a finitely aligned LCSC, then every $f\in \Semi_\Lambda$ is the supremum of a finite family of elements of the form $\elmap{\alpha}{\beta}$ with $\alpha,\beta\in \Lambda$ and $s(\alpha)=s(\beta)$. Moreover, if such a decomposition is irredundant, then it is unique (up to permutation).
\end{lemma}

Notice that, given any finite family $\{\alpha_1,\ldots,\alpha_n\}\subset \Lambda$, the elements $\{\elmap{\alpha_i}{\alpha_i}\}_{i=1}^n\subset \Semi_\Lambda$ are pairwise compatible, so that $\bigvee_{i=1}^n\elmap{\alpha_i}{\alpha_i}\in \cI(\Lambda)$, but not necessarily to $\Semi_\Lambda$. Thus, in order to do some essential arguments, we need to consider a new object. 

\begin{definition}
	Let $\Lambda$ be a finitely aligned LCSC. We define 
	\[\SemiT_\Lambda=\left\lbrace\bigvee_{i=1}^n \elmap{\alpha_i}{\beta_i} \mid \{\elmap{\alpha_i}{\beta_i}\}_{i=1}^n\subset \Semi_\Lambda\text{ are pairwise compatible}\right\rbrace .\]
\end{definition}

By Lemma \ref{lemma1_2_4}, $\Semi_\Lambda\subseteq \SemiT_\Lambda\subset \cI(\Lambda)$. Moreover, by \cite[Proposition 1.4.20 \& Proposition 1.4.17]{L}, $\SemiT_\Lambda$ is closed by composition and inverses, and moreover, is finitely distributive. Thus,

\begin{lemma}[{\cite[Lemma 2.8]{OP1}}]
	Let $\Lambda$ be a finitely aligned LCSC. Then, $\SemiT_\Lambda$ is the smallest finitely complete, finitely distributive, inverse semigroup containing $\Semi_\Lambda$.  
\end{lemma}

Now, we proceed to understand who are the elements in $\Idem[\SemiT_\Lambda]$ and the order relation, in order to be able to work with filters on $\Idem$.

\begin{lemma}[{\cite[Lemma 2.9]{OP1}}]\label{lemma1_2_6}
	Let $\Lambda$ be a finitely aligned LCSC. Then $e\in \Idem[\SemiT_\Lambda]$ if and only if $e=\bigvee_{i=1}^n\elmap{\alpha_i}{\alpha_i}$ for some $\alpha_1,\ldots,\alpha_n\in\Lambda$.
\end{lemma} 

\begin{proposition}[{\cite[Proposition 2.10]{OP1}}]\label{propo1_2_9}
	Let $\Lambda$ be a finitely aligned LCSC, and let $e=\bigvee_{i=1}^n\elmap{\alpha_i}{\alpha_i}$, $f=\bigvee_{j=1}^m\elmap{\beta_j}{\beta_j}$ be idempotents of either $\Semi_\Lambda$ or $\SemiT_\Lambda$ (written in irredundant form). Then, the following statements are equivalent:
	\begin{enumerate}
		\item $e\leq f$,
		\item for each $1\leq k\leq n$, there exists $1\leq l\leq m$ such that $\beta_l\leq \alpha_k$.
	\end{enumerate}
\end{proposition}

\subsection{A space of filters}
Now, we describe the spaces of filters on $\widehat{\mathcal{E}}_0$ we need to define our groupoid. But, instead of using the classical description, we produce a new model, associated to the order properties of our category $\Cat$. This model should behave in some way as paths in a graph.\vspace{.2truecm}

To guarantee that every filter has such a path model, we introduce a restriction on $\Lambda$. 

\begin{definition}
	Let $\Lambda$ be a finitely aligned LCSC. We say that a filter $\eta \in \widehat{\mathcal{E}}_0$ enjoys condition $(*)$ if given $\bigvee_{i=1}^n\elmap{\alpha_i}{\alpha_i}\in\eta$, then there exists $1\leq j\leq n$ such that $\elmap{\alpha_j}{\alpha_j}\in\eta$\footnote{This property corresponds to the filter being \emph{prime}.}.
\end{definition} 

Notice that, if $\Lambda$ is singly aligned, then every filter enjoy condition $(*)$ (see e.g. \cite[Proposition 3.5]{DGKMW}). Before showing how to construct the path model of $\eta$, let us remark that there exist filters where this property always holds. 

\begin{lemma}[{\cite[Lemma 3.10 \& Corollary 3.11]{OP1}}]\label{lemma2_2_2} 
Let $\Lambda$ be a finitely aligned LCSC. Then, any $\eta\in\FiltT$ satisfies condition $(*)$. In particular, every $\eta\in\UlFilt$ satisfies condition $(*)$.
\end{lemma}

Now, we proceed to introduce a set of paths, inspired in a idea of Spielberg \cite{S1}.

\begin{definition}
	Let $\Lambda$ be a finitely aligned LCSC. A nonempty subset $F$ of $\Lambda$ is:\begin{enumerate}
		\item \emph{Hereditary}, if $\alpha\in\Lambda$, $\beta\in F$ and $\alpha\leq \beta$ implies $\alpha\in F$,
		\item \emph{(upwards) directed},  $\alpha,\beta\in F$ implies that there exists $\gamma\in F$ with $\alpha,\beta\leq \gamma$.
	\end{enumerate}
We denote by $\Lambda^*$ the set of nonempty, hereditary, directed subsets of $\Lambda$.
 \end{definition}
 
The elements of $\Cat^*$ are our paths. We remark here that, in the context of partially ordered abelian groups (with interpolation), this kind of sets are called \emph{intervals}. However, we prefer to use the term \emph{path} to refer to them in our context.

\begin{definition}
	 Given $\eta\in\Filt$, we define 
	 \[\Delta_\eta:=\{\alpha\in\Lambda \mid \elmap{\alpha}{\alpha}\in\eta\} .\]
\end{definition} 

\begin{lemma}[{\cite[Lemma 3.14]{OP1}}]\label{lemma2_2_6}
	Let $\Lambda$ be a finitely aligned LCSC. For every $\eta\in \Filt$  satisfying condition $(*)$ we have that $\Delta_\eta\in \Lambda^*$.
\end{lemma}

\begin{definition}\label{Def: Cond_Ast}
If $\Lambda$ is a finitely aligned LCSC, we define
\[\widehat{\mathcal{E}}_{*}=\{ \eta \in \Filt \mid \eta \text{ satisfies condition } (*)\}.\]
\end{definition}

Thus,

\begin{corollary}[{\cite[Corollary 3.16]{OP1}}]\label{corol2_2_7}
If $\Lambda$ is a finitely aligned LCSC, then 
\[\begin{array}{rl} \Phi:\widehat{\mathcal{E}}_{*} & \longrightarrow \Lambda^* \\
\eta & \longmapsto \Delta_\eta\,,
\end{array}\]
is a well-defined map.
\end{corollary}

Now, we construct an inverse for this map.

\begin{definition}
	Given $F\in\Lambda^*$, we define 
	\[\eta_F:=\{f\in\EGE\mid f\geq \elmap{\alpha}{\alpha}\text{ for some }\alpha\in F\} .\]
\end{definition}

\begin{lemma}[{\cite[Lemma 3.18]{OP1}}]\label{lemma2_2_9} If $\Lambda$ is a finitely aligned LCSC, then for every $F\in\Lambda^*$ we have that $\eta_F\in\widehat{\mathcal{E}}_{*}$.
\end{lemma}

\begin{corollary}[{\cite[Corollary 3.19]{OP1}}]\label{corol2_2_10}
	If $\Lambda$ is a finitely aligned LCSC, then 
	\[\begin{array}{rl} \Psi:\Lambda^* & \longrightarrow \widehat{\mathcal{E}}_{*} \\
	F & \longmapsto \eta_F\,,
	\end{array},\]
	is a well-defined map.
\end{corollary}

\begin{lemma}[{\cite[Lemma 3.20]{OP1}}]\label{lemma2_2_11}
If $\Lambda$ is a finitely aligned LCSC, then $\Phi$ and $\Psi$ are naturally inverse bijections. 
\end{lemma}

Before tracking $\UlFilt$ and $\FiltT$ through $\Psi$, we need to consider a suitable topology defined on $\Lambda^*$:

\begin{definition}\label{top_sets}
	Let $\Lambda$ be a finitely aligned LCSC. Then, given $X,Y\subset \Lambda$ finite sets, we define 
	\[\Mtop^{X,Y}=\{F\in\Lambda^*\mid X\subseteq F\text{ and }Y\cap F=\emptyset \} .\]
	We endow a topology on $\Lambda^*$, with a basis of open sets
	\[\{\Mtop^{X,Y}\mid X,Y\subseteq \Lambda\text{ finite sets}\} .\]
\end{definition}

Since $\widehat{\mathcal{E}}_{*}\subseteq \Filt$, we can equip $\widehat{\mathcal{E}}_{*}$ with the induced topology. To simplify, we also use $\mathcal{U}(X,Y)$ to denote the restriction of the basic open sets of the topology of $\Filt$ to $\widehat{\mathcal{E}}_{*}$.

\begin{lemma}[{\cite[Lemma 3.22]{OP1}}]\label{lemma2_3_2}
Let $\Lambda$ be a finitely aligned LCSC. Then,
\[\{\cU(X,Y)\mid X=\{\elmap{\alpha}{\alpha}\}\,,Y=\{\elmap{\beta_i}{\beta_i}\}_{i=1}^n \} ,\]
is a basis for the topology of $\widehat{\mathcal{E}}_{*}$.
\end{lemma}

Hence, both $\widehat{\mathcal{E}}_{\infty}$ and $\widehat{\mathcal{E}}_{\text{tight}}$ are topological subspaces of $\widehat{\mathcal{E}}_{*}$, in particular, the closure of $\widehat{\mathcal{E}}_{\infty}$ in $\Filt$ coincides with the closure of $\widehat{\mathcal{E}}_{\infty}$ in $\widehat{\mathcal{E}}_{*}$. As a consequence, we have the following: 

\begin{lemma}[{\cite[Lemma 3.23]{OP1}}]\label{lemma2_3_3} Let  $\Lambda$ be a finitely aligned LCSC. Then $\Phi$ and $\Psi$ are homeomorphisms.
\end{lemma}

Now, we identify both $\Phi(\UlFilt)$ and $\Phi(\FiltT)$ intrinsically. 

\begin{definition} Given $\Lambda$ a LCSC, we say that $C\in\Lambda^*$ is \emph{maximal} if whenever $C\subset D$ with $D\in \Lambda^*$ we have that $D=\Lambda$.  we denote $\Lambda^{**}:=\{C\in\Lambda^*\mid C\text{ is maximal}\}$.
\end{definition}

\begin{lemma}[{\cite[Lemma 3.26]{OP1}}]\label{lemma2_3_6} Let $\Lambda$ be a countable, finitely aligned LCSC. Then given $\eta\in \Filt$ the following statements are equivalent:
\begin{enumerate}
\item $\eta\in\UlFilt$,
\item $\Delta_\eta\in \Lambda^{**}$. 
\end{enumerate}
\end{lemma}

This means that $\Phi(\UlFilt)=\Lambda^{**}$. Since $\Phi$ is continuous and $\FiltT=\overline{\UlFilt}^{\|\cdot\|_{\widehat{\mathcal{E}}_{*}}}$ we have that
%\begin{align*}\label{image_tight} 
\[ \Cat_{\text{tight}}:= \Phi\left( \FiltT\right)  =\Phi\left( \overline{\UlFilt}^{\|\cdot\|_{\widehat{\mathcal{E}}_{*}}}\right)=\overline{\Phi(\UlFilt)}^{\|\cdot\|_{\Lambda^*}}= \overline{\Lambda^{**}}^{\|\cdot\|_{\Lambda^*}}.\]
%\end{align*}

\subsection{The action} Since we have a homeomorphism $\Psi:\widehat{\mathcal{E}}_{*}\to \Lambda^*$ with inverse $\Psi$ (Lemma \ref{lemma2_3_3}), we can transfer the action of $\Semi_\Lambda$ on $\widehat{\mathcal{E}}_{*}$ to $\Lambda^*$. First, we fix the domain and range.

\begin{definition}
	Let $s=\elmap{\alpha}{\beta}\in \Semi_\Lambda$. Then, we define 
	\[E_\alpha:=E_{ss^*}=\Phi(D_{ss^*})=\Phi(\cU(\{\alpha\},\emptyset))=\{C\in\Lambda^*\mid\alpha\in C \}\,\]
	and
\[E_\beta:=E_{s^*s}=\Phi(D_{s^*s})=\Phi(\cU(\{\beta\},\emptyset))=\{C\in\Lambda^*\mid \beta\in C \}\,.\]
	Given $s=\bigvee_{i=1}^n\elmap{\alpha_i}{\beta_i}$, we define 
	\[E_{s^*s}=\bigcup_{i=1}^n E_{\beta_i}=\bigcup_{i=1}^n\Phi(D_{\elmap{\beta_i}{\beta_i}})=\Phi(D_{s^*s})\,.\]
\end{definition}

The sets $E_{s^*s}$ and $E_{ss^*}$ are the natural candidates for being the domain and range of the partial action of $s\in \Semi_\Lambda$ on $\Lambda^*$. 

The next step is to define the action. We start by defining the action in the particular case of $s=\elmap{\alpha}{\beta}$. In this case, $F\in E_\beta$ if and only if $\beta\in F$. Then, we define 
\[\elmap{\alpha}{\beta}\cdot F=\bigcup_{\substack{\beta\leq \gamma\\ \gamma\in F}}[\alpha\sigma^\beta (\gamma)]\,,\]
where $[\delta]$ is the set of initial segments of $\delta\in\Lambda$, which clearly belong to $\Lambda^*$. 

Moreover, $\elmap{\alpha}{\beta}\cdot F\in E_\alpha$. Thus, in this case, we have that
\[\elmap{\alpha}{\beta}\cdot :E_\beta\to E_\alpha\,,\qquad F\mapsto \elmap{\alpha}{\beta}\cdot F\,,\]
is a well-defined map. 

Now, set $s=\bigvee_{i=1}^n\elmap{\alpha_i}{\beta_i}\in \Semi_\Lambda$, with 
\[E_{s^*s}=\bigcup_{i=1}^n E_{\beta_i}\qquad\text{and}\qquad E_{ss^*}=\bigcup_{i=1}^n E_{\alpha_i}\,.\]
By Lemma \ref{lemma1_1_2_5}, $\elmap{\alpha_i}{\beta_i}$ and $\elmap{\alpha_j}{\beta_j}$ are compatible for $1\leq i,j\leq n$, and then $\elmap{\alpha_i}{\beta_i}\elmap{\beta_j}{\alpha_j}$ and $\elmap{\beta_i}{\alpha_i}\elmap{\alpha_j}{\beta_j}$ are idempotents in $\Semi_\Lambda$.

For every $1\leq i\leq n$  and every $F\in E_{\beta_i}$, we define 
\[s\cdot F=\elmap{\alpha_i}{\beta_i}\cdot F\,.\]
The point is that, if $F\in E_{\beta_i}\cap E_{\beta_j}$ for $1\leq i\ne j\leq n$, then $\elmap{\alpha_i}{\beta_i}\cdot F=\elmap{\alpha_j}{\beta_j}\cdot F$. 

Because of this fact, we define the map
\[s\cdot :E_{s^*s}\to E_{ss^*}\,,\]
as follows: given $F\in E_{s^*s}=\bigcup_{i=1}^n E_{\beta_i}$ we can assume, after re-indexing, that $\beta_1,\ldots,\beta_k\in F$ and $\beta_{k+1},\ldots,\beta_{n}\notin F$. Thus,
\[s\cdot F=\bigcup_{\substack{\bigvee_{i=1}^k\beta_i\leq \gamma\\ \gamma \in F}}[\alpha_i\sigma^{\beta_i}(\gamma)]\,.\]

Now, we prove a result that is essential to fix the dictionary.

\begin{lemma}[{\cite[Lemma 4.6]{OP1}}]\label{lemma3_3_1}
Let $\Lambda$ be a finitely aligned LCSC. Then, for every $s\in \Semi_\Lambda$ and any $\eta\in D_{s^*s}\cap \widehat{\mathcal{E}_{*}}$ we have that $s\cdot \Delta_\eta=\Delta_{s\cdot \eta}$.
\end{lemma}

\begin{corollary}[{\cite[Corollary 4.7]{OP1}}]\label{corol3_3_2}
Let $\Lambda$ be a countable, finite aligned  LCSC. Then, given $s\in \Semi_\Lambda$:
\begin{enumerate}
\item $s\cdot$ restricts to an action on $\Lambda^{**}$,
\item $s\cdot$  restricts to an action on $\Lambda_{tight}$.
\end{enumerate}
\end{corollary}

\subsection{The tight groupoid} We are ready to give a nice description of $\Grpd_{tight}(\Semi_\Lambda)$.

\begin{lemma}[{\cite[Lemma 4.8]{OP1}}]\label{lemma_simpli}
Let $\Lambda$ be a finitely aligned LCSC, let $s=\bigvee_{i=1}^n\elmap{\alpha_i}{\beta_i}\in\Semi_{\Lambda}$, and let $\xi\in D_{s^*s}$. Suppose that $\beta_k\in \Delta_\xi$ for some $1\leq k\leq n$. Then $[\elmap{\alpha_k}{\beta_k},\xi]=[s,\xi]\in \Grpd_{tight}(\Semi_\Lambda)$.  In particular, 
\[\Grpd_{tight}(\Semi_\Lambda)=\{[\elmap{\alpha}{\beta},\xi]\mid s(\alpha)=\beta,\, \beta\in \Delta_\xi \}\,.\]
\end{lemma}

Thus, we have a picture for our groupoid:

\begin{lemma}[{\cite[Lemma 4.9]{OP1}}]\label{lemma3_4_1}
Let $\Lambda$ be a finite aligned LCSC. Then, $\Grpd_{tight}(\Semi_\Lambda)$ is topologically isomorphic to the germ groupoid $\Semi_\Lambda\rtimes \Lambda_{tight}$.
\end{lemma}

The next result is a refined version of \cite[Lemma 4.12]{S1}, which is basic for next results.

\begin{lemma}[{\cite[Lemma 5.1]{OP1}}]\label{lem_eq_sets}
	Let $\Lambda$ be a finitely aligned LCSC, let $F,G\in\Lambda^*$ and $\alpha,\beta\in\Lambda$ be such that $\tau^\alpha\cdot F=\tau^\beta\cdot G$. Then, there exist $\delta\in F$ and $\gamma\in G$ such that $\alpha\delta=\beta\gamma$.
\end{lemma}

Let us compare this with Spielberg's groupoid. First, we recall the definition of Spielberg's groupoid associated to a small category: we start by defining an equivalence relation on $\Lambda \times\Lambda \times \Lambda^*$ by saying that $(\alpha,\beta, F)\sim (\alpha',\beta',F')$ if there exist $G\in\Lambda^*$, $\gamma,\gamma'\in\Lambda$ such that $F=\tau^\gamma\cdot G$, $F'=\tau^{\gamma'}\cdot G$, $\alpha\gamma=\alpha'\gamma'$ and $\beta\gamma=\beta'\gamma'$. Denote $\mathcal{G}=\Lambda \times\Lambda \times \Lambda^*/\sim$. Now, we define a partial operation on $\mathcal{G}$. To this end, fix the set of composable pairs 
\[\Grpd^{(2)}:=\{ ([\alpha,\beta,F],[\gamma,\delta,G])\mid \tau^\beta\cdot F=\tau^\gamma\cdot G)  \}\,,\]
and define $[\alpha,\beta, F]^{-1}=[\beta,\gamma, F]$. Given a pair $([\alpha,\beta, F],[\gamma,\delta, G])\in\Grpd^{(2)}$, we define the multiplication by 
\[[\alpha,\beta, F][\gamma,\delta, G]=[\alpha \xi,\delta\eta, H]\,,\]
where $\xi\in F$ and $\eta\in G$ are the elements given in Lemma \ref{lem_eq_sets} such that $\beta \xi=\gamma \eta$, and $H=\sigma^\xi\cdot F=\sigma^\eta\cdot G$. Finally, the sets $[\alpha,\beta, U]:=\{[\alpha,\beta,F]:F\in U\}$ for $U$ an open subset of $\Lambda^*$ form a basis for the topology of $\Grpd$, under which $\Grpd$ is an \'etale groupoid. By Corollary \ref{corol3_3_2}, we have that $\Grpd_{|\partial \Lambda}=\{[\alpha,\beta,F]\in \Grpd\mid F\in \Lambda_{tight} \}$.

\begin{proposition}[{\cite[Proposition 5.2]{OP1}}]\label{iso_group_sp}
Let $\Lambda$ be a countable, finitely aligned LCSC. Then the map 
\[\begin{array}{cccc}
\Phi : & G_{\vert \partial\Lambda} & \rightarrow & \Semi_\Lambda\rtimes \Lambda_{\text{tight}}\\ & [\alpha,\beta, F] & \mapsto & [\elmap{\alpha}{\beta}, \tau^\beta\cdot F]
\end{array}
\]
is an isomorphism of topological groupoids.
\end{proposition}

With respect to the algebras associated to $\Cat$, we have the following results.

\begin{proposition}[{\cite[Proposition 4.15]{OP1}}]\label{Prop:tighness}
Let $\Lambda$ be a (countable) finitely aligned LCSC. Then:
\begin{enumerate}
\item The natural semigroup homomorphism
\[
\begin{array}{cccc}
\pi : & \Semi_\Lambda  & \rightarrow  &  C^*(\Lambda) \\
 & \elmap{\alpha}{\beta} &  \mapsto &   T_\alpha T_\beta^*
\end{array}
\]
is a universal tight representation of $\Semi_\Lambda$ in the category of $C^*$-algebras.
\item For any unital commutative ring $R$, the natural semigroup homomorphism
\[
\begin{array}{cccc}
\pi : & \Semi_\Lambda  & \rightarrow  &  R\Lambda \\
 & \elmap{\alpha}{\beta} &  \mapsto &   T_\alpha T_\beta^*
\end{array}
\]
is a universal tight representation of $\Semi_\Lambda$ in the category of $R$-algebras.
\end{enumerate}
\end{proposition}

Hence, because of \cite[Theorem 13.3]{Exel1}  and \cite[Corollary 5.3]{St2}, we have

\begin{theorem}[{\cite[Theorem 4.16]{OP1}}]\label{Th: GroupoidRep}
Let $\Lambda$ be a (countable) finitely aligned LCSC. Then:
\begin{enumerate}
\item $C^*(\Lambda)\cong C^*(\mathcal{G}_{\text{tight}}(\Semi_\Lambda))$.
\item For any unital commutative ring $R$, $R\Lambda\cong A_R(\mathcal{G}_{\text{tight}}(\Semi_\Lambda))$.
\end{enumerate}
\end{theorem}

Thus, given any LCSC $\Cat$ with FA, we can use the results in \cite{EP1} to describe the properties of the algebras in terms of combinatorial properties of $\Cat$, as we did with self-similar graphs. 

As a final remark, Spielberg $C^*$-algebras model all (reduced) groupoid $C^*$-algebras associated to second countable, ample groupoids. To see this, notice that:
\begin{enumerate}
\item Exel \cite{ExelDual} showed that every second countable, ample groupoid $\Grpd$ is topologically isomorphic to $\Grpd_{\text{tight}}(\Semi)$ for a suitable, countable inverse subsemigroup $\Semi \subset \Grpd$.
\item Donsig et al. \cite{DGKMW} showed that every (countable) inverse semigroup $\Semi$ is Morita equivalent to the inverse semigroup $\Semi_\Cat$ associated to a \emph{singly aligned} left cancellative small category $\Cat$. In particular, $\Grpd_{\text{tight}}(\Semi)$ and $\Grpd_{\text{tight}}(\Semi_\Cat)$ are strongly Morita equivalent, and then so are their associated algebras.
\end{enumerate}
This means that Spielberg's construction provides the most general combinatorial model possible for $C^*$-algebras.

\subsection{Simplicity} In \cite[Section 10]{S1} conditions are given in a category of paths $\Lambda$ for $\Grpd_{|\partial \Lambda}$ being topologically free, minimal and locally contractive, but right-cancellation of $\Lambda$ is crucial in the proofs therein. We use the isomorphism in Proposition \ref{iso_group_sp} and the characterization of these properties given in \cite{EP1}, to extend the Spielberg results in \cite[Section 10]{S1} to finitely aligned LCSC. 

\begin{definition}
	Let $\Semi$ be an inverse semigroup, and let $s\in \Semi$. Given an idempotent $e\in \EGE$ such that $e\leq s^*s$, we say that:
	\begin{enumerate}
		\item $e$ is \emph{fixed} under $s$, if $se=e$,
		\item $e$ is \emph{weakly-fixed} under $s$, if $(sfs^*) f\neq 0$, for every non-zero idempotent $f\leq e$.
	\end{enumerate}
\end{definition}

\begin{definition}
	Given an action $\alpha:\Semi\curvearrowright X$, let $s\in\Semi$, and let $x\in D_{s^*s}$.
	\begin{enumerate}
		\item $\alpha_s(x)=x$, we say that $x$ is a \emph{fixed point} for $s$. We denote by $F_s$ the set of fix points for $s$.
		\item If there exists $e\in \EGE$, such that $e\leq s$, and $x\in D_{e}$, we say that $x$ is a \emph{trivially fixed point} for $s$.
		\item We say that $\alpha$ is a \emph{topologically free} action, if for every $s$ in $\Semi$, the interior of the set of fixed points for $s$ consists of trivial fixed points.  
	\end{enumerate}	 
\end{definition} 

Given an action $\alpha:\Semi\curvearrowright X$, the groupoid $\Semi\rtimes X$ is effective if and only if the action $\alpha$ is topologically free \cite[Theorem 4.7]{EP1}. 

\begin{remark}\label{remark_top}
	Let $\Lambda$ be a finitely aligned LCSC, and let $\Semi_\Lambda\curvearrowright \FiltT$ be the associated action. Let $s\in\Semi_\Lambda$, and $\xi\in D_{s^*s}\cap \FiltT$. If $s=\bigvee_{i=1}^n\elmap{\alpha_i}{\beta_i}$, then $s^*s=\bigvee_{i=1}^n\elmap{\beta_i}{\beta_i}\in\xi$. Since $\xi$ satisfies condition $(*)$, there exists $1\leq j\leq n$ such that $\elmap{\beta_j}{\beta_j}\in\xi$. Let $C\in \Lambda_{\text{tight}}$ be such that $\xi=\eta_C$. Then we have that $\beta_j\in C$. By the definition of the action $\Semi_{\Lambda}\curvearrowright \Lambda_{\text{tight}}$ we have that $s\cdot C=\elmap{\alpha_j}{\beta_j}\cdot C$, and hence $s\cdot \xi=\elmap{\alpha_j}{\beta_j}\cdot \xi$. Thus, without loss of generality, we can assume that $s=\elmap{\alpha_j}{\beta_j}$.
\end{remark}

\begin{theorem}[{\cite[Theorem 6.4]{OP1}}]\label{eq_top_free}
	Let $\Lambda$ be a countable, finitely aligned LCSC. If either $\mathcal{G}_{\text{tight}}(\Semi_\Lambda)$ is Hausdorff or $\hat{\mathcal{E}}_\infty =\hat{\mathcal{E}}_{\text{tight}}$, then the following statements are equivalent:
	\begin{enumerate}
		\item $\Grpd_{tight}(\Semi_\Lambda)$ is effective.
		\item For every $s\in \Semi_\Lambda$, and for every $e\in \EGE_\Lambda$ which is weakly-fixed under $s$, there exists a finite cover for $e$ consisting of fixed idempotents. 
		\item Given $\alpha,\beta\in \Lambda$ with $r(\alpha)=r(\beta)$ and $s(\alpha)=s(\beta)$, if $\alpha\delta\Cap \beta\delta$ for every $\delta\in s(\alpha)\Lambda$ then there exists $F\in\mathsf{FE}(s(\alpha))$ such that $\alpha\gamma=\beta\gamma$ for every $\gamma\in F$.
	\end{enumerate}
\end{theorem}

\begin{remark}
Observe that if $\Lambda$ has right cancellation, condition $(3)$ in Theorem \ref{eq_top_free} reduces to aperiodicity as defined in \cite[Definition 10.8]{S1}.
\end{remark}
 
\begin{theorem}[{\cite[Theorem 6.6]{OP1}}]\label{eq_minimal}
If $\Lambda$ is a countable, finitely aligned LCSC, then the following statements are equivalent: 
\begin{enumerate}
	\item $\Grpd_{tight}(\Semi_\Lambda)$ is minimal.
	\item For every nonzero $e,f\in \EGE_\Lambda$, there are $s_1,\ldots,s_n\in \Semi_\Lambda$, such that $\{s_ifs_i^*\}_{i=1}^n$ is an outer cover for $e$.
	\item For every $\alpha,\beta\in \Lambda$ there exists $F\in \mathsf{FE}(\alpha)$ such that for each $\gamma\in F$, $s(\beta)\Lambda s(\gamma)\neq \emptyset$.
\end{enumerate}
\end{theorem}

Then, we have the following result:

\begin{theorem}[{\cite[Theorem 6.7]{OP1}}]\label{theorem_simpleLambda}
Let $\Lambda$ be a countable, finitely aligned LCSC. If $\mathcal{G}_{\text{tight}}(\Semi_\Lambda)$ is Hausdorff and amenable, then the following statements are equivalent: 
\begin{enumerate}
\item $C^*(\Lambda)$ is simple.
\item The following properties hold:
\begin{enumerate}
\item Given $\alpha,\beta\in \Lambda$ with $r(\alpha)=r(\beta)$ and $s(\alpha)=s(\beta)$, if $\alpha\delta\Cap \beta\delta$ for every $\delta\in s(\alpha)\Lambda$ then there exists $F\in\mathsf{FE}(s(\alpha))$ such that $\alpha\gamma=\beta\gamma$ for every $\gamma\in F$.
\item For every $\alpha,\beta\in \Lambda$ there exists $F\in \mathsf{FE}(\alpha)$ such that for each $\gamma\in F$, $s(\beta)\Lambda s(\gamma)\neq \emptyset$.
\end{enumerate}
\end{enumerate}
\end{theorem}

The same result holds for the $K$-algebra $K\Lambda$ over any field $K$, without requiring amenability for $\mathcal{G}_{\text{tight}}(\Semi_\Lambda)$.

\section{Zappa-Sz\'ep products of groups acting on LCSC categories}

B\'edos, Kaliszewski, Quigg, and Spielberg \cite{BKQS} extended the notion of self-similar action of a group on a graph to \emph{left cancellative small categories}. But instead of giving a generalization of the original definition, they chose an approach related to a more general construction, called the Zappa-Sz\'ep product, introduced by Brin \cite{Br}; Lawson used it to describe self-similar groups via Zappa-Sz\'ep products of groups on free monoids \cite{Lawson_SSI}. The relevant fact is that, in this particular case, the new construction turns out to be a left-cellative small category, too. Then, they defined an algebra associated to their self-similar object by using a Cuntz-Pimsner construction.

Ortega and Pardo \cite{OP1} studied the $C^*$-algebras associated to Zappa-Sz\'ep products of groups on left cancellative small categories, using the groupoid approach introduced in Section \ref{Sect:LCSC}. Subsequently, they extended the construction to Zappa-Sz\'ep products of \emph{groupoids} acting on left cancellative small categories \cite{OP2}, using ideas similar to those of \cite{LRRW18}; we will present this construction in the next section.\vspace{.2truecm}

In this section we analyze the notion of Zappa-Sz\'ep products of groups acting on left cancellative small categories. 

If $G$ is a discrete group (with unit $\ideG$), we use multiplicative notation for the group operation. We say that the group $G$ acts on a small category $\Lambda$ by permutations when 
\[r(g\cdot\alpha)=g\cdot r(\alpha) \text{ and } s(g\cdot \alpha)=g\cdot s(\alpha) \text{ for every }\alpha\in\Lambda,\,g\in G\,.\]
For the remainder of the section, we assume that $G$ acts by permutations on $\Lambda$.

A \emph{1-cocyle} for the action of $G$ on a small category $\Lambda$ is a function $\varphi:G\times \Lambda \to G$ that satisfies the 1-cocyle identity 
\[\varphi(gh,\alpha)=\varphi(g,h\cdot \alpha)\varphi(h,\alpha) \text{ for all }g,h\in G,\,\alpha\in\Lambda\,.\]
In particular, the 1-cocycle identity says that $\varphi(\ideG,\alpha)=\ideG$ for every $\alpha\in\Lambda$. 

\begin{definition}\label{definition_CategoryCocycle1}
 A 1-cocycle $\varphi:G\times \Lambda \to \Lambda$ on a small category $\Cat$ is a \emph{category cocycle} if for all $g\in G$, $v\in\Lambda^0$, and $\alpha,\beta\in\Lambda$ with $s(\alpha)=r(\beta)$ we have: 
	\begin{enumerate}
		\item $\varphi(g,v)=g$.
		\item $\varphi(g,\alpha)\cdot r(\alpha)=g\cdot r(\alpha)$.
		\item $g\cdot(\alpha\beta)=(g\cdot \alpha)(\varphi(g,\alpha)\cdot \beta)$.
		\item $\varphi(g,\alpha\beta)=\varphi(\varphi(g,\alpha),\beta)$.
	\end{enumerate}
We call $(\Lambda,G,\varphi)$ a \emph{category system}.
\end{definition}

The next definition comes from \cite{BKQS}:

\begin{definition}\label{definition_zappa-szed1}
	Let $(\Lambda,G,\varphi)$ be a category system. We denote by $\Lambda\rtimes^\varphi G$ the small category with 
	\[\Lambda\rtimes^\varphi G: = \Lambda\times G \text{ and } (\Lambda\rtimes^\varphi G)^0=\Lambda\times\{\ideG\}\,,\]
	and $r,s:\Lambda\rtimes^\varphi G\to (\Lambda\rtimes^\varphi G)^0$ maps defined by 
	\[r(\alpha,g)=(r(\alpha),\ideG) \text{ and }s(\alpha,g)=(g^{-1}\cdot s(\alpha),\ideG)\,.\]
	Moreover, for $(\alpha,g),(\beta,h)$ with $s(\alpha,g)=r(\beta,h)$, the composition is defined by
	\[(\alpha,g)(\beta,h)=(\alpha(g\cdot\beta),\varphi(g,\beta)h)\,.\]
	We call $\Lambda\rtimes^\varphi G$ the \emph{Zappa-Sz\'ep product of $(\Lambda,G,\varphi)$}.
\end{definition}

It was proven that $\Lambda\rtimes^\varphi G$ is left cancellative whenever so is $\Lambda$ \cite[Proposition 3.5]{BKQS}, and as observed in \cite[Remark 3.9]{BKQS}, the elements of the form $(v,g)$ (where $v\in \Lambda^0$ and $g\in G$) are units of $\Lambda\rtimes^\varphi G$. Then, given $(\alpha,g)\in\Lambda\rtimes^\varphi G$ and $h\in G$, we have that 
\[(\alpha,g)(g^{-1}\cdot s(\alpha), g^{-1}h)=(\alpha,h)\,,\]
so $(\alpha,g)\approx (\alpha,h)$. 
Moreover,  $\Lambda\rtimes^\varphi G$ is finitely aligned (singly aligned) whenever $\Lambda$ is finitely aligned (singly aligned) \cite[Proposition 3.12]{BKQS}. In particular, 
\[(\alpha,g)\vee (\beta,h)=(\alpha\vee \beta)\times\{\ideG\}\,.\] 

\begin{definition}
A category system $(\Lambda,G,\varphi)$ is called \emph{pseudo free} if, whenever $g\cdot \alpha=\alpha$ and $\varphi(g,\alpha)=\ideG$, then $g=\ideG$. 
\end{definition}

\begin{proposition}[{\cite[Proposition 5.6]{EP2}}]\label{prop_pseudo}
Let $\Lambda$ be a LCSC, and let $(\Lambda,G,\varphi)$ be a pseudo-free category system. Then, for all $g_1,g_2\in G$, and $\alpha\in\Lambda$, one has that 
	\[g_1\cdot \alpha=g_2\cdot \alpha\text{ and }\varphi(g_1,\alpha)=\varphi(g_2,\alpha)\Rightarrow g_1=g_2\,.\] 
\end{proposition}

\begin{remark}
	Given a category system $(\Lambda,G,\varphi)$ on $\Lambda$ LCSC, where $\Lambda$ is a right cancellative, it may happen that $\Lambda\rtimes^\varphi G$ does not satisfy right cancelation. Given $(\alpha,a),(\beta,b)$ and $(\gamma,g)$ in $\Lambda\rtimes^\varphi G$ we have that
	$(\alpha,a)(\gamma,g)=(\beta,b)(\gamma,g)$ if and only if $\alpha(a\cdot \gamma)=\beta(b\cdot \gamma)$ and $\varphi(a,\gamma)=\varphi(b,\gamma)$. In particular, the system is pseudo-free if and only if $\Lambda\rtimes^\varphi G$ is right cancellative.  
\end{remark}

\begin{remark} Let $\Lambda$ be a finitely aligned LCSC, let $(\Lambda,G,\varphi)$ a category system, and let 
\[F=\{(\gamma_1,h_1),\ldots,(\gamma_n,h_n)\}\subseteq \Lambda\rtimes^\varphi G.\] 
Then, given $(\alpha,g)\in \Lambda\rtimes^\varphi G$, $F\in \mathsf{FE}(\alpha,g)$ of $\Lambda\rtimes^\varphi G$ if and only if $\{\gamma_1,\ldots,\gamma_n\}\in \mathsf{FE}(\alpha)$ of $\Lambda$. 
\end{remark} 

Since $\Lambda\rtimes^\varphi  G$ is a finitely aligned LCSC when so is $\Cat$, the following results are a direct translation of Theorem \ref{eq_top_free} and \ref{eq_minimal}. 

\begin{proposition}[{\cite[Proposition 7.7]{OP1}}]
Let $\Lambda$ be a finitely aligned LCSC, let $G$ be a countable discrete group acting on $\Cat$, and let $(\Lambda,G,\varphi)$ be a category system. If either $\mathcal{G}_{\text{tight}}(\Semi_{\Lambda\rtimes^\varphi  G})$ is Hausdorff or $\hat{\mathcal{E}}_\infty =\hat{\mathcal{E}}_{\text{tight}}$, then the following statements are equivalent:
	\begin{enumerate}
		\item $\Grpd_{tight}(\Semi_{\Lambda\rtimes^\varphi  G})$ is effective.
		\item Given $(\alpha,a),(\beta,b)\in \Lambda\rtimes^\varphi G$ with $r(\alpha,a)=r(\beta,b)$ and $s(\alpha,a)=s(\beta,b)$, if $(\alpha,a)(\delta,d)\Cap (\beta,b)(\delta,d)$ for every $(\delta,d)\in s((\alpha,a))(\Lambda\rtimes^\varphi G)$, then there exists $F\in\mathsf{FE}(s(\alpha,a))$ such that $(\alpha,a)(\gamma,d)=(\beta,b)(\gamma,d)$ for every $(\gamma,d)\in F$.
		\item Given $\alpha,\beta\in \Lambda$, $a,b\in G$ with $r(\alpha)=r(b)$ and $a^{-1}\cdot s(\alpha)=b^{-1}\cdot s(\beta)$, if $\alpha (a\cdot \delta ) \Cap \beta(b\cdot \delta)$ for every $\delta\in (a^{-1}\cdot s(\alpha))\Lambda$ then there exists $F\in\mathsf{FE}(a^{-1}\cdot s(\alpha))$ such that $\alpha (a\cdot \gamma)=\beta(b\cdot\gamma)$ and $\varphi(a,\gamma)=\varphi(b,\gamma)$ for every $\gamma\in F$.

	\end{enumerate}
\end{proposition}

\begin{proposition}[{\cite[Proposition 7.8]{OP1}}]
Let $\Lambda$ be a finitely aligned LCSC, let $G$ be a countable discrete group acting on $\Cat$, and let $(\Lambda,G,\varphi)$ be a category system. Then, the following statements are equivalent: 
	\begin{enumerate}
		\item $\Grpd_{tight}(\Semi_{\Lambda\rtimes^\varphi  G})$ is minimal.
		\item For every $(\alpha,a),(\beta,b)\in \Lambda\rtimes^\varphi G$ there exists $F\in \mathsf{FE}((\alpha,a))$ such that for each $(\gamma,g)\in F$, $s(\beta,b)(\Lambda\rtimes^\varphi G) s(\gamma,g)\neq \emptyset$.
		\item For every $\alpha,\beta\in \Lambda$ there exists $F\in \mathsf{FE}(\alpha)$ such that for each $\gamma\in F$, there exists $g\in G$ with $s(\beta)\Lambda  (g\cdot s(\gamma))\neq \emptyset$.

	\end{enumerate}
\end{proposition}

Then, by an analogous argument to that of Theorem \ref{theorem_simpleLambda}, we have the following result:

\begin{theorem}[{\cite[Theorem 7.9]{OP1}}]\label{theorem_simpleSystem}
Let $\Lambda$ be a countable, finitely aligned LCSC, let $G$ be a countable discrete group acting on $\Cat$, and let $(\Lambda,G,\varphi)$ be a category system. If $\mathcal{G}_{\text{tight}}(\Semi_{\Lambda\rtimes^\varphi  G})$ is Hausdorff and amenable, then the following statements are equivalent: 
\begin{enumerate}
\item $C^*(\Semi_{\Lambda\rtimes^\varphi  G})$ is simple.
\item The following properties hold:
\begin{enumerate}
\item Given $\alpha,\beta\in \Lambda$, $a,b\in G$ with $r(\alpha)=r(b)$ and $a^{-1}\cdot s(\alpha)=b^{-1}\cdot s(\beta)$, if $\alpha (a\cdot \delta ) \Cap \beta(b\cdot \delta)$ for every $\delta\in (a^{-1}\cdot s(\alpha))\Lambda$ then there exists $F\in\mathsf{FE}(a^{-1}\cdot s(\alpha))$ such that $\alpha (a\cdot \gamma)=\beta(b\cdot\gamma)$ and $\varphi(a,\gamma)=\varphi(b,\gamma)$ for every $\gamma\in F$.
\item For every $\alpha,\beta\in \Lambda$ there exists $F\in \mathsf{FE}(\alpha)$ such that for each $\gamma\in F$, there exists $g\in G$ with $s(\beta)\Lambda  (g\cdot s(\gamma))\neq \emptyset$.
\end{enumerate}
\end{enumerate}
\end{theorem}
The same result holds for the $K$-algebra $K\Semi_{\Lambda\rtimes^\varphi  G}$ over any field $K$, without requiring amenability for $\mathcal{G}_{\text{tight}}(\Semi_\Lambda)$. There are results determining necessary and sufficient conditions on $\Cat$ and $\Lambda\rtimes^\varphi  G$ for $\Grpd_{\text{tight}}(\Semi_\Cat)$ and $\Grpd_{\text{tight}}(\Semi_{\Lambda\rtimes^\varphi  G})$ to be amenable \cite[Section 8]{OP1}.\vspace{.2truecm}

To end this section, we will have a look on the case of $\Lambda=E^*$, where $E$ is a countable graph. In order to fix the relation between $\Grpd_{tight}(\Semi_{G,E})$ and $\Grpd_{tight}(\Semi_{{E^*}\rtimes^\varphi  G})$, we first need to state the relation between $\Semi_{G,E}$ and $\Semi_{{E^*}\rtimes^\varphi  G}$. On the one hand, we have
\[\Semi_{G,E}=\{(\alpha, g, \beta) : \alpha,\beta\in E^*, g\in G, s(\alpha)=g\cdot s(\beta)\}.\]
On the other hand,
\[\Semi_{{E^*}\rtimes^\varphi  G}=\langle \elmap{(\alpha, g)}{(\beta,h)} : \alpha,\beta\in E^*, g,h\in G, g^{-1}\cdot s(\alpha)=h^{-1}\cdot s(\beta)\rangle.\]
Since $(x,g)\in ({{E^*}\rtimes^\varphi  G})^{-1}$ for all $x\in E^0$, $g\in G$, we have that 
\[\elmap{(\alpha, g)}{(\beta,h)} =\elmap{(\alpha, g)}{(h^{-1}\cdot s(\beta),h^{-1})}\elmap{(s(\beta), h)}{(\beta, \ideG)}=\elmap{(\alpha, gh^{-1})}{(\beta,\ideG)}.\]
Moreover, since $E^*$ is singly aligned, so is ${{E^*}\rtimes^\varphi  G}$ \cite[Proposition 3.12(ii)]{BKQS}. Thus, by \cite[Theorem 3.2]{DGKMW}, 
\[\Semi_{{E^*}\rtimes^\varphi  G}=\{ \elmap{(\alpha, g)}{(\beta,\ideG)} : \alpha,\beta\in E^*, g\in G, s(\alpha)=g\cdot s(\beta)\}\,.\]
Hence, the map
\[
\begin{array}{cccc}
\pi: & \Semi_{G,E} &  \rightarrow &  \Semi_{{E^*}\rtimes^\varphi  G} \\
 & (\alpha, g, \beta) & \mapsto  &   \elmap{(\alpha, g)}{(\beta,\ideG)}
\end{array}
\]
is a well-defined, onto $\ast$-semigroup homomorphism. Let us characterize when $\pi$ is injective. To this end, take $(\alpha, g, \beta), (\gamma, h, \delta)\in \Semi_{G,E}$ such that
\[ \elmap{(\alpha, g)}{(\beta,\ideG)}=\pi(\alpha, g, \beta)=\pi(\gamma, h, \delta)= \elmap{(\gamma, h)}{(\delta,\ideG)}.\]
Being both equal functions, they must have the same domain, i.e. $(\beta, \ideG)({E^*}\rtimes^\varphi  G)=(\delta, \ideG)({E^*}\rtimes^\varphi  G)$. Since $({{E^*}\rtimes^\varphi  G})^{-1}=E^0\times G$, we conclude that $\beta=\delta$. Moreover, $\tau^{(\alpha, g)}=\tau^{(\gamma, h)}$ on their common domain, so that for every $\lambda \in (g^{-1}\cdot s(\alpha))E^*$ and for every $\ell \in G$ we have
\[(\alpha(g\cdot \lambda), \varphi(g,\lambda)\ell)=\tau^{(\alpha, g)}(\lambda, \ell)=\tau^{(\gamma, h)}(\lambda, \ell)=(\gamma(h\cdot\lambda), \varphi(h,\lambda)\ell).\]
Since the self-similar action of $G$ on $E^*$ preserves lengths of paths, we conclude that $\alpha=\gamma$, and that for every $\lambda \in (g^{-1}\cdot s(\alpha))E^*$ we have $g\cdot \lambda=h\cdot\lambda$ and $\varphi(g,\lambda)=\varphi(h,\lambda)$. Thus, the existence of nontrivial kernel for $\pi$ is equivalent to the existence of $g\in G$, $\alpha \in E^*$ such that for all $\lambda \in s(\alpha)E^*$ satisfies $g\cdot\lambda=\lambda$ and $\varphi(g,\lambda)=\ideG$; in other words, injectivity of $\pi$ is equivalent to the fact that the self-similar action of $G$ on $E^*$ is faithful on vertex-based trees of $E$. Notice that if $(E,G, \varphi)$ is pseudo-free, then the above condition is trivially fulfilled, so that $\pi$ will be an isomorphism in this case. Moreover, being ${{E^*}\rtimes^\varphi  G}$ singly aligned, we see that it is right cancellative exactly when $(E,G, \varphi)$ is pseudo free. In this case, not only $\Semi_{G,E}\cong \Semi_{{E^*}\rtimes^\varphi  G}$, but also they are weak semilattices, so their associated tight groupoids are Hausdorff.

Now, we look at the relation between the corresponding tight groupoids $\Grpd_{tight}(\Semi_{G,E})$ and $\Grpd_{tight}(\Semi_{{E^*}\rtimes^\varphi  G})$. First, notice that the idempotent semilattices of $\Semi_E$, $\Semi_{G,E}$ and $\Semi_{{E^*}\rtimes^\varphi  G}$ coincide, so that the spaces of filters, ultrafilters and tight filters are the same (up to natural isomorphism). Furthermore, the partial actions $\Semi_{G,E}\curvearrowright \hat{\mathcal{E}}_0$ and $\Semi_{{E^*}\rtimes^\varphi  G}\curvearrowright \hat{\mathcal{E}}_0$ are $\pi$-equivariant. Moreover, the germ relationship is compatible with $\pi$. Thus, $\pi$ induces a continuous, open, onto groupoid homomorphism
\[
\begin{array}{cccc}
\Phi: & \Grpd_{tight}(\Semi_{G,E}) & \rightarrow  & \Grpd_{tight}(\Semi_{{E^*}\rtimes^\varphi  G})  \\
 & [\alpha, g, \beta; \eta] &  \mapsto &   [\elmap{(\alpha, g)}{(\beta, \ideG)}; \eta] .
\end{array}
\]
Using the Morita equivalence reduction shown in Theorem \ref{Thm:Desin} (\cite[Theorem 3.2]{EPS}), we can assume that $E$ is row-finite with no sources or sinks, whence $\hat{\mathcal{E}}_\infty=\hat{\mathcal{E}}_{\text{tight}}=E^\infty$; let us reduce to this case, to simplify the computations. We now show that $\Phi$ is injective. To do this, let $[\alpha, g, \beta; \eta] \in \ker \Phi$. This means that $\elmap{(\alpha, g)}{(\beta, \ideG)}$ is an idempotent. According to the calculations done before, this happens exactly when $\alpha=\beta$ and for every $\lambda \in s(\alpha)E^*$ we have that $g\cdot\lambda=\lambda$ and $\varphi(g,\lambda)=\ideG$. Pick $\lambda$ any initial segment in $\eta\in E^\infty$. Notice that $\lambda \in s(\alpha)E^*$, and thus $(\alpha\lambda, \ideG, \alpha\lambda)\in \alpha\eta$ (seen as a filter), while
\[(\alpha, g, \alpha)\cdot (\alpha\lambda, \ideG, \alpha\lambda)= (\alpha(g\lambda), \varphi (g, \lambda), \alpha\lambda)=(\alpha\lambda, \ideG, \alpha\lambda).\]
Hence, by the germ relation, if $\eta=\lambda\hat{\eta}$, then 
\[[\alpha, g, \alpha; \alpha\eta)]=[\alpha\lambda, \ideG, \alpha\lambda; \alpha\lambda\hat{\eta})]\in \Grpd_{tight}(\Semi_{G,E})^{(0)}.\]
Thus, $\Phi$ is a homeomorphism and an isomorphism of groupoids. This guarantees that, independent of the choice for representing the self-similar graph system $(G,E,\varphi)$, their associated tight groupoids -and hence their algebras- are the same.

\section{Zappa-Sz\'ep products of groupoids acting on LCSC categories}

An obvious generalization of the previous construction, as well as that of \cite{LRRW18}, is to consider Zappa-Sz\'ep products of \emph{groupoids} acting on left cancellative small categories. This require to state the right definition of what a groupoid action means.

By technical reasons (associated with some applications of the results in this section), it is convenient to introduce a notion of \emph{length function} defined on a  left cancellative small category.

\begin{definition}
	Let $\Cat$ be a LCSC and let $\Gamma\subseteq  Q$ be a submonoid of a group $Q$ with unit element $\ideQ$, and such that $\Gamma\cap\Gamma^{-1}=\{\ideQ\}$. A map $\dmap:\Lambda\to\Gamma$ is called  a \emph{length function} if:
	\begin{enumerate}
		\item $\dmap(\alpha\beta)=\dmap(\alpha)\dmap(\beta)$ for every $\alpha,\beta\in\Lambda$ with $s(\alpha)=r(\beta)$.
\end{enumerate}
A length function is said to satisfy the \emph{weak factorization property} \textbf{(WFP)} if for every $\alpha\in \Lambda$ and $\gamma_1,\gamma_2\in \Gamma$ with $\dmap(\alpha)=\gamma_1\gamma_2$, there are $\alpha_1,\alpha_2\in \Lambda$ with $s(\alpha_1)=r(\alpha_2)$, $\dmap(\alpha_i)=\gamma_i$ for $i=1,2$, such that $\alpha=\alpha_1\alpha_2$, and moreover for every $\beta_1,\beta_2\in \Lambda$ with $s(\beta_1)=r(\beta_2)$, $\dmap(\beta_i)=\gamma_i$ for $i=1,2$, such that $\alpha=\beta_1\beta_2$, there exist $g_1,g_2\in \Cat^{-1}$ such that $\beta_1=\alpha_1g_1$ and $\beta_2=g_2\alpha_2$. 	 
\end{definition}
Observe that, given a LCSC $\Cat$, there always exists what we call the \emph{trivial length function} $\dmap:\Cat \to \Gamma$, defined by $\dmap(\alpha)=\ideQ$ for every $\alpha\in \Cat$. 

The above definition is motivated by the one given in \cite[Section 3]{LV2}, where the authors only consider the case when $\Gamma=\N^k$.

We recall that a length function $\dmap:\Cat \to \Gamma$ satisfies \emph{unique factorization property} \textbf{(UFP)} if, for all $\alpha\in \Cat$ and $\gamma_1,\gamma_2\in \Gamma$ with $\dmap(\alpha)=\gamma_1\gamma_2$, there are unique $\alpha_1,\alpha_2\in \Cat$ with $\alpha=\alpha_1\alpha_2$ and $\dmap(\alpha_i)=\gamma_i$ for $i=1,2$ (see e.g. \cite{KP}). 

\begin{remark}\label{rema_WFP} 
	Let $\Cat$ be a LCSC  and let $\Gamma\subseteq Q$ be a submonoid of a group $Q$  with unit element $\ideQ$ such that $\Gamma\cap\Gamma^{-1}=\{\ideQ\}$. Then, given  a length function   $\dmap:\Lambda\to\Gamma$  satisfying the WFP, we have:
	\begin{enumerate}
		\item  $\dmap^{-1}(\ideQ)=\Cat^{-1}$.		
		\item 	$\dmap:\Cat \to \Gamma$ satisfies the UFP if and only if it satisfies the WFP and $\Cat$ has no inverses. 
	\end{enumerate}
\end{remark}

\begin{definition}
	A LCSC $\Cat$ is called \emph{action-free} if the action of $\Lambda^{-1}$ on $\Lambda$ is free: whenever $g\gamma=\gamma$ for some $\gamma\in \Cat$ and $g\in \Cat^{-1}$, then $g=r(\gamma)$.
\end{definition}

Observe that if $\Cat$ has no inverses, then $\Cat$ is action-free.

\begin{lemma}[{\cite[Lemma 2.10]{OP2}}]\label{pseudofree_cat_right_cancellative}
	Let $\Cat$ be a LCSC and let $\dmap:\Cat\to \Gamma$ be a length function satisfying the WFP.  Then, $\Cat$ is right cancellative if and only if $\Cat$ is action-free.
\end{lemma}

Now, we define the action of a groupoid on a LCSC, essentially following the same schema as in \cite{LRRW18}. The first step is to define the notion of partial isomorphism of a small category $\Cat$, inspired by \cite[Definition 3.1]{LRRW18}. 

\begin{definition}\label{Def: Partial_Isomorphism}
Let $\Cat$ be a LCSC  with length function $\dmap:\Lambda\to\Gamma$. Let $\PisoG$ be the set of partial isomorphisms of $\Cat$, $f\in \PisoG$ such that satisfies: 
\begin{enumerate}
\item $f:v\Cat \rightarrow w\Cat$ is a bijection  for some  $v,w\in \Cat^0$,
\item $f(\alpha\Cat)=f(\alpha)\Cat$ for every $\alpha\in v\Cat$,
\item $f(v)=w$, and
\item $\dmap(f(\alpha))=\dmap(\alpha)$ for every $\alpha\in v\Cat$.
\end{enumerate}
We denote by $v=:d(f)$ (the \emph{domain of $f$}) and  by $w=:c(f)$ (the \emph{codomain of $f$}).
\end{definition}

\begin{remas}\label{Rem:ConditionsAction}$\mbox{ }$\vspace{.1truecm}

\begin{enumerate}
	\item Observe that condition $(3)$ does not follow from the previous conditions. 
	Indeed, given $g\in \Cat^{-1}$, the map $f:s(g)\Cat\to r(g)\Cat$ given by $f(\gamma)=g\gamma$ for every $\gamma\in s(g)\Cat$ satisfies conditions $(1)$-$(2)$ but not $(3)$, whenever $g\in \Cat^{-1}\setminus \Cat^0$.
	\item Given $f:v\Cat\to w\Cat$ in $\PisoG$ and $g\in v\Cat^{-1}$, we have that $f(gg^{-1})=w$ by condition $(3)$, and hence by condition $(2)$ there exists $h\in \Cat$ such that $f(g)h=w$, so $f(g)\in \Cat^{-1}$.
\end{enumerate}
\end{remas}

\begin{lemma}[{\cite[Lemma 3.3]{OP2}}]\label{cocycle}
	Let $\Cat$ be a LCSC  with length function $\dmap:\Lambda\to\Gamma$. Then, given $f\in \PisoG$ and $\alpha\in d(f)\Cat$, there exists a unique function $f\vert_{\alpha}:s(\alpha)\Cat\to s(f(\alpha))\Cat$ in $\PisoG$ such that $f(\alpha\beta)=f(\alpha)f\vert_{\alpha}(\beta)$ for every $\beta\in s(\alpha)\Cat$.
\end{lemma}

\begin{lemma}[{\cite[Lemma 3.4]{OP2}}]\label{cocycle_prop}
	Let $\Cat$ be a LCSC with length function $\dmap:\Lambda\to\Gamma$. Then, given $\alpha,\beta\in \Cat$ with $s(\alpha)=r(\beta)$ and $f\in \PisoG$ with $r(\alpha)=d(f)$, we have that:
	\begin{enumerate}
		\item $f\vert_{d(f)}=f$.
		\item $s(f(\alpha))=f\vert_{\alpha}(s(\alpha))=c(f\vert_{\alpha})$.
		\item $\id_{s(\alpha)}=(\id_{r(\alpha)})\vert_{\alpha}$.
		\item $f\vert_{\alpha\beta}=(f\vert_{\alpha})\vert_{\beta}$.
		%\item $f(\alpha\beta)=f(\alpha)f_{|\alpha}(\beta)$.
	\end{enumerate} 
\end{lemma}
\vspace{.2truecm}

The next step is to show that $\PisoG$ has a natural groupoid structure.

\begin{lemma}[{\cite[Lemma 3.5]{OP2}}]\label{Piso_comp}
Let $\Cat$ be a LCSC with length function $\dmap:\Lambda\to\Gamma$. Then, given $f,g\in \PisoG$ such that $d(f)=c(g)$,  the composition $f\circ g\in \PisoG$.	
\end{lemma}

For each $v\in \Cat^0$, we define $\id_v: v\Cat \rightarrow v\Cat$ as the identity map on $v\Cat$, so that $\id_v\in \PisoG$.

\begin{lemma}[{\cite[Lemma 3.6]{OP2}}]\label{Piso_inv}
	Let $\Cat$ be a LCSC  with length function $\dmap:\Lambda\to\Gamma$, and let $f\in \PisoG$. Then $f^{-1}\in \PisoG$.	
\end{lemma}

Now, using Lemmas \ref{Piso_comp} and \ref{Piso_inv}, we can define a groupoid structure in $\PisoG$:
 
 \begin{definition}\label{proposition:PIso-groupoid}
 	If $\Cat$ is  a LCSC with length function $\dmap:\Lambda\to\Gamma$, then $\PisoG$ is a discrete groupoid where, given $f,g\in \PisoG$, the product $fg$ is defined as the composition $f\circ g$ whenever $d(f)=c(g)$, and $f^{-1}$ is the set-theoretical inverse of $f$.  Moreover, we can identify the unit space of $(\PisoG)^{(0)}=\{\id_v:v\in \Cat^0 \}$ with $\Cat^0$.
 \end{definition}

With this in mind, we can define the notion of action of a groupoid on a LCSC that we need:

\begin{definition}\label{definition_GroupoidAction}
Let $\Lambda$ be a LCSC with length function $\dmap:\Lambda\to\Gamma$, and let $\Grpd$ be a discrete groupoid. An action of $\Grpd$ on $\Cat$ is a groupoid homomorphism $\phi: \Grpd  \rightarrow   \PisoG$. Given $g\in \Grpd$ and $\alpha\in d(\phi(g))\Cat$, we write the action of $g$ on $\alpha$ by $g\cdot \alpha:=\phi(g)(\alpha)$.
\end{definition}

The identification $(\PisoG)^{(0)}=\Cat^0$ implies that $\phi(\Grpdu)\subseteq\Cat^0$, and $\phi(s(g))=d(\phi(g))$ and $\phi(r(g))=c(\phi(g))$ for any $g\in \Grpd$. Thus,  
\[\phi(g): \phi(s(g))\Cat \rightarrow \phi(r(g))\Cat\,.\] 

\begin{remas}\label{Rem:2.9} Let $\Cat$ be a LCSC  with length function $\dmap:\Lambda\to\Gamma$, and let $\phi:\Grpd\to \PisoG$ be a groupoid action on $\Cat$. Then:
\begin{enumerate} 
	\item  Suppose that $\Delta^0:=\phi(\Grpdu)\subsetneqq \Cat^0$. Let us define  
	\[\Delta:=\{\alpha\in \Cat\mid s(\alpha),r(\alpha)\in \Delta^0 \}\,.\] 
	Then $\Delta$ is a LCSC, and $\phi(\Grpd)\subseteq \text{PIso}(\Delta,\dmap)$.
	However, $\Delta$ is not necessarily finitely aligned if so is $\Cat$.  Thus, such a restriction will affect the arguments. Hence, we always assume that, given an action of $\Grpd$ on $\Cat$, either $\phi(\Grpdu)=\Cat^0$, or that the restricted subcategory $\Delta$ inherits the finitely aligned property. In practice, this means that, after restricting if necessary, we assume that $\phi(\Grpdu)= \Cat^0$.
	
	\item  Suppose that $\phi_{|\Grpdu}$ is not an injective map. Then it can happen that there exists $(g,h)\notin \Grpdc$ (so that $s(g) \neq r(h)$), but $(\phi(g),\phi(h))\in (\PisoG)^{(2)}$, (and thus $\phi(s(g))=d(\phi(g))=c(\phi(h))=\phi(r(h))$). Hence, $\phi(g)\phi(h)\in \PisoG\setminus \phi(\Grpd)$. This will affect defining composability of elements in the corresponding Zappa-Sz\'ep product. Since our model should include non-faithful self-similar actions of a groupoid $\Grpd$ on $\Cat$, we cannot skip that case. Then, for a non-injective action, we need to include the condition that for every $g,h\in \Grpd$, then $s(g)=r(h)$ if and only if $d(\phi(g))=\phi(s(g))=\phi(r(h))=c(\phi(h))$, or equivalently, that $\phi_{|\Grpdu}$ is injective. 
	\end{enumerate}
\end{remas}

According to Remarks \ref{Rem:2.9}, we assume that an action of a discrete groupoid $\Grpd$ on a small category $\Cat$ through a (not necessarily injective) groupoid homomorphism $\phi: \Grpd\rightarrow \PisoG$, satisfies that $\phi_{|\Grpdu}$ is a bijection. Therefore, we identify $\Grpdu=\Cat^{(0)}=\PisoG^{(0)}$, omitting $\phi$.\vspace{.2truecm}

Now, we are ready to define the Zappa-Sz\'ep product of a left cancellative small category. This is inspired by the construction of the self-similar graph of a groupoid on a finite graph \cite[Section 3]{LRRW18} and the Zappa-Sz\'ep product of a group acting on a left cancellative small category \cite{BKQS,OP1}. This has also been done recently in \cite{LV2}, where they construct Zappa-Sz\'ep products of groupoids acting on higher-rank graphs. 

First, we state the abstract notion of self-similar action, in order to boil up the exact definition of a 1-cocycle that we need.

\begin{definition}\label{Def:SelfSimilarAction}
Let $\Grpd$ be a discrete groupoid acting on a LCSC $\Cat$ with length function $\dmap:\Lambda\to\Gamma$. We say that the action is \emph{self-similar} if, for every $g\in \Grpd$ and for every $\alpha\in s(g)\Cat$, there exists $h\in \Grpd$ such that
\[g\cdot (\alpha\mu)=(g\cdot \alpha)(h\cdot\mu)\]
for every $\mu\in s(\alpha)\Cat$.
\end{definition}

Clearly, in the above definition $h$ depends on $g$ and $\alpha$. A natural question is to decide whether such a $h$ is unique. Suppose that $h_1, h_2\in \Grpd$ satisfies $g\cdot (\alpha\mu)=(g\cdot \alpha)(h_i\cdot \mu)$ ($i=1,2)$ for every $\mu\in s(\alpha)\Cat$. Then
\[(g\cdot \alpha)(h_1\cdot \mu) =(g\cdot \alpha)(h_2\cdot \mu) .\]
Since  $\Lambda$ is left cancellative, then 
\[(h_1\cdot \mu)=(h_2\cdot \mu)\]
for every $\mu\in s(\alpha)\Lambda$. 
Hence,  we can only guarantee that $\phi(h_1)=\phi(h_2)$ and $h_1=h_2$ hold only when $\phi$ is injective.
\vspace{.2truecm}
\ 
We use this fact to define a suitable notion of a cocycle for such an action. We essentially follow \cite[Section 4]{BKQS}, taking care of the fact that the action is partial.

\begin{definition}\label{Def:GroupoidCocycle}
Let $\Cat$ be a LCSC with length function $\dmap:\Lambda\to\Gamma$, and let $\Grpd$ be a discrete groupoid acting on $\Cat$. Consider the set
\[\Grpd{}_s\times_r\Cat:=\{ (g,\alpha)\in \Grpd\times \Cat \mid s(g)=r(\alpha)\}\,.\]
A \emph{(partial) cocycle of the action of $\Grpd$ on $\Cat$} is a function $\varphi :  \Grpd{}_s\times_r\Cat \rightarrow  \Grpd$ satisfying the cocycle identity 
\[\varphi (gh,\alpha)=\varphi(g, h\cdot\alpha)\varphi(h,\alpha)\] 
for any $(g,h)\in \Grpdc$ 
 and any $\alpha \in \Cat$ such that $s(h)=r(\alpha)$. Observe, in particular, that $\varphi(r(\alpha),\alpha)\in \Grpdu$ for every $\alpha\in \Cat$.
\end{definition}

Now, we state the properties that a cocycle must have (in a similar list to that of \cite[Section 7]{OP1}). But instead of imposing them, we will try to deduce from the definition, as in \cite[Section 3]{LRRW18}.\vspace{.2truecm}

The first step is to fix the requirement to guarantee that the action of $\Grpd$ is compatible with the composition in $\Cat$. We introduce minimal requirements to guarantee this fact when constructing our ``self-similar'' actions, and also to ensure that we can associate a suitable small category to the action and the cocycle.

\begin{definition}\label{Def:CategoryCocycle}
Let $\Cat$ be a LCSC with length function $\dmap:\Lambda\to\Gamma$, and let $\Grpd$ be a discrete groupoid acting on $\Cat$. A (partial) cocycle $\varphi$ for the action of $\Grpd$ on $\Cat$ is said to be a \emph{category cocycle} if for every $g\in \Grpd$, $\alpha\in\Cat$ with $s(g)=r(\alpha)$ and every $\beta\in s(\alpha)\Cat$ we have that:
\begin{enumerate}
\item[(CC1)] $\varphi(g,d(\phi(g)))=g$.	
\item[(CC2)] $\varphi(g,\alpha\beta)=\varphi(\varphi(g,\alpha), \beta)$.
\item[(CC3)] $g\cdot (\alpha\beta)=(g\cdot \alpha)(\varphi(g,\alpha)\cdot \beta)$.
\end{enumerate}
If we have an action of $\Grpd$ on $\Cat$, and $\varphi$ is a category cocycle for this action, we say that $(\Cat,\dmap, \Grpd, \varphi)$ is a \emph{category system}.
\end{definition}

\begin{remark}\label{remark_ss}
	Let $\Cat$ be a LCSC with length function $\dmap:\Lambda\to\Gamma$, let $\Grpd$ be a discrete groupoid acting on $\Cat$, and let $(\Cat,\dmap, \Grpd, \varphi)$ be a category system. Then, by condition $(CC3)$ in Definition \ref{Def:CategoryCocycle},
the action of $\Grpd$ on $\Cat$ is self-similar. Moreover, given $g\in \Grpd$ and $\alpha\in s(g)\Cat$, we have that $\phi(\varphi(g,\alpha))=\phi(g)\vert_{\alpha}$. 
\end{remark}

\begin{lemma}[{\cite[Lemma 4.5]{OP2}}]
Let $\Cat$ be a LCSC with length function $\dmap:\Lambda\to\Gamma$, and let $\Grpd$ be a discrete groupoid acting on $\Cat$. Let $\varphi$ be a category cocycle for the action of $\Grpd$ on $\Cat$. Then  for every $g\in \Grpd$, $\alpha\in\Cat$ with $s(g)=r(\alpha)$ we have that:
\begin{enumerate}
	\item[(CC4)] $s(g\cdot \alpha)=\varphi(g,\alpha)\cdot s(\alpha)=r(\varphi(g,\alpha))$.
	\item[(CC5)] $s(\alpha)=\varphi(r(\alpha),\alpha)$.
\end{enumerate}
\end{lemma}

\begin{remark}
	Similarly, $(CC1)$ and $(CC2)$ will follow from $(CC3)$, modulo the action homomorphism $\phi$. In particular, \underline{when the action homomorphism is injective}, $(CC1)$ and $(CC2)$ are direct consequences of $(CC3)$.
 \end{remark}

\begin{example}
	If $\Cat$ is a LCSC with length function $\dmap:\Lambda\to\Gamma$, then $\varphi:\PisoG {}_d\times_r \Cat\to \PisoG$ defined by $\varphi(f,\alpha)=f\vert_{\alpha}$ is a category cocyle (Lemmas \ref{cocycle} \& \ref{cocycle_prop}), and thus $(\Cat,\dmap,\PisoG,\varphi)$ is a category system. 
\end{example}
The property of being a category cocycle is the correct version of \cite[(2.3)]{EP2}. Let us determine which are the basic properties satisfied by a category cocycle. These properties, jointly with Definition \ref{Def:CategoryCocycle}, are analogous to those proved in \cite[Lemma 3.4 \& Proposition 3.6]{LRRW18}.

\begin{proposition}[{\cite[Proposition 4.8]{OP2}}]\label{Prop:left cancellative small category}
Let $\Cat$ be a LCSC with length function $\dmap:\Lambda\to\Gamma$, let $\Grpd$ be a discrete groupoid acting on $\Cat$, and let $(\Cat,\dmap, \Grpd, \varphi)$ be a category system. Then, for every $g\in \Grpd$, every $\alpha, \beta\in \Cat$ with $s(g)=r(\alpha)$, $s(\alpha)=r(\beta)$, we have that:
\begin{enumerate}
\item $r(g\cdot \alpha)=g\cdot r(\alpha)$.
\item $s(\varphi(g,\alpha))=s(\alpha)$.
\item $\varphi(g,\alpha)^{-1}=\varphi(g^{-1}, g\cdot \alpha)$.
\end{enumerate}
\end{proposition}

Now, we fix the definition of Zappa-Sz\'ep product for a category system $(\Cat,\Grpd,\varphi)$, similar to that in \cite[Section 4]{BKQS}. 

\begin{definition}\label{definition_zappa-szed2}
Let $\Cat$ be an LCSC with length function $\dmap:\Lambda\to\Gamma$, let $\Grpd$ be a discrete groupoid acting on $\Cat$, and let $(\Cat,\dmap,\Grpd,\varphi)$ be a category system. We denote by $\Cat \Join^\varphi \Grpd$ the set
\[\Cat \Join^\varphi \Grpd:=\Cat{}_s\times_r\Grpd=\{ (\alpha, g)\in \Cat \times \Grpd \mid s(\alpha)=r(g)\}\,,\]
with distinguished elements
\[(\Cat \Join^\varphi \Grpd)^{(0)}:= \Cat^0 \Join^\varphi \Grpdu\,.\]
We equip this pair of sets with range and source maps $r,s:\Lambda \Join^\varphi G \rightarrow (\Lambda \Join^\varphi G)^{(0)}$ defined by 
\[r(\alpha, g):=(r(\alpha), r(\alpha))\text{ and } s(\alpha, g):=(s(g),s(g))=(g^{-1}\cdot s(\alpha), s(g)),\]
for every $(\alpha, g)\in \Cat \Join^\varphi \Grpd$. Observe that the latter equality holds because $g\cdot s(g)=r(g)=s(\alpha)$, whence $g^{-1}\cdot s(\alpha)=g^{-1}\cdot r(g)=s(g)$.

Moreover, given $(\alpha, g), (\beta, h)\in \Cat \Join^\varphi \Grpd$ with $s(\alpha, g)=r(\beta, h)$, we define the composition of two elements as follows:
\[(\alpha, g)(\beta, h):=(\alpha(g\cdot \beta), \varphi(g,\beta)h).\]
\end{definition}

We now adapt the arguments of \cite[Proposition 4.6]{BKQS} and \cite[Proposition 4.13]{BKQS} to prove the next results. 

\begin{proposition}[{\cite[Proposition 4.10]{OP2}}]
	Let $\Cat$ be a LCSC with length function $\dmap:\Lambda\to\Gamma$, let $\Grpd$ be a discrete groupoid acting on $\Cat$, and let $(\Cat,\dmap,\Grpd,\varphi)$ be a category system. If given $(\alpha, g), (\beta, h)\in \Cat \Join^\varphi \Grpd$ with $s(\alpha, g)=r(\beta, h)$, we define the composition of these two elements by
	\[(\alpha, g)(\beta, h):=(\alpha(g\cdot \beta), \varphi(g,\beta)h)\,,\]
	then $\Cat \Join^\varphi \Grpd$ is a small category.
\end{proposition}

\begin{lemma}[{\cite[Lemma 4.11]{OP2}}]\label{lemma_inv}
	Let $\Cat$ be a LCSC with length function $\dmap:\Lambda\to\Gamma$, let $\Grpd$ be a discrete groupoid acting on $\Cat$, and let $(\Cat, \dmap,\Grpd,\varphi)$ be a category system. Then, $(\Cat \Join^\varphi \Grpd)^{-1}=\Cat^{-1}{}_s\times_r\Grpd$.
\end{lemma}

\begin{proposition}[{\cite[Proposition 4.12]{OP2}}]\label{proposition_LC and FA}
Let $\Cat$ be a LCSC with length function $\dmap:\Lambda\to\Gamma$, let $\Grpd$ be a discrete groupoid acting on $\Cat$, and let $(\Cat,\dmap,\Grpd,\varphi)$ be a category system. Then:
\begin{enumerate}
\item $\Cat\Join^\varphi \Grpd$ is left cancellative.
\item If $\Cat$ is finitely (singly) aligned, then:
\begin{enumerate}
\item $(\alpha, g)(\Cat\Join^\varphi \Grpd) \cap (\beta,h)(\Cat\Join^\varphi \Grpd)=(\alpha\Cat \cap \beta\Cat)_s\times _r \Grpd$.
\item $\Cat\Join^\varphi \Grpd$ is finitely (singly) aligned.
\end{enumerate}
\end{enumerate}
\end{proposition}

In addition, we can prove an analog of \cite[Lemma 4.15]{BKQS}

\begin{lemma}[{\cite[Lemma 4.13]{OP2}}]\label{lemma_FiniteCovers}
Let $\Cat$ be a LCSC with length function $\dmap:\Lambda\to\Gamma$, let $\Grpd$ be a discrete groupoid acting on $\Cat$, and let $(\Cat,\dmap,\Grpd,\varphi)$ be a category system. Take any $(v, v)\in (\Cat\Join^\varphi \Grpd)^0$, and let $F\subseteq (v, v)(\Cat\Join^\varphi \Grpd)$. Set
\[H:=\{ \alpha \in v\Cat \mid \text{ there is } g\in \Grpd \text{ such that } (\alpha, g)\in F\}.\]
Then, $F$ is exhaustive at $(v, v)$ if and only if $H$ is exhaustive at $v$.
\end{lemma}

\begin{proposition}[{\cite[Proposition 4.14]{OP2}}]\label{ZP_WFP}
Let $\Cat$ be a LCSC with length function $\dmap:\Lambda\to\Gamma$ satisfying the WFP, let $\Grpd$ be a discrete groupoid acting on $\Cat$, and let $(\Cat,\dmap,\Grpd,\varphi)$ be a category system.  Then the map $\dmap: \Cat\Join^\varphi \Grpd\to \Gamma$ defined by $\dmap(\alpha,g)=\dmap(\alpha)$ for every $(\alpha,g)\in \Cat\Join^\varphi \Grpd$ is a length function satisfying the WFP.
\end{proposition}

\begin{remark}
	Observe that, in the proof of Proposition \ref{ZP_WFP}, if $\Cat$ satisfies the WFP as defined in \cite{LV2}, then so does the Zappa-Sz\'ep product. 
\end{remark}

\begin{definition}
Let $\Cat$ be a LCSC with length function $\dmap:\Lambda\to\Gamma$, and let $\Grpd$ be a discrete groupoid acting on $\Cat$. We say that a category system $(\Cat,\dmap,\Grpd,\varphi)$ is \emph{pseudo-free} if $\Cat$ is an action-free category, and whenever $g\cdot \alpha=f\alpha$ and $\varphi(g,\alpha)=s(\alpha)$ for some $\alpha\in \Cat$, $g\in \Grpd$ and $f\in \Cat^{-1}$, then $g\in \Grpdu$.
\end{definition}

\begin{proposition}[{\cite[Proposition 4.17]{OP2}}]\label{ZP_RC}
Let $\Cat$ be a LCSC with length function $\dmap:\Lambda\to\Gamma$ satisfying the WFP, let $\Grpd$ be a discrete groupoid acting on $\Cat$, and let $(\Cat,\dmap,\Grpd,\varphi)$ be a category system.  Then, $\Cat\Join^\varphi \Grpd$ is right cancellative if and only if $(\Cat,\dmap,\Grpd,\varphi)$ is pseudo-free.
\end{proposition}
  
 \begin{remark}\label{ZP_comp} Let $\Cat$ be a LCSC with length function $\dmap:\Lambda\to\Gamma$, let $\Grpd$ be a discrete groupoid acting on $\Cat$, let $(\Cat,\dmap,\Grpd,\varphi)$ be a category system, and let $\Cat \Join^\varphi \Grpd$ be the associated Zappa- Sz\'ep product.
 	\begin{enumerate}
 		\item Let $(\alpha,g),(\beta,f)\in \Cat \Join^\varphi \Grpd$,  and suppose that $(\alpha,g)\leq(\beta,f)$. Then, there exist $(\delta,h)\in \Cat \Join^\varphi \Grpd $ with $r(\delta,h)=s(\alpha,g)$, that is, $r(\delta)=g^{-1}\cdot s(\alpha)$ such that 
 	\[(\beta,f)=(\alpha,g)(\delta,h)=(\alpha(g\cdot \delta),\varphi(g,\delta)h)\,.\]
 	Hence,  $\alpha(g\cdot \delta)= \beta$ and $f=\varphi(g,\delta)h$. Thus, $g\cdot \delta=\sigma^\alpha(\beta)$ by left cancellation, and so $\delta=g^{-1}\cdot \sigma^\alpha(\beta)$ and $h=\varphi(g,g^{-1}\cdot\sigma^\alpha(\beta))^{-1}f=\varphi(g^{-1},\sigma^\alpha(\beta))f$ because of the cocycle identity. Therefore,  
 	\[(\alpha,g)\leq(\beta,f) \text{ if and only if } \alpha\leq \beta\,,\] 
	and then we have that 
 	\[\sigma^{(\alpha,g)}(\beta,f)=(g^{-1}\cdot \sigma^\alpha(\beta),\varphi(g^{-1},\sigma^\alpha(\beta))f)\,.\]
 	The above observation, together with  Proposition \ref{proposition_LC and FA}(2)(a), shows that the map 
 	\[F\mapsto F\Join^\varphi \Grpd:=\{(\alpha,g)\in \Cat \Join^\varphi \Grpd:  \alpha\in F \} \]
 	is a bijection between $\Cat^*$ and $(\Cat \Join^\varphi \Grpd)^*$. This bijection clearly restricts to a bijection of their tight (maximal) hereditary upper-directed subsets, so we identify $(\Cat \Join^\varphi \Grpd)_{tight}$ with $\Cat_{tight}$. 
 	\item  Let $(\alpha,g),(\beta,f)\in \Cat \Join^\varphi \Grpd$ with $s(\alpha,g)=s(\beta,f)$, and let $\elmap{(\alpha,g)}{(\beta,f)}\in \Semi_{\Cat \Join^\varphi \Grpd}$. Observe that, since  $(\beta,f)(\Cat \Join^\varphi \Grpd)=(\beta,s(f))(\Cat \Join^\varphi \Grpd)$, we have that 
\[\elmap{(\alpha,g)}{(\beta,f)}=\elmap{(\alpha,gf^{-1})}{(\beta,s(f))}.\] 
 	
 	\item Given $F\in \Cat_{tight}$ and $(\alpha,g),(\beta,f)\in \Cat \Join^\varphi \Grpd$ with $s(\alpha,g)=s(\beta,f)$, we have that 
 	\[\elmap{(\alpha,g)}{(\beta,f)}\cdot (F\Join^\varphi \Grpd)=F'\Join^\varphi \Grpd\,,\]
where $F':=\bigcup\limits_{\substack{\beta\leq \gamma\\ \gamma\in F}}[\alpha((gf^{-1})\cdot \sigma^\beta (\gamma))]\in \Cat_{tight}$.
 	 	\end{enumerate}
 \end{remark}  

Now, we are ready to characterize the properties of $\Grpd_{tight}(\Cat \Join^\varphi \Grpd)$, following the same process as in the two previous sections. First, we give necessary and sufficient conditions for $\Grpd_{tight}(\Cat \Join^\varphi \Grpd)$ to be a Hausdorff groupoid. 
 
 \begin{proposition}[{\cite[Proposition 6.10]{OP2}}]\label{ZP_hausddorff}
 Let $\Cat$ be a LCSC with length function $\dmap:\Lambda\to\Gamma$, let $\Grpd$ be a discrete groupoid acting on $\Cat$, let $(\Cat,\dmap,\Grpd,\varphi)$ be a category system, and let $\Cat \Join^\varphi \Grpd$ be the associated Zappa-Sz\'ep product.  Then, $\Grpd_{tight}(\Cat \Join^\varphi \Grpd)$ is a Hausdorff groupoid if and only if, given $\alpha,\beta\in \Cat$ and $g\in\Grpd$ with $g^{-1}\cdot s(\alpha)=s(\beta)$, there exists a finite subset  $H\subset s(\beta)\Cat$ such that 
 \begin{enumerate}
 \item $\alpha(g\cdot\gamma)=\beta\gamma$ and $\varphi(g,\gamma)=s(\gamma)$ for every $\gamma\in H$, and
 \item for every $\delta\in s(\beta)\Cat$ with $\alpha(g\cdot \delta)=\beta\delta$ and $\varphi(g,\delta)=s(\gamma)$, then $ \delta\Cap \gamma$ for some $\gamma\in H$.
 \end{enumerate}
\end{proposition}

\begin{corollary}[{\cite[Corollary 6.11]{OP2}}]\label{ZP_pseudofree_Hausdorff}
 	Let $\Cat$ be a LCSC with length function $\dmap:\Lambda\to\Gamma$ satisfying the WFP, let $\Grpd$ be a discrete groupoid acting on $\Cat$, let $(\Cat,\dmap,\Grpd,\varphi)$ be a category system, and let $\Cat \Join^\varphi \Grpd$ be the associated Zappa-Sz\'ep product. If $(\Cat,\dmap,\Grpd,\varphi)$ is pseudo-free, then $\Grpd_{tight}(\Cat \Join^\varphi \Grpd)$ is a Hausdorff groupoid. 
 \end{corollary}  

Now, we are ready to characterize minimality and effectiveness for the corresponding tight groupoid. The following results are straightforward translations of the results in \cite[Section 6 \& 7]{OP1}.

\begin{proposition}[{\cite[Proposition 6.12]{OP2}}]\label{ZP_TP}
Let $\Cat$ be a LCSC with length function $\dmap:\Lambda\to\Gamma$, let $\Grpd$ be a discrete groupoid acting on $\Cat$, let $(\Cat,\dmap,\Grpd,\varphi)$ be a category system, and let $\Cat \Join^\varphi \Grpd$ be the associated Zappa-Sz\'ep product.  If $\Grpd_{tight}(\Cat \Join^\varphi \Grpd)$ is Hausdorff or $\Cat_{tight}=\Cat^{**}$, then the following statements are equivalent:
\begin{enumerate}
\item $\Grpd_{tight}(\Cat \Join^\varphi \Grpd)$ is topologically free.
\item Given $(\alpha,a),(\beta,b)\in \Cat \Join^\varphi \Grpd$ with $r(\alpha,a)=r(\beta,b)$ and $s(\alpha,a)=s(\beta,b)$, if $(\alpha,a)(\delta,d)\Cap (\beta,b)(\delta,d)$ for every $(\delta,d)\in s(\alpha,a)(\Cat \Join^\varphi \Grpd)$ then there exists $C\in\mathsf{FE}(s(\alpha,a))$ such that $(\alpha,a)(\gamma,d)=(\beta,b)(\gamma,d)$ for every $(\gamma,d)\in C$.
\item Given $\alpha,\beta\in \Cat$, $a,b\in \Grpd$ with $r(a)=s(\alpha)$, $r(b)=s(\beta)$, $r(\alpha)=r(\beta)$ and $a^{-1}\cdot s(\alpha)=b^{-1}\cdot s(\beta)$, if $\alpha (a\cdot \delta ) \Cap \beta(b\cdot \delta)$ for every $\delta\in (a^{-1}\cdot s(\alpha))\Lambda$ then there exists $C\in\mathsf{FE}(a^{-1}\cdot s(\alpha))$ such that $\alpha (a\cdot \gamma)=\beta(b\cdot\gamma)$ and $\varphi(a,\gamma)=\varphi(b,\gamma)$ for every $\gamma\in F$.
\end{enumerate}
\end{proposition}

\begin{proposition}[{\cite[Proposition 6.13]{OP2}}]\label{ZP_Min}
Let $\Cat$ be a LCSC with length function $\dmap:\Lambda\to\Gamma$, let $\Grpd$ be a discrete groupoid acting on $\Cat$, let $(\Cat,\dmap,\Grpd,\varphi)$ be a category system, and let $\Cat \Join^\varphi \Grpd$ be the associated Zappa-Sz\'ep product.   If $\Grpd_{tight}(\Cat \Join^\varphi \Grpd)$ is Hausdorff or $\Cat_{tight}=\Cat^{**}$, then the following statements are equivalent:
\begin{enumerate}
\item $\Grpd_{tight}(\Cat \Join^\varphi \Grpd)$ is minimal.
\item For every $(\alpha,a),(\beta,b)\in \Lambda\Join^\varphi \Grpd$ there exists $C\in \mathsf{FE}(\alpha,a)$ such that for each $(\gamma,g)\in C$, $s(\beta,b)(\Lambda\Join^\varphi G) s(\gamma,g)\neq \emptyset$.
\item For every $\alpha,\beta\in \Lambda$ there exists $C\in \mathsf{FE}(\alpha)$ such that for each $\gamma\in C$,  there exists $g\in \Grpd$  with $s(\beta)\Lambda  (g\cdot s(\gamma))\neq \emptyset$.
\end{enumerate}
\end{proposition}

Thus, we conclude with the following result:

\begin{theorem}[{\cite[Theorem 6.14]{OP2}}]\label{theorem_simpleSystemnou}
Let $\Cat$ be a LCSC with length function $\dmap:\Lambda\to\Gamma$, and let $\Grpd$ be a discrete groupoid acting on $\Cat$. Let $(\Cat,\dmap,\Grpd,\varphi)$ be a category system and $\Cat \Join^\varphi \Grpd$ be the associated Zappa-Sz\'ep product.   If $\Grpd_{tight}(\Cat \Join^\varphi \Grpd)$ is Hausdorff and amenable, then the following statements are equivalent:
\begin{enumerate}
\item $C^*_r(\Grpd_{tight}(\Cat \Join^\varphi \Grpd))$ (the reduced groupoid $C^*$-algebra) is simple.
\item The following properties hold:
\begin{enumerate}
\item Given $\alpha,\beta\in \Lambda$, $a,b\in \Grpd$ with $r(\alpha)=r(\beta)$ and $a^{-1}\cdot s(\alpha)=b^{-1}\cdot s(\beta)$, if $\alpha (a\cdot \delta ) \Cap \beta(b\cdot \delta)$ for every $\delta\in (a^{-1}\cdot s(\alpha))\Lambda$ then there exists $C\in\mathsf{FE}(a^{-1}\cdot s(\alpha))$ such that $\alpha (a\cdot \gamma)=\beta(b\cdot\gamma)$ and $\varphi(a,\gamma)=\varphi(b,\gamma)$ for every $\gamma\in C$.
\item For every $\alpha,\beta\in \Lambda$ there exists $C\in \mathsf{FE}(\alpha)$ such that for each $\gamma\in C$,  there exists $g\in \Grpd$  with $s(\beta)\Lambda  (g\cdot s(\gamma))\neq \emptyset$.
\end{enumerate}
\end{enumerate}
\end{theorem}

The same result holds for the $K$-algebra $K(\Lambda\Join^\varphi  G)$ over any field $K$, without requiring amenability for $\mathcal{G}_{\text{tight}}(\Semi_\Lambda)$. There are results determining necessary and sufficient conditions on $\Cat$ and $\Lambda\rtimes^\varphi  G$ for $\Grpd_{\text{tight}}(\Semi_{\Lambda\rtimes^\varphi  G})$ be amenable \cite[Section 7]{OP2}. 

\section{Homology and K-theory for self-similar graph algebras}

Homology and K-theory play a central role in understanding the structure of $C^*$-algebras, as well as in their classification. In the particular case of groupoid $C^*$-algebras, the groupoid homology of an ample groupoid $\Grpd$ shares many similarities with the K-theory of its reduced $C^*$-algebra, which is closely related to the K-theory when the isotropy groups are torsion-free~\cite{PY22, miller2024isomorphisms}. Matui had conjectured~\cite{Matui16} a connection $K_i(C^*_\lambda(\Grpd)) \cong \bigoplus_{q \geq 0} H_{2q+i}(\Grpd)$ for $i=0,1$, known as the HK property. However, torsion in the isotropy becomes an obstruction~\cite{Scarparo20} for the HK property, even for principal groupoids~\cite{Deeley}.\vspace{.2truecm}

When specializing to self-similar graph $C^*$-algebras, the homology of the groupoid $\Grpd_{(G,X)}$ associated to a self-similar group action $(G,X,\sigma)$, and the K-theory of its Ne\-kra\-she\-vych algebra $\mathcal O_{(G,X)}$, provide fundamental invariants for the self-similar group. The K-theory is particularly relevant in light of the Kirchberg--Phillips Theorem~\cite{Kirchberg, Phillips}, as Nekrashevych algebras are often purely infinite and simple. Nekrashevych computed the K-theory of Nekrashevych algebras in the case of certain iterated monodromy groups~\cite{N2}. He used a two-step approach: first, he related the K-theory of $\mathcal O_{(G,X)}$ with that of its gauge-invariant subalgebra, and then described it as an inductive limit of matrix amplifications of $C^*(G)$. A similar approach, at the groupoid level, was used by Ortega and S\'anchez~\cite{OS} to study the homology of the groupoid associated to a certain self-similar action of the infinite dihedral group. Deaconu also outlined a general strategy along these lines for both K-theory and homology~\cite{DeaconuSelfSimilar}.\vspace{.2truecm}

Miller and Steinberg compute homology and K-theory in much greater generality, by relating the groupoid $\Grpd_{(G,X)}$ and the $C^*$-algebra $\mathcal O_{(G,X)}$ of a self-similar group action $(G,X,\sigma)$ directly to $G$ and $C^*(G)$, respectively. This provides a way to compute the homology and K-theory in entirely group-theoretic terms. Their strategy is to use that self-similar groupoid actions are exactly self-correspondences of discrete groupoids~\cite{AKM22}, and to take advantage of their induced mappings in homology~\cite{miller2023ample} and K-theory. This allows them to extend the results from self-similar groups to self-similar graphs, and even to self-similar groupoid actions on graphs, without extra effort. Once this is done, they show that Matui's computation of the homology of graph groupoids~\cite{Matui} (see also~\cite{ON}) and Nyland and Ortega's computation of the homology of Exel--Pardo--Katsura groupoids \cite{O}) follow directly from their results. \vspace{.2truecm}

Let us have a quick look at their results. First, let us recall the definitions for \'etale correspondences and their composition~\cite{AKM22}, which are basic for understanding this work.

\begin{definition}
Let $\Grpd$ and $\mathcal{H}$ be ample groupoids. An \textit{\'etale correspondence} $\Omega \colon \Grpd \to \mathcal{H}$ is a $\Grpd$-$\mathcal{H}$-bispace $\Omega$ such that the right action $\Omega \curvearrowleft \mathcal{H}$ is free, proper and \'etale. We write $\ranT \colon \Omega \to \Grpd^0$ and $\sour \colon \Omega \to \mathcal{H}^0$ for the anchor maps. That $\Omega \curvearrowleft \mathcal{H}$ is \'etale means that $\sour \colon \Omega \to \mathcal{H}^0$ is \'etale, while free and proper together mean that the map $\ranT \times \sour \colon \Omega \rtimes \mathcal{H} \to \Omega \times \Omega$ is a closed embedding. We say $\Omega$ is \textit{proper} if the induced map $\Omega/\mathcal{H} \to \Grpd^0$ is proper.
\end{definition}

For an ample groupoid $\mathcal{H}$ and a free, proper, \'etale right $\mathcal{H}$-space $\Omega$ and $(\omega_1,\omega_2) \in \Omega \times_{\Omega/\mathcal{H}} \Omega$, we write $\langle \omega_1, \omega_2 \rangle$ for the unique $h \in \mathcal{H}$ satisfying $\omega_2 = \omega_1 h$. The map 
\[\langle - , - \rangle \colon \Omega \times_{\Omega/\mathcal{H}} \Omega \to \mathcal{H}\] 
is continuous~\cite[Lemma~3.4]{AKM22}.

One formulation of the notion of self-similar groups is via a proper correspondence from a group to itself (see Nekrashevych~\cite{N-book} where the term ``covering bimodule'' is used).

\begin{definition}\label{def:v.corr}
A \textit{self-similar groupoid action} $(\Grpd,\mathcal X)$ consists of a discrete groupoid $\Grpd$ and an \'etale correspondence $\mathcal X\colon \Grpd\to \Grpd$ with anchor maps $\ranT,\sour\colon \mathcal X \to \Grpd^0$. We may refer to $\mathcal X$ as a \textit{self-similarity} of $\Grpd$. The self-similar groupoid action is called row-finite if $\mathcal X$ is proper. An object/vertex $v\in \Grpd^0$ is called regular if $0<\vert\ranT^{-1}(v)/\Grpd\vert<\infty$; otherwise $v$ is called singular. One says that a singular vertex $v$ is a source if $\ranT^{-1}(v)=\emptyset$, and an infinite receiver otherwise.
\end{definition}

Regular objects form an invariant subset $\Grpd^0_{\reg}$ of the unit space and thus determine an invariant subgroupoid ${\Grpd_{\reg}} = \Grpd \vert_{\Grpd^0_{\reg}}$. To be consistent with~\cite{EP2}, we say that $\mathcal X$ is pseudo-free if the left action of $\Grpd$ is free.   

There is a reformulation of this notion in terms of actions on graphs that can be found in~\cite[Example~4.4]{AKM22}.  

\begin{definition}\label{def:v.graph}
A \textit{self-similar groupoid action} $(\Grpd,E,\sigma)$ is a discrete groupoid $\Grpd$ whose unit space is the vertex set $E^0$ of a directed graph $\ranT, \sour \colon E \to E^0$, with a left action $\Grpd \curvearrowright E$ (not necessarily an action by graph partial automorphisms, see \cite{AKM22}) with anchor $\ranT \colon E \to E^0$, written $(g,e) \mapsto g(e)$, and a $1$-cocycle $\sigma \colon \Grpd \times_{E^0} E \to \Grpd$, written $(g,e) \mapsto g \vert_e$, such that $\ranT(g \vert_e) = \sour(g(e))$, $\sour(g \vert_e) = \sour(e)$ and $(hg)\vert_e = h \vert_{g(e)} g \vert_e$ (equivalently, $\sigma \colon \Grpd \ltimes E \to \Grpd$ is a groupoid homomorphism with unit space map $\sigma \vert_E = \sour \colon E \to E^0$). The element $g\vert_e$ is called the section of $g$ at $e$.
\end{definition}

The translation between the two definitions is as follows:
\begin{enumerate}
\item If $(\Grpd,E,\sigma)$ is a self-similar groupoid action in the sense of Definition~\ref{def:v.graph}, then the associated correspondence is $\mathcal X_{(\Grpd,E)} = E\times_{E^0} \Grpd$ with $\ranT(e,g) = \ranT(e)$, $\sour(e,g)=\sour(g)$, right action $(e,h)g =(e,hg)$ and left action $g(e,h) = (g(e),g|_eh)$.  Moreover, $E$ is a set of representatives for $\mathcal X_{(\Grpd,E)}/\Grpd$.
\item If $(\Grpd,\mathcal X)$ is a self-similar groupoid action in the sense of Definition~\ref{def:v.corr}, then we can choose a transversal $E \subseteq \mathcal X$ to $\mathcal X/\Grpd$.
By freeness of the right action, each $x\in \mathcal X$ can be uniquely written in the form $eg$ with $e\in E$, $g\in G$.  Putting $E^0=\Grpd^0$, we can define $\ranT,\sour\colon E\to E^0$ by restriction.  If $g\in \Grpd$ and $\sour(g)=\ranT(e)$, then $ge=g(e)g\vert_e$ for a unique $g(e)\in E$ and $g\vert_e\in \Grpd$.  This defines a left action of $\Grpd$ with anchor $\ranT$ isomorphic to the action of $\Grpd$ on $\mathcal X/\Grpd$, and $\sigma(g,e) = g\vert_e$ defines a $1$-cocycle with  $\ranT(g\vert_e) = \sour(g(e))$, $\sour(g\vert_e) = \sour(e)$ and $(hg)\vert_e = h\vert_{g(e)}g\vert_e$. Thus, $(\Grpd,E,\sigma)$ is a self-similar groupoid action in the sense of definition~\ref{def:v.graph} and $\mathcal X_{(\Grpd,E)}\to \mathcal X$ given by $(e,g)\mapsto eg$ is an isomorphism of correspondences.  It is easy to see that $\mathcal X_{(\Grpd,E)}$ is proper if and only if $E$ is row-finite. Moreover, a vertex $v\in E^0$ is regular if and only if it is not a source and not an infinite receiver.
\end{enumerate}

By a \emph{self-similar group action} we shall always mean in these notes $(G,\mathcal X)$ with $\mathcal X$ a nontrivial proper correspondence over the group $G$, equivalently $(G,E,\sigma)$ with finite alphabet $E$ of size $|E| \geq 2$ (although the case of nonproper correspondences was considered in~\cite{SS2}).\vspace{.2truecm}

Now, we define the corresponding $C^*$-algebra. Given an \'etale correspondence $\mathcal X\colon \Grpd\to \mathcal{H}$ of a discrete groupoids, we can extend the map $\langle - , - \rangle \colon \mathcal X \times_{\mathcal X/\mathcal{H}} \mathcal X\to \mathcal{H}$ to a map $\langle -,-\rangle\colon \mathcal X\times \mathcal X\to \mathcal{H}\cup \{0\}$ by 
\[\langle x,y\rangle = \begin{cases} h, & \text{if}\ xh=y,\\ 0, & \text{else.}\end{cases}\] 
Trivially, $\langle x,yh\rangle = \langle x,y\rangle h$, $\langle y,x\rangle = \langle x,y\rangle^{-1}$, $\langle x,x\rangle = \sour(x)$ and $\langle gx,y\rangle = \langle x,g^{-1} y \rangle$.  
For example, if $(\Grpd,E,\sigma)$ is a self-similar groupoid action in the sense of Definition~\ref{def:v.graph}, then for $\mathcal X_{(\Grpd,E)}$, one has $\langle (e,g),(f,h)\rangle = g^{-1} h\delta_{e,f}$.

We write $M_{\mathcal X} \colon C^*(\Grpd) \to C^*(\Grpd)$ for the $C^*$-correspondence of a groupoid self-similarity $\mathcal X \colon \Grpd \to \Grpd$. The space $M_{\mathcal X}$ is densely spanned by elements $m_x$ for $x \in \mathcal X$, and $C^*(\Grpd)$ is densely spanned by the partial isometries $u_g \in C^*(\Grpd)$ for $g \in \Grpd \cup \{0\}$ with $u_0 = 0$, and thus the $C^*$-correspondence $M_{\mathcal X}$ is determined by the relations:
\begin{enumerate}
\item $u_g \cdot m_x = m_{gx} \delta_{\sour(g), \ranT(x)}$,
\item $m_x \cdot u_g = m_{xg} \delta_{\sour(x), \ranT(g)}$, and
\item $\langle m_x, m_y \rangle = u_{\langle x,y \rangle}$
\end{enumerate}
for $x,y \in \mathcal X$ and $g \in \Grpd$. The regular objects determine an ideal $C^*(\Grpd_{\reg})$ in $C^*(\Grpd)$. Moreover, this ideal acts by compact operators on $M_{\mathcal X}$ because for each $v \in {\Grpd^0}_{\reg}$ the projection $u_v$ acts as the compact operator $\sum_{x\Grpd \in \ranT_{\mathcal X/\Grpd}^{-1} (v)} m_xm_x^*$. The Hilbert module $M_{\mathcal X}$ is full if and only if there are no sinks. 

\begin{definition}
The \textit{$C^*$-algebra} $\mathcal O_{\mathcal X}$ of a self-similar groupoid action $(\Grpd,\mathcal X)$ is the relative Cuntz--Pimsner algebra of the $C^*$-correspondence $M_{\mathcal X}$ over $C^*(\Grpd)$ with respect to the ideal $C^*(\Grpd_{\reg})$, and the \textit{Toeplitz algebra} $\mathcal T_{\mathcal X}$ of $(\Grpd,\mathcal X)$ is the Pimsner--Toeplitz algebra of $M_{\mathcal X}$. Concretely, this means that 
 $\mathcal T_{\mathcal X}$ is the universal $C^*$-algebra generated by elements $u_g$ for $g \in \Grpd$ and $m_x$ for $x \in \mathcal X$ with the groupoid relations $u_gu_h = u_{gh} \delta_{\sour(g),\ranT(h)}$ on $(u_g)_{g \in \Grpd}$ and the relations:
\begin{enumerate}
\item[(T1)] $u_g m_x = m_{gx} \delta_{\sour(g),\ranT(x)}$,
\item[(T2)] $m_x u_g = m_{xg} \delta_{\sour(x),\ranT(g)}$, and
\item[(T3)] $m_{x}^* m_{y} = u_{\langle x, y \rangle}$, where $u_0=0$,
\end{enumerate}
for $x,y \in \mathcal X$ and $g \in \Grpd$. The $C^*$-algebra $\mathcal O_{\mathcal X}$ has the additional relations:
\begin{itemize}
\item[(CK)] $u_v = \sum_{x\Grpd \in \ranT_{\mathcal X/\Grpd}^{-1} (v)} m_{x} m_{x}^*$ 
\end{itemize}
for $v \in \Grpd^0_{\reg}$. 
\end{definition}

\begin{remark}
Katsura's nonrelative Cuntz--Pimsner algebra is not used because it is acceptable that $C^*(\Grpd_{\reg})$ may act non-faithfully on $M_{\mathcal X}$.
\end{remark}

If $(\Grpd,E,\sigma)$ is a self-similar groupoid action, then for a right $\Grpd$-transversal $E \subseteq \mathcal X$ the elements $(m_e)_{e \in E} \subseteq M_{\mathcal X}$ generate $M_{\mathcal X}$ as a Hilbert $C^*(\Grpd)$-module, so the above presentation reduces to the following:

\begin{proposition}[{\cite[Proposition 2.6]{MS}}]
Let $(\Grpd,E,\sigma)$ be a self-similar groupoid action with correspondence $\mathcal X = \mathcal X_{\Grpd,E} \colon \Grpd \to \Grpd$. The Toeplitz algebra $\mathcal T_{(\Grpd,E)} = \mathcal T_{\mathcal X}$ of $(\Grpd,E,\sigma)$ is the universal $C^*$-algebra generated by elements $u_g$ for $g \in \Grpd$ and $m_e$ for $e \in E$ with the groupoid relations $u_gu_h = u_{gh} \delta_{\sour(g),\ranT(h)}$ on $(u_g)_{g \in \Grpd}$ and the relations:
\begin{enumerate}
\item[(T1)] $u_g m_e = \delta_{\sour(g),\ranT(e)} m_{g(e)} u_{g \vert_e}$, and
\item[(T3)] $m_{e}^* m_{f} = \delta_{e,f} u_{\sour(f)}$
\end{enumerate}
for $e,f \in E$ and $g \in \Grpd$. The $C^*$-algebra $\mathcal O_{(\Grpd,E)} = \mathcal O_{\mathcal X}$ has the additional relations:
\begin{itemize}
\item[(CK)] $u_v = \sum_{\ranT(e) = v} m_{e} m_{e}^*$ for each regular vertex $v \in E^0_{\reg}$.
\end{itemize}
\end{proposition}

As in~\cite{EP2,N2}, they construct a groupoid model for $\mathcal O_{\mathcal X}$ (or, equivalently, $\mathcal O_{(G,E)}$), which is moreover the tight groupoid of an inverse semigroup, although we will not explore this construction.

Recall the definition of \emph{groupoid homology}. If $X$ is a space with a basis of compact (Hausdorff) open sets and $A$ is an abelian group (written additively), then $AX$ denotes the abelian group of mappings $X\to A$ spanned by elements of the form $a1_U$ where $a \in A$ and $1_U$ is the characteristic function of a compact open set $U \subseteq X$. If $X$ is Hausdorff, these are precisely the compactly supported locally constant mappings $X\to A$.

The construction $X\mapsto AX$ is functorial with respect to \'etale maps and contravariantly functorial with respect to proper maps. If $f\colon X\to Y$ is \'etale, then $f_\ast\colon AX\to AY$ is given by $f_*(g)(y)= \sum_{x\in f^{-1}(y)}g(x)$. If $p\colon X\to Y$ is a proper map, then $p^*\colon AY \to AX $ is given by $p^*(f) =f\circ p$. Notice that an open inclusion $i$ is \'etale with $i_*$ extension by $0$,   and a closed inclusion $i$ is proper with $i^*$ restriction of functions.

We recall now the definition of the homology of an ample groupoid $\Grpd$ with coefficients in an abelian group $A$ via the formulation of Matui~\cite{Matui}.  There are \'etale maps $d_i\colon \Grpd^n\to \Grpd^{n-1}$ for $n\geq 2$ and $i=0,\ldots, n$ given by
\begin{equation}\label{eq:face.maps}
d_i(g_1,\ldots,g_n) = \begin{cases}(g_2,\ldots, g_n), & \text{if}\ i=0,\\ (g_1,\ldots,g_ig_{i+1},\ldots,g_{n}), & \text{if}\ 1\leq i\leq n-1,\\ (g_1,\ldots, g_{n-1}), & \text{if}\ i=n.\end{cases}
\end{equation}
We define $d_0,d_1\colon \Grpd^1\to \Grpd^0$ by $d_0(g) = \sour(g)$ and $d_1(g)=\ranT(g)$, which are again \'etale.
It is well known that these maps satisfy the semisimplicial identities, 
and so we can define a chain complex $C_\bullet(\Grpd,A)$ with $C_n(\Grpd,A) =  A \Grpd^n$ for $n\geq 0$, and $\partial_n\colon C_n(\Grpd,A)\to C_{n-1}(\Grpd,A)$ given by $\partial_n=\sum_{i=0}^n(-1)^i(d_i)_*$ for $n\geq 1$.  As usual, we take $\partial_0=0$.  The homology of this chain complex is denoted $H_\bullet(\Grpd,A)$.  When $A=\mathbb Z$, we often write $C_\bullet(\Grpd)$ and $H_\bullet(\Grpd)$. There is a picture of groupoid homology in terms of the Tor functor with $H_\bullet(\Grpd,A) \cong \text{Tor}^{\Z \Grpd}_\bullet(\Z \Grpd^0, A \Grpd^0)$, and for general $\Z \Grpd$-modules $M$ we set $H_\bullet(\Grpd ; M) = \text{Tor}^{\Z \Grpd}_\bullet(\Z \Grpd^0, M)$;  this is the left derived functor of the $\Grpd$-coinvariants $M\mapsto M_{\Grpd}=M/\langle 1_Um-1_{\sour(U)}m\rangle$, where $U$ ranges over compact open bisections of $\Grpd$ and $m$ ranges over $M$. It is known that groupoid homology is invariant under Morita equivalence.\vspace{.2truecm}

Now, we are ready to present the main results of this section.

\begin{theorem}[{\cite[Theorem 5.1]{MS}}]\label{t:main.homology}
Let $(\Grpd,\mathcal X)$ be a self-similar groupoid action. Then, for each abelian group $A$, there is a long exact sequence in homology 
\[ \cdots \to H_{n+1}(\Grpd_{\mathcal X},A) \to H_n(\Grpd_{\reg},A) \to H_n(\Grpd,A) \to H_n(\Grpd_{\mathcal X},A) \to \cdots \]
where the middle map is $\id - H_n(\mathcal X_{\reg},A)$, with $\mathcal X_{\reg}\colon\Grpd_{\reg}\to \Grpd$ the restriction of $\mathcal X$ to $\Grpd_{\reg}$.
\end{theorem}

Applying \cite[Propositions 4.6 \& 4.8]{MS} to $H_n(\mathcal X_{\reg})$ in Theorem~\ref{t:main.homology}, yields the following ``groups only'' description of the long exact sequence:

\begin{corollary}[{\cite[Corollary 5.2]{MS}}]\label{c:homology.transversal}
Let $(\Grpd,E,\sigma)$ be a self-similar groupoid action on a graph $E$ with cocycle $\sigma$. Let $T^0 \subseteq \Grpd^0$ be a transversal for $\Grpd$ and set $T^0_{\reg} = T^0 \cap E^0_{\reg}$. Then, for each abelian group $A$, there is a long exact sequence in homology 
\[ \cdots \to H_{n+1}(\Grpd_{(\Grpd,E)},A) \to H_n(\Grpd^v_v,A) \to H_n(\Grpd^w_w,A) \to H_n(\Grpd_{(\Grpd,E)},A) \to \cdots \]
where the middle map is $\id-\Phi_n$, and $\Phi_n$ admits the following description: fix, for each $v\in T^0_{\reg}$, a left $\Grpd^v_v$-transversal $T_v$ to $vE$. Consider, for $e \in T_v$, the virtual homomorphism $\sigma_e \colon \Grpd_e \to G^{\sour(e)}_{\sour(e)}$, $\sigma_e(g)=g\vert_e$. 
For each $w \in \sour(\bigsqcup_{v\in T^0_{\reg}}T_v)$, choose $h_w \in \Grpd$ with $\sour(h_w) = w$ and $\ranT(h_w) = t(w) \in T^0$ and set $c_w \colon \Grpd^w_w \to \Grpd^{t(w)}_{t(w)}$ to be the homomorphism $g \mapsto h_w^{-1} g h_w$. 
Then \[ \Phi_n = \bigoplus_{v\in T^0_{\reg}}\sum_{e \in T_v} H_n(c_{\sour(e)},A)\circ H_n(\sigma_e,A)\circ \mathrm{tr}^{\Grpd_v^v}_{\Grpd_e},\] and it is induced by the chain map
\begin{align*}
\bigoplus_{v \in T^0_{\reg}} C_\bullet(\Grpd^v_v) & \to \bigoplus_{w \in T^0} C_\bullet(\Grpd^w_w) \\
(g_1,\dots,g_m) & \mapsto \sum_{\ranT(e)=\sour(g_m)} (c_{\sour(e)}(g_1\vert_{g_2\cdots g_m(e)}),\ldots, c_{\sour(e)} (g_m\vert_e)).
\end{align*}
In particular, $(\Phi_0)_{w,v} =\vert\ranT^{-1}(v)\cap \sour^{-1}(\Grpd w)\vert \colon H_0(\Grpd^v_v,A) \to H_0(\Grpd^w_w,A)$.
\end{corollary}

When applying to actions of groups we are in the following situation, where the statement becomes considerably simpler.

\begin{corollary}[{\cite[Corollary 5.4]{MS}}]\label{c:transitive.transfer}
Let $(G,E,\sigma)$ be a self-similar group action on a finite alphabet $E$ of cardinality at least $2$ with cocycle $\sigma$. For $e\in E$ and let $\sigma_e\colon G_e\to G$ be the virtual endomorphism $\sigma_e(g)=g\vert_e$. Then there is a long exact sequence
\[ \cdots \to H_{n+1}(\Grpd_{(G,E)},A) \to H_n(G,A) \to H_n(G,A) \to H_n(\Grpd_{(G,E)},A) \to \cdots \]
where the middle map is $\id-\Phi_n$, and $\Phi_n = \sum_{e \in T} H_n(\sigma_e,A)\circ \mathrm{tr}^{G}_{G_e}$ for any $G$-transversal $T \subseteq E$ and is induced by the chain map 
\begin{align*}
C_\bullet(G) & \to  C_\bullet(G) \\
(g_1,\dots,g_m) & \mapsto \sum_{e \in E} (g_1\vert_{g_2\cdots g_m(e)},\ldots, g_m\vert _e ).
\end{align*}
In particular, $H_0(\Grpd_{(G,E)})\cong \mathbb Z/(\vert E\vert-1)\mathbb Z$ and $H_1(\Grpd_{(G,E)}))\cong \coker (\id-\Phi_1)$ where 
\[
\begin{array}{crcl}
\Phi_1\colon & {G}^\text{ab} & \to  & {G}^\text{ab}  \\
 & g[G,G] & \mapsto  &   \sum_{e\in E} g\vert_e[G,G]
\end{array}.
\] 
\end{corollary}

With respect to K-theory, they obtain the following result:

\begin{theorem}[{\cite[Theorem 5.5]{MS}}]\label{t:main.k}
Let $(\Grpd,E,\sigma)$ be a self-similar groupoid action on a graph $E$ with cocycle $\sigma$. Let $T^0 \subseteq \Grpd^0$ be a transversal for $\Grpd$ and set $T^0_{\reg} = T^0 \cap E^0_{\reg}$. Then, there is a six-term sequence in K-theory
\[\begin{tikzcd}[arrow style=math font]
	{\bigoplus\limits_{v \in T^0_{\reg}} K_0(C^*(\Grpd^v_v))} & {\bigoplus\limits_{w \in T^0} K_0(C^*(\Grpd^w_w))} & {K_0(\OGE)} \\
	{K_1(\OGE)} & {\bigoplus\limits_{w \in T^0} K_1(C^*(\Grpd^w_w))} & {\bigoplus\limits_{v \in T^0_{\reg}} K_1(C^*(\Grpd^v_v))}
	\arrow["{1-\Phi_0}", from=1-1, to=1-2]
	\arrow[from=1-2, to=1-3]
	\arrow[from=1-3, to=2-3]
	\arrow[from=2-1, to=1-1]
	\arrow[from=2-2, to=2-1]
	\arrow["{1-\Phi_1}", from=2-3, to=2-2]
\end{tikzcd}\]
where $\Phi_i$ admits the following description for $i=0,1$: fix, for each $v\in T^0_{\reg}$, a left $\Grpd^v_v$-transversal $T_v$ to $vE$. Consider for $e \in T_v$ the virtual homomorphism $\sigma_e \colon \Grpd_e \to \Grpd^{\sour(e)}_{\sour(e)}$, $\sigma_e(g)=g\vert_e$. 
For each $w \in \sour(\bigsqcup_{v\in T^0_{\reg}}T_v)$, choose $h_w \in G$ with $\sour(h_w) = w$ and $\ranT(h_w) = t(w) \in T^0$ and set $c_w \colon \Grpd^w_w \to \Grpd^{t(w)}_{t(w)}$ to be the homomorphism $g \mapsto h_w^{-1} g h_w$. 
Then, 
\[ \Phi_i = \bigoplus_{v\in T^0_{\reg}}\sum_{e \in T_v} K_i(c_{\sour(e)}) \circ K_i(\sigma_e)\circ \mathrm{tr}^{\Grpd}_{\Grpd_e}. \]
\end{theorem}

For group actions, we have this simpler form:

\begin{corollary}[{\cite[Corollary 5.6]{MS}}]\label{c:transitive.transfer.Ktheory}
Let $(G,E,\sigma)$ be a self-similar group action on a finite alphabet $E$ with cocycle $\sigma$.  For $e\in E$, let $\sigma_e\colon G_e\to G$ be the virtual endomorphism $\sigma_e(g)=g\vert_e$.  Then there is a long exact sequence
\[\begin{tikzcd}[arrow style=math font]
	{K_0(C^*(G))} & {K_0(C^*(G))} & {K_0(\OGE)} \\
	{K_1(\OGE)} & {K_1(C^*(G))} & {K_1(C^*(G))}
	\arrow["{1-\Phi_0}", from=1-1, to=1-2]
	\arrow[from=1-2, to=1-3]
	\arrow[from=1-3, to=2-3]
	\arrow[from=2-1, to=1-1]
	\arrow[from=2-2, to=2-1]
	\arrow["{1-\Phi_1}", from=2-3, to=2-2]
\end{tikzcd}\]
where $\Phi_i = \sum_{e \in T} K_i(\sigma_e)\circ \mathrm{tr}^{G}_{G_e}$ for $i = 0,1$ and any $G$-transversal $T \subseteq E$.
\end{corollary}

Corollary \ref{c:transitive.transfer.Ktheory} reduces the computation of $K_i(\OGE)$ to compute $K_i(C^*(G))$ for an arbitrary discrete group $G$. Unfortunately, this is an open problem related to the Baum-Connes Conjecture; this conjecture has been verified for a large number of groups, including all amenable groups.

Also note that \cite{MS} generalizes the result of Nyland and Ortega~\cite{ON} on homology of groupoids associated to Katsura algebras. Their results allow arbitrary cardinality row-finite graphs and sources, while previous results restrict to countable row-finite graphs and no sources.

Finally, for Katsura algebras, we have the following result:

\begin{corollary}[{\cite[Corollary 6.3]{MS}}]\label{c:Katsura}
Let $A,B$ be integer matrices over some index set $J$ with $A$ the adjacency matrix of a row-finite graph such that $A_{ij}=0$ implies $B_{ij}=0$.  Let $J'\subseteq J$ be the set of indices of zero rows. Let $A',B'$ be the matrices obtained from $A,B$, respectively, by removing the rows corresponding to indices in $J'$. Then we have:
\begin{enumerate}
  \item $H_0(\Grpd_{A,B})\cong \coker (\id -(A')^T)$,
  \item $H_1(\Grpd_{A,B})\cong \ker(\id -(A')^T)\oplus \coker (\id -(B')^T)$,
  \item $H_2(\Grpd_{A,B})\cong \ker (\id -(B')^T)$,
\end{enumerate}
and $H_n(\Grpd_{A,B})=0$ for $n\geq 3$.  Moreover, $K_0(\mathcal O_{A,B})\cong \coker (\id -(A')^T)\oplus \ker (\id -(B')^T)$ and $K_1(\mathcal O_{A,B)}) \cong \ker(\id -(A')^T)\oplus \coker (\id -(B')^T)$.
\end{corollary}

\end{document}